\documentclass{article}

\usepackage{microtype}
\usepackage{graphicx}
\usepackage{float}
\usepackage{caption}
\usepackage{subcaption}

\usepackage{booktabs} 

\usepackage{hyperref}

\usepackage[accepted]{icml2024}

\usepackage{amsmath}
\usepackage{amssymb}
\usepackage{mathtools}
\usepackage{amsthm}
\usepackage{physics}
\usepackage{breqn}

\usepackage[capitalize,noabbrev]{cleveref}

\theoremstyle{plain}
\newtheorem{theorem}{Theorem}[section]

\newtheorem{lemma}[theorem]{Lemma}

\theoremstyle{definition}

\theoremstyle{remark}

\newcommand{\paren}[1]{\left(#1\right)}

\usepackage[textsize=tiny]{todonotes}

\icmltitlerunning{UGrid: An Efficient-And-Rigorous Neural Multigrid Solver for Linear PDEs}

\begin{document}

\twocolumn[
\icmltitle{UGrid: An Efficient-And-Rigorous Neural Multigrid Solver for Linear PDEs}



\icmlsetsymbol{equal}{*}

\begin{icmlauthorlist}
\icmlauthor{Xi Han}{sbu}
\icmlauthor{Fei Hou}{iscas,ucas}
\icmlauthor{Hong Qin}{sbu}
\end{icmlauthorlist}

\icmlaffiliation{sbu}{Department of Computer Science, Stony Brook University (SUNY), Stony Brook, NY 11794, USA.}
\icmlaffiliation{iscas}{Key Laboratory of System Software (Chinese Academy of Sciences) and State Key Laboratory of Computer Science, Institute of Software, Chinese Academy of Sciences, Beijing, 100190, China.}
\icmlaffiliation{ucas}{University of Chinese Academy of Sciences, Beijing, 100049, China}

\icmlcorrespondingauthor{Xi Han}{xihan1@cs.stonybrook.edu}
\icmlcorrespondingauthor{Fei Hou}{houfei@ios.ac.cn}
\icmlcorrespondingauthor{Hong Qin}{qin@cs.stonybrook.edu}

\icmlkeywords{Partial Differential Equations, Numerical Solver, Multigrid Method, Machine Learning, ICML}

\vskip 0.3in
]



\printAffiliationsAndNotice{}  

\begin{abstract}
{
Numerical solvers of Partial Differential Equations (PDEs) 
are of fundamental significance to science and engineering. 
To date, the historical reliance on legacy techniques 
has circumscribed possible integration of big data knowledge 
and exhibits sub-optimal efficiency for certain PDE formulations, 
while data-driven neural methods typically lack 
mathematical guarantee of convergence and correctness. 
This paper articulates a mathematically rigorous neural solver 
for linear PDEs. 
The proposed \textbf{\textit{UGrid}} solver,
built upon the principled integration of 
\textbf{\textit{U}}-Net and Multi\textbf{\textit{Grid}},
manifests a mathematically rigorous proof of both convergence 
and correctness, and showcases high numerical accuracy, 
as well as strong generalization power to various 
input geometry/values and multiple PDE formulations. 
In addition, we devise a new residual loss metric, 
which enables self-supervised training and affords more stability 
and a larger solution space over the legacy losses. 
}
\end{abstract}

\section{Introduction}
\label{sec:introduction}

\textbf{Background and Major Challenges}. 
{
PDEs are quintessential to various computational problems 
in science, engineering, and relevant applications
in simulation, modeling, and scientific computing. 
Numerical solutions play an irreplaceable role in common practice 
because in rare cases do PDEs have analytic solutions, 
and many general-purpose numerical methods have been made available. 
Iterative solvers \cite{Saad03} are one of the most-frequently-used
methods to obtain a numerical solution of a PDE. Combining iterative
solvers with the multigrid method \cite{Briggs00} significantly
enhances the performance for large-scale problems. 
Yet, the historical
reliance on legacy generic numerical solvers has circumscribed
possible integration of big data knowledge and exhibits sub-optimal
efficiency for certain PDE formulations. 
In contrast, recent deep neural methods 
have the potential to learn such knowledge from big data
and endow numerical solvers with compact structures and high
efficiency, and have achieved impressive results \cite{Marwah21}. 
However, many currently available neural methods 
treat deep networks as black boxes. 
Other neural methods are typically
trained in a fully supervised manner on loss functions that directly
compare the prediction and the ground truth solution, confining the 
solution space and resulting in numerical oscillations in the relative 
residual error even after convergence. 
These methods generally have challenges unconquered including, 
a lack of sound mathematical backbone, 
no guarantee of correctness or convergence, and low accuracy, 
thus unable to handle complex, unseen scenarios. 
}

\begin{figure}[!htb]
    \centering
    \includegraphics[width=\linewidth]{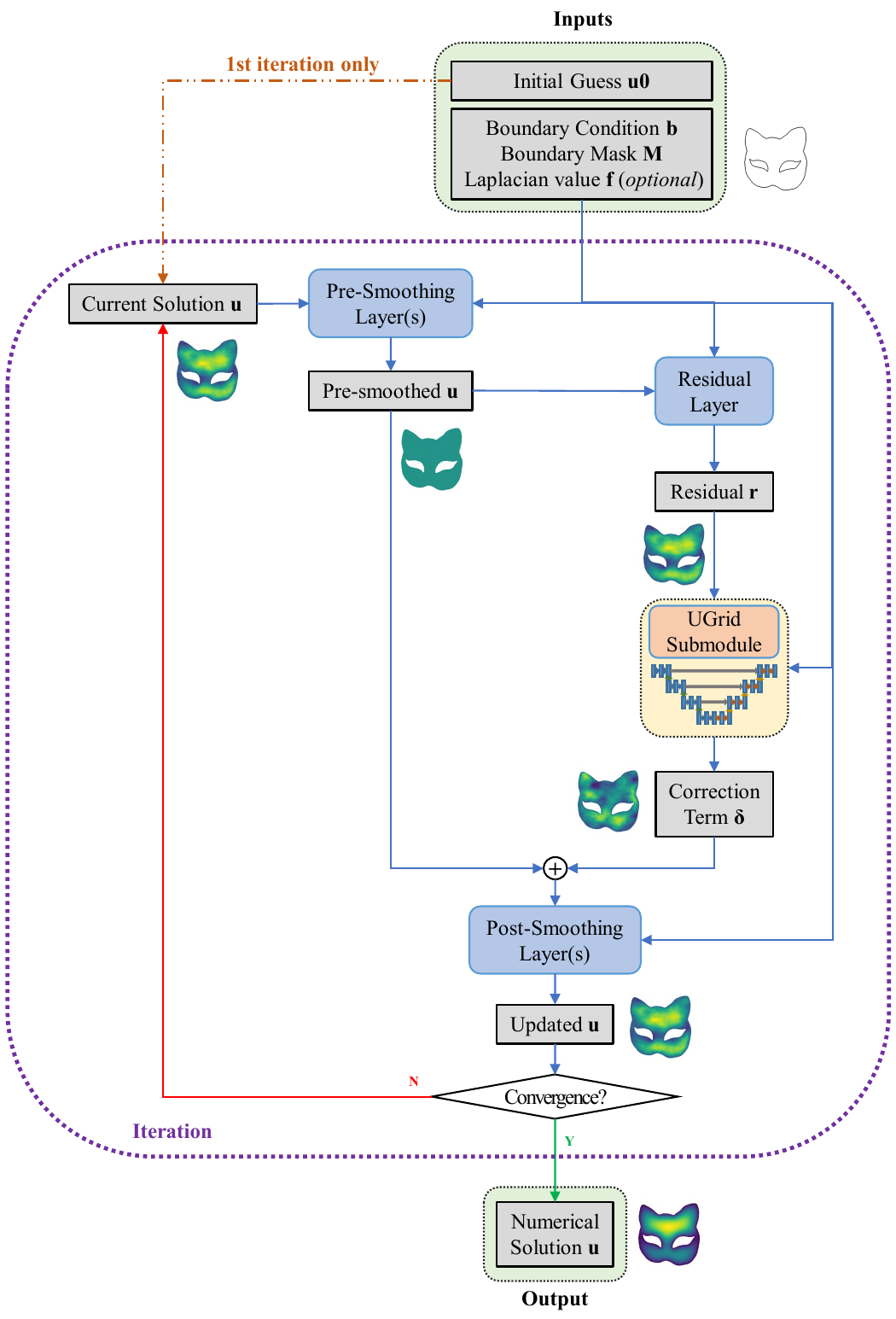}
    \captionsetup{belowskip=0pt}
    \caption{
    Overview of our novel method. 
    Given PDE parameters and its current numerical estimation, 
    the smoothing operations are applied multiple times first. 
    Then, the current residual is fed into our UGrid submodule 
    (together with the boundary mask). 
    Next, the solution is corrected by the correction term
    and post-smoothed. 
    Collectively, it comprises one iteration of 
    the neural solver. 
    The UGrid submodule (detailed in 
    Fig.~\ref{fig:ugrid_submodule})
    aims to mimic the multigrid V-cycle, 
    and its parameters are learnable, 
    so as to endow our framework with 
    the ability to learn from data.}
    \label{fig:pipeline}
\end{figure}

\textbf{Motivation and Method Overview}. 
Inspired by \cite{Hsieh19}'s prior work on 
integrating the structure of multigrid V-cycles \cite{Briggs00}
and U-Net \cite{Ronneberger15} with convergence guarantee, 
and to achieve high efficiency and strong robustness,
we aim to fully realize neural methods' 
modeling and computational potential 
by implanting the legacy numerical methods' mathematical backbone 
into neural methods in this paper. 
In order to make our new framework fully explainable,
we propose the UGrid framework 
(illustrated in Fig.~\ref{fig:pipeline})
based on the structure of multigrid V-cycles 
for learning the functionality of multigrid solvers. 
We also improve the convolutional operators
originating from \cite{Hsieh19} 
to incorporate arbitrary boundary conditions 
and multiple differential stencils without modifying 
the \textit{overall structure} of the key iteration process,
and transform the iterative update rules and 
the multigrid V-cycles into a concrete 
Convolutional Neural Network (CNN) structure. 

\textbf{Key Contributions}. 
The salient contributions of this paper comprise: 
(1) \textit{Theoretical insight}. 
We introduce a novel explainable neural PDE solver 
founded on a solid mathematical background, 
affording high efficiency, high accuracy, 
and strong generalization power to linear PDEs; 
(2) \textit{New loss metric}. 
We propose a residual error metric as the loss function,  
which optimizes the residual of the prediction.  
Our newly-proposed error metric 
enables self-supervised learning and facilitates 
the unrestricted exploration of the solution space. 
Meanwhile, it eliminates the numerical oscillation
on the relative residual error upon convergence, 
which has been frequently observed on the
legacy mean-relative-error-based loss metrics; 
and 
(3) \textit{Extensive experiments}. 
We demonstrate our method's capability to 
numerically solve PDEs by learning multigrid operators 
of various linear PDEs subject to arbitrary 
boundary conditions of complex geometries and topology, 
whose patterns are unseen during the training phase. 
Extensive experiments and comprehensive evaluations 
have verified all of our claimed advantages, 
and confirmed that our proposed method outperforms 
the SOTA. 


%
%
\section{Related Work}
\label{sec:related_work}

\textbf{Black-box-like Neural PDE Solvers}. 
Much research effort has been devoted
to numerically solve PDEs 
with neural networks and deep learning techniques. 
However, most of the previous work treats neural networks
as black boxes and thus come with no mathematical proof 
of convergence and correctness. 
As early as the 1990s, 
\cite{Wang90, WangMendel90, Wang91} 
applied simple neural networks to solve linear equations. 
Later, more effective neural-network-based methods 
like 
\cite{Polycarpou91, Cichocki92, Lagaris98} 
were proposed to solve the Poisson equations. 
On the other hand, 
\cite{Wu94, Xia99, Takala03, Liao10, Li17}
used Recurrent Neural Networks (RNNs) 
in solving systems of linear matrix equations.
Most recently, the potential of CNNs and 
Generative Adversarial Networks (GANs) on 
solving PDEs was further explored by
\cite{Tompson17, Tang17, Farimani17, Sharma18, Ozbay21}. 
Utilities used for neural PDE solvers also include 
backward stochastic differential equations \cite{Han18}
and PN junctions \cite{Zhang19}. 
On the contrary, the proposed UGrid mimics the multigrid solver, and all its contents are explainable and have corresponding counterparts in an MG hierarchy. 

\textbf{Physics-informed Neural PDE Solvers}. 
Physics-informed Neural Network (PINN)-based solvers 
effectively optimize the residual of the solution.
Physical properties, including pressure, velocity \cite{Yang16}
and non-locality \cite{Pang20} are also used to articulate 
neural solvers. 
Mathematical proofs on the minimax optimal bounds \cite{LuChen22}
and structural improvements \cite{Lu21, Lu21DeepXDE}
are also made on the PINN architecture itself, 
endowing physics-informed neural PDE solvers 
with higher efficiency and interpretability. 
Hinted by these, we propose the residual loss metric, 
which enables self-supervised training, 
enlarges the solution space 
and enhances numerical stability. 

\textbf{Neural PDE Solvers with Mathematical Backbones}. 
\cite{Zhou09}
proposed a NN-based linear system and
its solving algorithm with a convergence guarantee. 
\cite{Hsieh19} 
modified the Jacobi iterative solver by 
predicting an additional correction term 
with a multigrid-inspired linear operator, 
and proposed a linear neural solver with
guarantee on correctness upon convergence. 
\cite{Greenfeld19} 
proposed to learn a mapping from 
a family of PDEs to the optimal prolongation operator 
used in the multigrid method, 
which is then extended to 
Algebraic Multigrids (AMGs) on non-square meshes
via Graph Neural Networks (GNNs) by \cite{Luz20}. 
On the other hand, 
\cite{Li21} 
proposed a Fourier neural operator 
that learns mappings between function spaces 
by parameterizing the integral kernel directly in Fourier space. 
In theory, 
\cite{Marwah21} proved that 
when a PDE's coefficients are representable 
by small NNs, 
the number of parameters needed will increase 
in a polynomial fashion with the input dimension. 

%
%
\vspace{-3mm}
\section{Mathematical Preliminary}
\label{sec:mathematical_preliminary}
\vspace{-1mm}



For mathematical completeness, we provide readers 
with a brief introduction to the concepts that are
frequently seen in this paper. 

\textbf{Discretization of $2$D Linear PDEs}. 
A linear PDE with Dirichlet boundary condition 
could be discretized with finite differencing 
techniques \cite{Saad03}, 
and could be expressed in the following form: 
\begin{equation}
    \left\{
    \begin{aligned}
        \mathcal{D} u(x, y) 
                & = f(x, y), 
                & (x, y) \in \mathcal{I} \\ 
        u(x, y) 
                & = b(x, y), 
                & (x, y) \in \mathcal{B}
    \end{aligned}
    \right. ,
    \label{equ:linear_pde}
\end{equation}
where 
$\mathcal{D}$ 
is a $2$D discrete linear differential operator, 
$\mathcal{S}$ 
is the set of all points on the discrete grid, 
$\mathcal{B}$ 
is the set of boundary points in the PDE, 
$\mathcal{I} = \mathcal{S} \setminus{\mathcal{B}}$ 
is the set of interior points in the PDE, 
$\partial \mathcal{S}$ $\subseteq \mathcal{B}$ 
is the set of \textit{trivial boundary points} of the grid. 

Using $\mathcal{D}$'s corresponding
finite difference stencil, 
Eq.~\ref{equ:linear_pde}
can be formulated into a sparse linear system 
of size $n^2 \times n^2$: 
\begin{equation}
    \left\{
    \begin{aligned}
        (\mathbf{I} - \mathbf{M}) \mathbf{A} \mathbf{u} & = (\mathbf{I} - \mathbf{M}) \mathbf{f} \\ 
        \mathbf{M} \mathbf{u} & = \mathbf{M} \mathbf{b}
    \end{aligned}
    \right. ,
    \label{equ:masked_linear_pde}
\end{equation}
where 
$\mathbf{A} \in \mathbb{R}^{n^2 \times n^2}$ 
is the $2$D discrete differential operator, 
$\mathbf{u} \in \mathbb{R}^{n^2}$ 
encodes the function values of 
the interior points and the non-trivial boundary points; 
$\mathbf{f} \in \mathbb{R}^{n^2}$ 
encodes the corresponding partial derivatives 
of the interior points; 
$\mathbf{b} \in \mathbb{R}^{n^2}$ 
encodes the non-trivial boundary values; 
$\mathbf{I}$ denotes the $n^2 \times n^2$ identity matrix; 
$\mathbf{M} \in \{0, 1\}^{n^2 \times n^2}$
is a diagonal binary boundary mask defined as
\begin{equation}
    \mathbf{M}_{k, k} = 
    \left\{
    \begin{aligned}
        1, \quad & (i, j) \in \mathcal{B} \setminus \partial \mathcal{S} \\
        0, \quad & (i, j) \in \mathcal{I}
    \end{aligned}
    \right. , 
    \,
    k = in + j, \,
    0 \le i, \ j < n .
\end{equation}
On the contrary of Eq.~\ref{equ:linear_pde}, 
both equations in Eq.~\ref{equ:masked_linear_pde} 
hold for \textit{all} grid points. 

\textbf{Error Metric And Ground-truth Solution}. 
When using numerical solvers, 
researchers typically substitute the boundary mask $\mathbf{M}$ 
into the discrete differential matrix $\mathbf{A}$
and the partial derivative vector $\mathbf{f}$, 
and re-formulate Eq.~\ref{equ:masked_linear_pde}
into the following generic sparse linear system: 
\begin{equation}
     \widetilde{\mathbf{A}} \,  \mathbf{u} 
   = \widetilde{\mathbf{f}} .
   \label{equ:general_sparse_linear_system}
\end{equation}
The \textit{residual} of a numerical solution $\mathbf{u}$ 
is defined as
\begin{equation}
    \mathbf{r}(\mathbf{u}) = 
    \widetilde{\mathbf{f}} - 
    \widetilde{\mathbf{A}} \, 
    \mathbf{u} .
    \label{equ:residual}
\end{equation}
In the ideal case, the \textit{absolute residual error} 
of an exact solution $\mathbf{u}^*$ should be 
$r_{\mathbf{u}^*} = \norm{ \mathbf{r}(\mathbf{u}^*) } = 0$. 
However, in practice, a numerical solution $\mathbf{u}$ 
could only be an approximation 
of the exact solution $\mathbf{u}^*$. 
The precision of $\mathbf{u}$ 
is evaluated by its \textit{relative residual error}, 
which is defined as
\begin{equation}
    \varepsilon_{\mathbf{u}} = 
    \norm{\widetilde{\mathbf{f}} - 
          \widetilde{\mathbf{A}} \, \mathbf{u}}
    \bigg/
    \norm{\widetilde{\mathbf{f}}} 
    .
    \label{equ:relative_residual_error}
\end{equation}
Typically, the ultimate goal of a numerical PDE solver
is to seek the optimization of the relative residual error. 
If we have
$\varepsilon_{\mathbf{u}} \le \varepsilon_{\mathrm{max}}$ 
for some small $\varepsilon_{\mathrm{max}}$, 
we would consider $\mathbf{u}$ to be a 
\textit{ground-truth solution}. 

\textbf{Linear Iterator}. 
A linear iterator (also called an iterative solver
or a smoother) 
for generic linear systems like
Eq.~\ref{equ:general_sparse_linear_system}
could be expressed as
\begin{equation}
      \mathbf{u}_{k + 1}
    = \paren{ \mathbf{I} - 
              \widetilde{\mathbf{P}}^{-1} \widetilde{\mathbf{A}} } 
      \mathbf{u}_{k} +   
      \widetilde{\mathbf{P}}^{-1} \widetilde{\mathbf{f}} ,
    \label{equ:regular_iterator}
\end{equation}
where $\widetilde{\mathbf{P}}$ 
is an easily invertible approximation 
to the system matrix $\widetilde{\mathbf{A}}$.

%
%
\section{Novel Approach}
\label{sec:novel_approach}

{
The proposed UGrid neural solver is built upon
the principled integration of 
the U-Net architecture and the multigrid method.
We observe that linear differential operators, 
as well as their approximate inverse in 
legacy iterative solvers, are analogous
to convolution operators. E.g., the discrete  
Laplacian operator is a $3 \times 3$ convolution kernel. 
Furthermore, the multigrid V-cycle is
also analogous to the U-Net architecture, 
with grid transfer operators
mapped to up/downsampling layers. 
Moreover, the fine-tuning process of multigrid's 
critical components on specific PDE 
formulations could be completed 
by learning from big data. 
These technical observations lead to
our neural implementation and optimization
of the multigrid routine. 
}

In spite of high efficiency, 
generalization power remains a major challenge 
for neural methods. 
Many SOTA neural solvers, e.g., \cite{Hsieh19}, 
fail to generalize to new scenarios unobserved 
during the training phase. 
Such new scenarios include: 
(1) New problem sizes; and
(2) New, complex boundary conditions
and right-hand sides, 
which includes geometries, topology,
and values (noisy inputs). 
Moreover, some of these methods are tested on 
Poisson equations only; 
neither mathematical reasoning nor empirical results show that they could trivially generalize to other PDEs
(with or without retraining). 
{
UGrid resolves
all problems above. 

UGrid is comprised of the following components: 
(1) 
The \textit{fixed} neural smoother,
which consists of our proposed convolutional operators
(Sec.~\ref{sec:novel_approach:conv_operators}); 
(2)
The \textit{learnable} neural multigrid,
which consists of our UGrid submodule
(Sec.~\ref{sec:novel_approach:ugrid_submodule});
(3)
A residual loss metric
(Sec.~\ref{sec:novel_approach:residual_loss})
which enables
the self-supervised training process. 
}

{
\subsection{Convolutional Operators}
\label{sec:novel_approach:conv_operators}

This subsection is organized as follows:
Sec.~\ref{sec:novel_approach:conv_operators:mask_conv_iter} 
introduces the masked operators, 
which mimic the smoothers in a legacy multigrid routine; 
and Sec.~\ref{sec:novel_approach:conv_operators:residual_operator}
introduces the masked residual operators 
for residual calculation. 

\subsubsection{Masked Convolutional Iterator}
\label{sec:novel_approach:conv_operators:mask_conv_iter}

A trivial linear iterator in the form of 
Eq.~\ref{equ:regular_iterator} does \textbf{not} 
fit in a neural routine. 
This is because in practice, 
the system matrix $\vb{A}$ for its
corresponding differential operator
encodes the \textit{boundary geometry} and 
turns into matrix $\widetilde{\mathbf{A}}$ in 
Eq.~\ref{equ:general_sparse_linear_system}. 
$\widetilde{\mathbf{A}}$'s elements are input-dependent, 
and is thus impossible to be expressed 
as a \textit{fix-valued} convolution kernel. 

We make use of the masked version of
PDEs (Eq.~\ref{equ:masked_linear_pde})
and their
\textit{masked convolutional iterators}, 
which are natural extensions of 
Eq.~\ref{equ:regular_iterator}: 
\begin{equation}
      \mathbf{u}_{k + 1}
    = \paren{ \mathbf{I} - \mathbf{M} }  
      \paren{ \paren{ \mathbf{I} - \mathbf{P}^{-1} \mathbf{A} } 
              \mathbf{u}_{k} +   
              \mathbf{P}^{-1} \mathbf{f} } + 
      \mathbf{M} \mathbf{b} ,
    \label{equ:masked_iterator}
\end{equation}
where $\mathbf{P}$ is an easily-invertible approximation
on the discrete differential operator $\mathbf{A}$. 
The correctness of Eq.~\ref{equ:masked_iterator} is guaranteed
by the following Theorem 
(proved as Theorem~\ref{thm:mask_iterator_convergence_sm} 
in the appendix): 
\begin{theorem}
    For a PDE in the form of 
    Eq.~\ref{equ:masked_linear_pde}, 
    the masked iterator Eq.~\ref{equ:masked_iterator} 
    converges to its ground-truth solution
    when its prototype Jacobi iterator converges
    and $\mathbf{P}$ is full-rank diagonal. 
    \label{thm:mask_iterator_convergence}
\end{theorem}

Matrix $\vb{A}$ has different formulations
for different linear PDEs. Without loss of generality, 
we consider the following problem formulations. 
Other possible PDE formulations are essentially
combinations of differential primitives
available in the three problems below, 
and our convolutional operators could be
trivially extended to higher orders of 
differentiation. 

\textbf{$2$D Poisson Problem}. 
Under Dirichlet boundary condition,
a Poisson problem could be expressed as follows: 
\begin{equation}
    \left\{
    \begin{aligned}
        \nabla^2 u(x, y) & = f(x, y), & (x, y) \in \mathcal{I} \\
        u(x, y) & = b(x, y), & (x, y) \in \mathcal{B}
    \end{aligned}
    \right.
    \ , 
    \label{equ:poisson}
\end{equation}
where $u$ is the unknown scalar field, 
$f$ is the Laplacian field, 
and $b$ is the Dirichlet boundary condition. 

Matrix $\vb{A}$ could be assembled by the five-point finite
difference stencil 
for $2$D Laplace operators \cite{Saad03}, 
and we could simply let $\mathbf{P} = -4 \mathbf{I}$, 
where $\mathbf{I}$ denotes the identity matrix. 
The update rule specified in Eq.~\ref{equ:masked_iterator}
thus becomes 
\begin{dmath}
         \vb{u}_{k + 1}(i, j)
    = \dfrac{1}{4} \qty( \vb{I} - \vb{M} )
      ( \vb{u}_{k}(i - 1, j) + 
            \vb{u}_{k}(i + 1, j) + 
            \vb{u}_{k}(i, j - 1) +
            \vb{u}_{k}(i, j + 1) - 
            \vb{f} ) + 
      \vb{M} \vb{b} . 
\end{dmath}
To transform the masked iterator into a convolution layer, 
we reorganize the column vectors
$\mathbf{u}$, $\mathbf{b}$, $\mathbf{M}$ and $\mathbf{f}$ 
into $n \times n$ matrices 
with their semantic meanings untouched. 
Then, the neural \textit{smoother} could be expressed as
\begin{dmath}
    \begin{aligned}
            \vb{u}_{k + 1} 
        & = \mathrm{smooth}(\vb{u}_k) \\
        & = \qty( \vb{1} - \vb{M} )
            \qty( \vb{u}_{k} * \vb{J} - 0.25 \vb{f} ) + 
            \vb{M} \vb{b} \ , \\
            \mathbf{J} 
        & = \begin{pmatrix}
            0    & 0.25 & 0 \\ 
            0.25 & 0    & 0.25 \\
            0    & 0.25 & 0
            \end{pmatrix} \ , 
    \end{aligned}
    \label{equ:masked_iterator_func}
\end{dmath}
where $\vb{1}$ is an $n \times n$ matrix whose elements
are all equal to $1$, 
and $*$ denotes $2$D convolution. 

\textbf{$2$D Helmholtz Problem}. 
Under Dirichlet boundary condition,
an inhomogeneous Helmholtz equation with 
spatially-varying wavenumber may be expressed as follows: 
\begin{equation}
    \left\{
    \begin{aligned}
        \nabla^2 u(x, y) + k^2(x, y) u(x, y) & = f(x, y), & (x, y) \in \mathcal{I} \\
        u(x, y) & = b(x, y), & (x, y) \in \mathcal{B}
    \end{aligned}
    \right.
    \ , 
    \label{equ:helmholtz}
\end{equation}
where $u$ is the unknown scalar field, 
$k^2$ is the spatially-varying wavenumber, 
$f$ is the (non-zero) right hand side, 
and $b$ is the Dirichlet boundary condition. 

For our proposed UGrid solver, we could naturally extend 
Eq.~\ref{equ:masked_iterator_func} into the following form
to incorporate Eq.~\ref{equ:helmholtz}: 
\begin{dmath}
    \vb{u}_{k + 1} = 
    \dfrac{1}{\vb{4} - \vb{k^2}}
    \qty(\vb{1} - \vb{M}) \qty(\vb{u}_k * 4\vb{J} - \vb{f}) + 
    \vb{M} \vb{b}
    \ ,
\end{dmath}
where all notations retain their 
meanings as in Eq.~\ref{equ:masked_iterator_func}. 

\textbf{$2$D Steady-state Convection-diffusion-reaction Problem}. 
Under Dirichlet boundary condition, an inhomogeneous 
steady-state convection-diffusion-reaction equation
may be expressed as follows: 
\begin{equation}
    \left\{
    \begin{aligned}
        \vb{v}(x, y) \cdot \grad{u(x, y)} - 
        \alpha \nabla^2 u(x, y) + 
        \beta u(x, y) & = f(x, y), \\
        & (x, y) \in \mathcal{I} \\
        u(x, y) & = b(x, y), \\
        & (x, y) \in \mathcal{B}
    \end{aligned}
    \right.
    \ , 
    \label{equ:diffusion}
\end{equation}
where $u$ is the unknown scalar field, 
$\vb{v} = (v_x, v_y)^{\top}$ is the vector velocity field, 
$\alpha$, $\beta$ are constants, 
$f$ is the (non-zero) right-hand side, 
and $b$ is the Dirichlet boundary condition. 

For our proposed UGrid solver, 
we could naturally extend 
Eq.~\ref{equ:masked_iterator_func} into the following form
to incorporate Eq.~\ref{equ:diffusion}: 
\begin{dmath}
    \vb{u}_{k + 1} = 
    \dfrac{1}{4\alpha + \beta}
    \qty(\vb{1} - \vb{M}) \\
    \qty(
        \alpha \vb{u}_k * 4 \vb{J} + 
        \vb{v_x} \qty(\vb{u}_k * \vb{J_x}) + 
        \vb{v_y} \qty(\vb{u}_k * \vb{J_y}) + 
        \vb{f}
    ) \\
    + 
    \vb{M} \vb{b}
    \ ,
\end{dmath}
\begin{align*}
    \vb{J_x} =     
    \begin{pmatrix}
      0 & 0 &    0 \\ 
    0.5 & 0 & -0.5 \\
      0 & 0 &    0
    \end{pmatrix} \ , 
    \quad
    \vb{J_y} =     
    \begin{pmatrix}
    0 & -0.5 & 0 \\ 
    0 &    0 & 0 \\
    0 &  0.5 & 0
    \end{pmatrix} \ , 
\end{align*}
where $\vb{J_x}$ and $\vb{J_y}$
are two convolution kernels introduced 
for the gradient operator in Eq.~\ref{equ:diffusion}, 
and all notations retain their 
meanings as in Eq.~\ref{equ:masked_iterator_func}.

\subsubsection{Convolutional Residual Operator}
\label{sec:novel_approach:conv_operators:residual_operator}

Except for the smoother, the multigrid method 
also requires the calculation of the 
residual in each iteration step. 
In practice, the \textit{residual operator} 
(Eq.~\ref{equ:residual}) 
can also be seamlessly implemented as a convolution layer. 
Because our masked iterator (Eq.~\ref{equ:masked_iterator})
guarantees that $\mathbf{u}$ satisfies 
$\mathbf{M} \mathbf{u} = \mathbf{M} \mathbf{b}$
at any iteration step, 
the residual operator 
for Poisson problems could be simplified into 
\begin{equation}
      \mathbf{r}(\mathbf{u})
    = \paren{ \mathbf{1} - \mathbf{M} } 
      \paren{ \mathbf{f} - \mathbf{u} * \mathbf{L} } 
    \ , \quad
    \mathbf{L} = 
    \begin{pmatrix}
    0 & 1 & 0 \\
    1 & -4 & 1 \\ 
    0 & 1 & 0
    \end{pmatrix} \ .
    \label{equ:residual_func}
\end{equation}

For Helmholtz problems, Eq.~\ref{equ:residual_func} 
could be naturally extended into
\begin{equation}
    \vb{r}(\vb{u}) = \qty(\vb{1} - \vb{M}) \qty(\vb{f} - \vb{u} * \vb{L} - \vb{k^2} \vb{u})
    \ ,
\end{equation}
where all notations retain their
meanings as in Eq.~\ref{equ:residual_func}. 

For steady-state convection-diffusion-reaction problems, 
we could extend Eq.~\ref{equ:residual_func} into
\begin{dmath}
    \vb{r}(\vb{u}) = 
    \qty(\vb{1} - \vb{M}) \\
    \qty(
        \vb{f} + 
        \vb{v_x} \qty(\vb{u} * \vb{J_x}) + 
        \vb{v_y} \qty(\vb{u} * \vb{J_y}) + 
        \alpha \vb{u} * \vb{L} - 
        \beta \vb{u}
    ) \ ,
\end{dmath}
where all notations retain their
meanings as in Eq.~\ref{equ:residual_func}. 
}

\subsection{Neural Network Design}
\label{sec:novel_approach:ugrid_submodule}

\textbf{UGrid Iteration}. 
We design the iteration step of our neural iterator 
as a sequence of operations as follows 
(which is illustrated in Fig.~\ref{fig:pipeline}): 
\begin{equation}
    \begin{aligned}
        \mathbf{u} & = \mathrm{smooth}^{\nu_1}(\mathbf{u}) 
        & \text{(Pre-smooth for $\nu_1$ times)}; \\ 
        \mathbf{r} & = \mathbf{r}(\mathbf{u}) 
        & \text{(Calculate the current residual)}; \\
        \delta & = 
        \mathrm{UGrid}(\mathbf{r}, \mathbf{1} - \mathbf{M}) 
        & \text{(UGrid submodule recursion)}; \\ 
        \mathbf{u} & = \mathbf{u} + \delta
        & \text{(Apply the correction term)}; \\ 
        \mathbf{u} & = \mathrm{smooth}^{\nu_2}(\mathbf{u})
        & \text{(Post-smooth for $\nu_2$ times)} .
    \end{aligned}
    \label{equ:ugrid_iterator}
\end{equation}
The entire iteration process is specifically designed to emulate 
the multigrid iteration \cite{Saad03}: 
We use the pre-smoothing and post-smoothing layers 
(as specified in Eq.~\ref{equ:masked_iterator_func}) 
to eliminate the high-frequency modes in the residual $\mathbf{r}$, 
and invoke the \textit{UGrid submodule} 
to eliminate the low-frequency modes. 

\begin{figure}[!htb]
    \centering
    \includegraphics[width=0.8\linewidth]{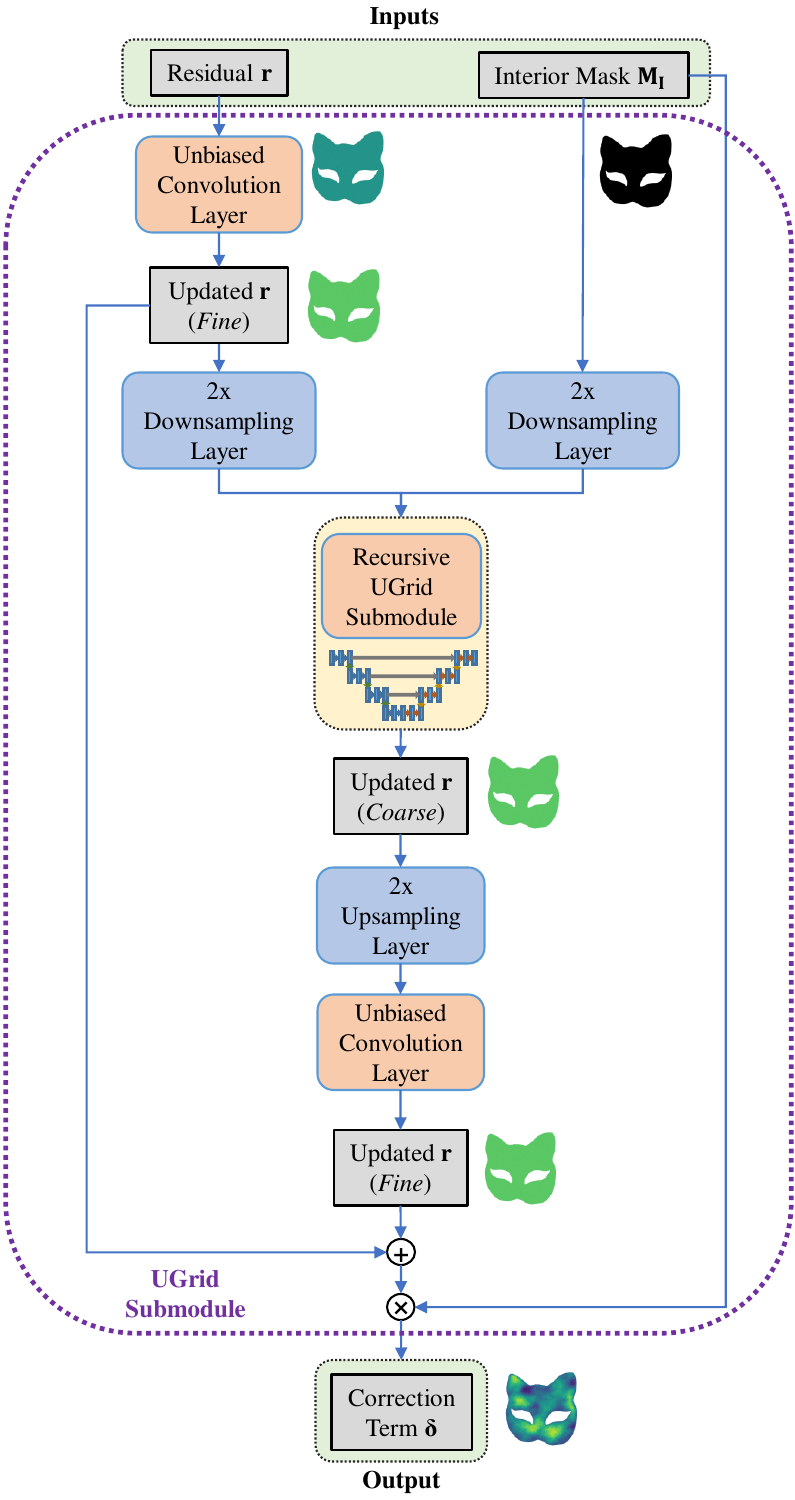}
    \caption{
    An overview of our recursive UGrid submodule. 
    The residual is smoothed by 
    \textit{unbiased} convolution layers,
    downsampled to be \textit{recursively updated} 
    by a $2$x-coarser UGrid submodule
    (note the orange recursive invocation
    of UGrid submodule in the middle), 
    then upsampled back to the fine grid, 
    smoothed, and added with the initial residual
    by skip-connection. 
    Boundary values are enforced by
    interior mask 
    via element-wise multiplication.
    The convolution layers (shown in orange) 
    are learnable; other layers (shown in blue) 
    are the fixed mathematical backbone.
    }
    \label{fig:ugrid_submodule}
\end{figure}

\textbf{UGrid Submodule}. 
Our UGrid submodule is also implemented as 
a fully-convolutional network, 
whose structure is highlighted in Fig.~\ref{fig:ugrid_submodule}. 
The overall structure of UGrid is built upon 
the principled combination of 
U-Net \cite{Ronneberger15} and multigrid V-cycle, 
and could be considered a ``V-cycle'' with skip connections. 
Just like the multigrid method, 
our UGrid submodule is also invoked recursively,
where each level of recursion would coarsen the mesh grid by $2$x. 

To approximate the linearity of the multigrid iteration
(note that we are learning to invert a system matrix), 
we implement the smoothing layers
in the legacy multigrid V-cycle 
(not to be confused with 
the pre-smoother and the post-smoother 
in Eq.~\ref{equ:ugrid_iterator}, 
which are \textit{outside} of the V-cycle hierarchy)
as learnable $2$D convolution layers \textit{without} any bias. 
For the same reason, and inspired by \cite{Hsieh19}, 
we also drop many commonly seen neural layers 
which would introduce non-linearity, 
such as normalization layers and activation layers. 

\subsection{Loss Function Design}
\label{sec:novel_approach:residual_loss}

\textbf{Legacy Loss Metric}. 
We refer the equivalents of the mean relative error 
between the prediction and a ground-truth solution  
as the \textit{legacy loss metric}. 
The following Theorem 
(proved in Sec.~\ref{sec:supplemental_material:legacy_loss})
shows that though intuitive, 
the legacy loss is unstable: 
\begin{theorem}
    When a neural network converges on a legacy loss metric
    such that its prediction $\mathbf{x}$ satisfies 
    $
        \mathcal{L}_{\mathrm{legacy}}(\mathbf{x}, \mathbf{y}) = 
        \mathrm{mean}
        \paren{ \abs{ \mathbf{x} - \mathbf{y} } /
                \abs{ \mathbf{y} }  }
        \le l_{\mathrm{max}}  
    $, 
    where $\mathbf{y}$ denotes the ground truth value, 
    the upper and lower bounds of
    $\mathbf{x}$'s relative residual error 
    are still dependent on the input. 
    \label{thm:legacy_loss}
\end{theorem}

Theorem~\ref{thm:legacy_loss} explains
our experimental observations that: 
(1) Optimizing the legacy loss metric 
does \textbf{not} increase the precision
in terms of the relative residual error; and
(2) The legacy loss metric restricts 
the solution space: A valid numerical solution with
low relative residual error may have a 
large relative difference from another valid solution 
(the one selected as the ground truth value). 
As a result, many valid solutions are 
\textit{unnecessarily rejected} 
by the legacy loss metric.

\textbf{Proposed Residual Loss Metric}. 
To overcome the shortcomings of the 
legacy loss metric, 
we propose to optimize the neural solver
directly to the residual in a self-supervised manner: 
\begin{equation}
    \mathcal{L}_{\vb{r}_{\mathrm{abs}}}(\vb{x}) = 
    \mathbb{E}_{\vb{x}}
    \qty[
    \norm{ \qty( \vb{1} - \vb{M} ) 
           \qty(\vb{f} - \vb{A} \, \vb{x} ) 
    }_2 
    ]
    \label{equ:loss}
    \ .
\end{equation}

{
We have conducted ablation studies on
our relative loss metric 
and the legacy loss (Sec.~\ref{sec:experiments} and Sec.~\ref{sec:appendix:ablation_study}). 
The results showcase that the proposed 
residual loss
performs better than the legacy loss
in terms of both efficiency and accuracy. 
Heuristically, the proposed residual loss metric is closer to
the fundamental definition of the precision 
of a PDE's solution, and is more robust and stable, 
because upon convergence, 
it guarantees an input-independent final accuracy. 

Moreover, the self-supervised training process 
endows our method with the following merits:
(1) Easier data generation (compared to other neural
routines which are trained on the legacy loss), 
and thus achieve better numerical performance;
(2) For a specific PDE formulation, we could 
easily get a decent neural multigrid solver
optimized for that specific formulation
(which outperforms existing general-purpose 
legacy routines), 
simply by training UGrid on the data generated. 
On the contrary, fine-tuning legacy solvers 
is a non-trivial task requiring a solid 
mathematical background 
as well as non-trivial programming effort. 
}

%
%
{

\section{Experiments and Evaluations}
\label{sec:experiments}

\textbf{Experiments Overview}. 
For each of the three 
PDEs mentioned in 
Sec.~\ref{sec:novel_approach:conv_operators}: 
(i) Poisson problem, 
(ii) inhomogeneous Helmholtz problem
with varying wave numbers, 
and (iii) inhomogeneous steady-state 
diffusion-convection-reaction problem, 
we train one UGrid model specialized 
for its formulation. 
We apply our model and the baselines 
to the task of $2.5$D freeform surface modeling. 
These surfaces are modeled by 
the three types of PDEs as $2$D height fields, 
with non-trivial geometry/topology. 
Each surface is discretized into: 
(1) Small-scale problem: 
A linear system of size $66,049 \times 66,049$; 
(2) Large-scale problem: 
A linear system of size $1,050,625 \times 1,050,625$;
(3) XL-scale problem:
A linear system of size $4,198,401 \times 4,198,401$; 
and (4) XXL-scale problem:
A linear system of size $16,785,409 \times 16,785,409$. 
UGrid is trained \textbf{on the large scale only}. 
Other problem sizes are designed to evaluate
UGrid's generalization power and scalability.

In addition, we have conducted an ablation
study on the residual loss metric
(v.s. legacy loss)
as well as the UGrid architecture itself
(v.s. vanilla U-Net). 

\textbf{Data Generation and Implementation Details}. 
Our new neural PDE solver is trained 
in a self-supervised manner on the residual loss. 
Before training, we synthesized a dataset with 
$16000$ $\qty(\mathbf{M}, \mathbf{b}, \mathbf{f})$ pairs. 
For Helmholtz and diffusion problems, we further
random-sample their unique coefficient fields
(more details available in Sec.~\ref{sec:supplemental_material:eval_helmholtz}
and Sec.~\ref{sec:supplemental_material:eval_diffusion}.)
To examine the generalization power of UGrid,
the geometries of boundary conditions in our training data 
are limited to ``Donuts-like'' shapes as shown in 
Fig.~\ref{fig:testcase} (h).
Moreover, all training data are restricted to
zero $\vb{f}$-fields \textit{only}, 
i.e., $\mathbf{f} \equiv \mathbf{0}$. 
Our UGrid model has $6$ recursive Multigrid submodules. 
We train our model and perform all experiments 
on a personal computer with $64$GB RAM, 
AMD Ryzen $9$ $3950$x $16$-core processor, 
and NVIDIA GeForce RTX $2080$ Ti GPU. 
We train our model for $300$ epochs with the Adam optimizer. 
The learning rate is initially set to $0.001$, 
and decays by $0.1$ for every $50$ epochs. 
We initialize all learnable parameters 
with PyTorch's default initialization policy. 
Our code is available as open-source at 
\url{https://github.com/AXIHIXA/UGrid}.

\begin{figure*}[!htb]
    \centering
    \begin{subfigure}[H]{0.19\linewidth}
        \centering
        \includegraphics[height=\linewidth]{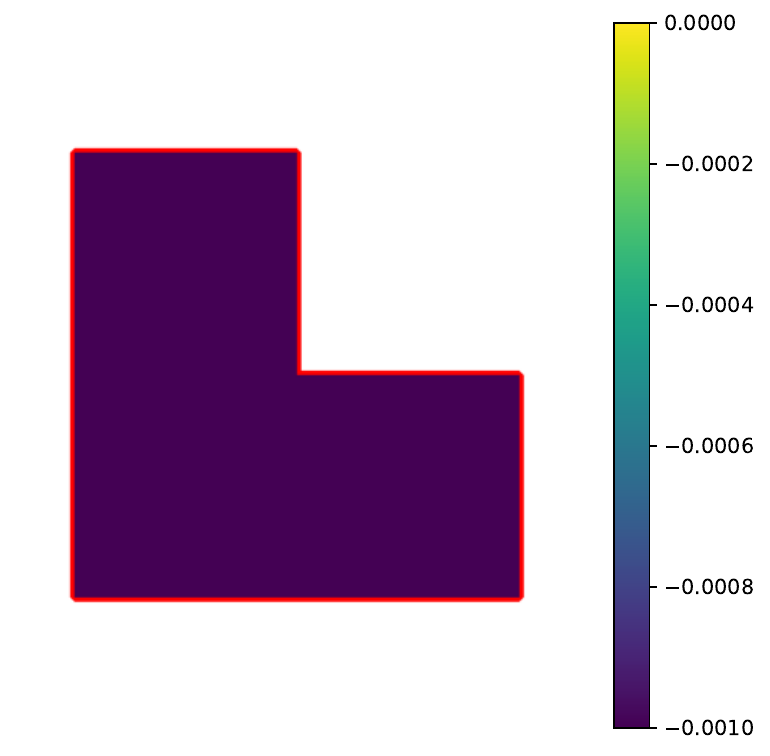}
        \caption{L-Shape}
    \end{subfigure}
    \begin{subfigure}[H]{0.19\linewidth}
        \centering
        \includegraphics[height=\linewidth]{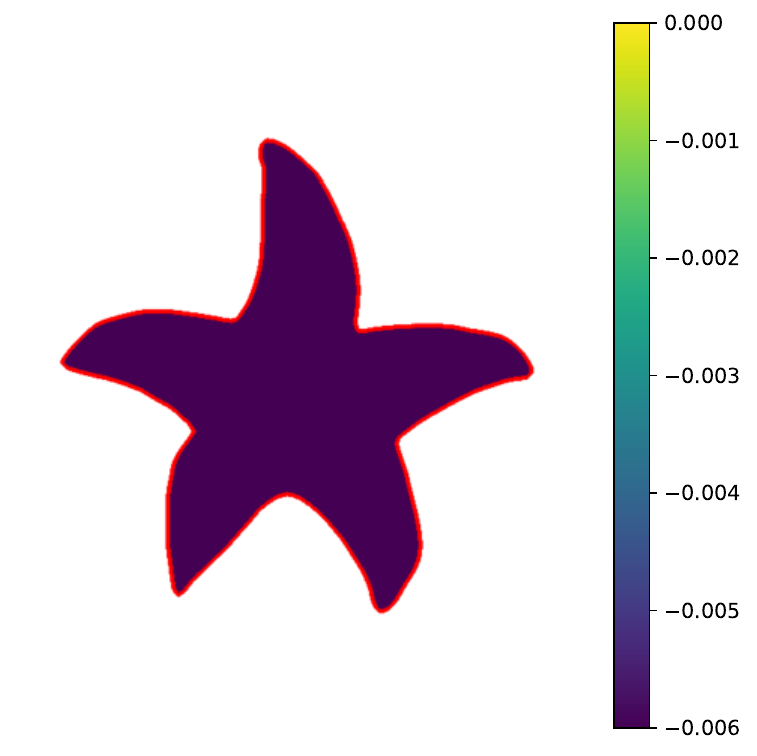}
        \caption{Star}
    \end{subfigure}
    \begin{subfigure}[H]{0.19\linewidth}
        \centering
        \includegraphics[height=\linewidth]{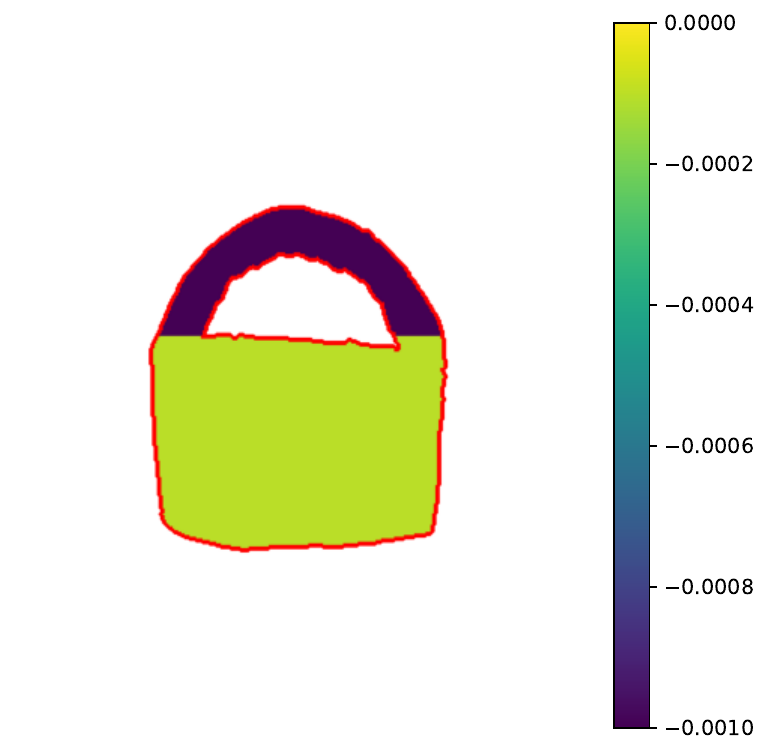}
        \caption{Lock}
    \end{subfigure}
    \begin{subfigure}[H]{0.19\linewidth}
        \centering
        \includegraphics[height=\linewidth]{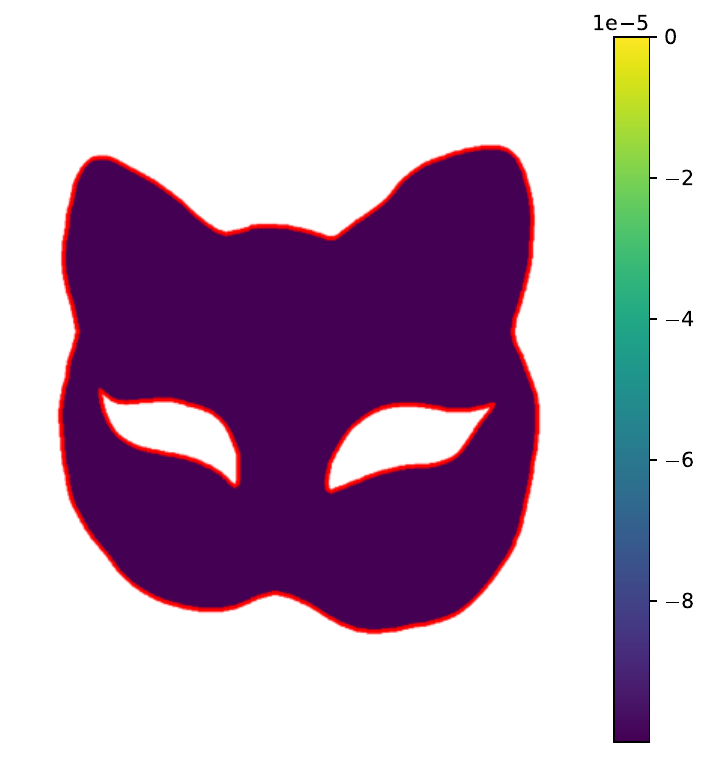}
        \caption{Cat}
    \end{subfigure}
    \begin{subfigure}[H]{0.19\linewidth}
        \centering
        \includegraphics[height=\linewidth]{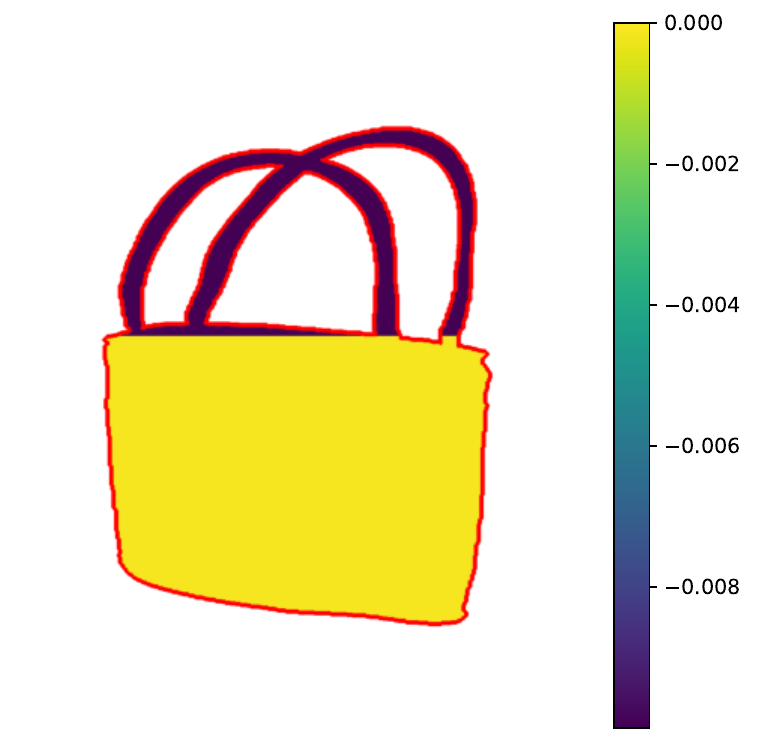}
        \caption{Bag}
    \end{subfigure}
    \begin{subfigure}[H]{0.19\linewidth}
        \centering
        \includegraphics[height=\linewidth]{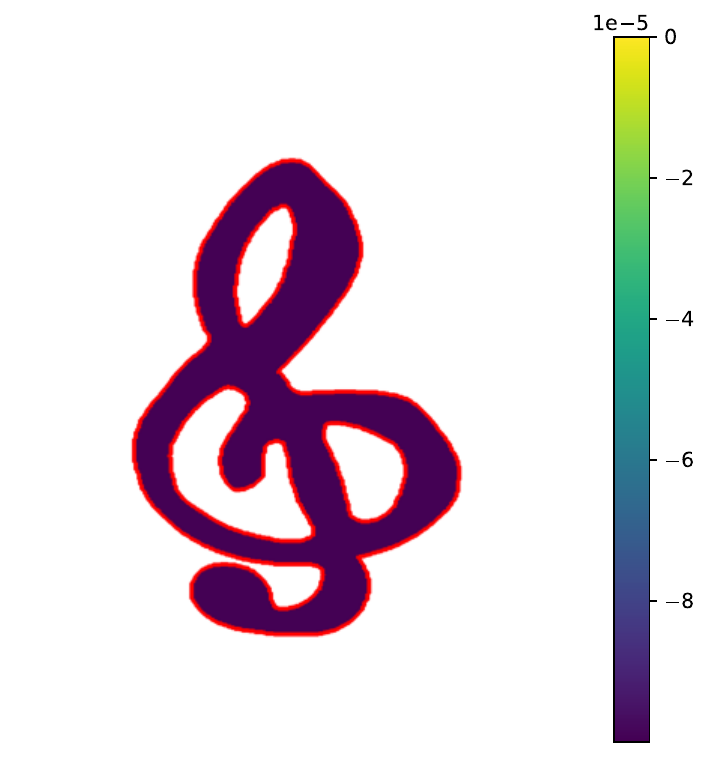}
        \caption{Note}
    \end{subfigure}
    \begin{subfigure}[H]{0.19\linewidth}
        \centering
        \includegraphics[height=\linewidth]{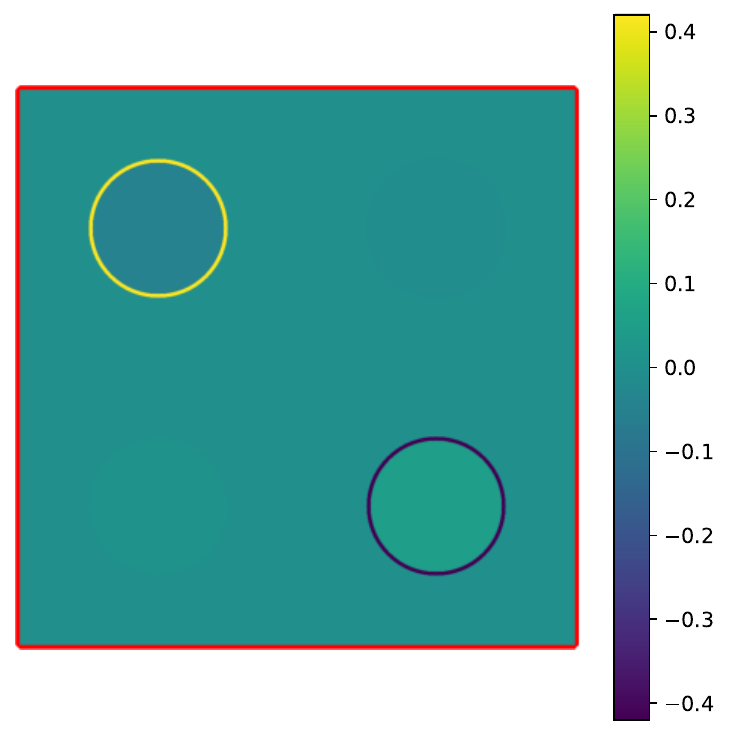}
        \caption{Sharp Feature}
    \end{subfigure}
    \begin{subfigure}[H]{0.19\linewidth}
        \centering
        \includegraphics[height=\linewidth]{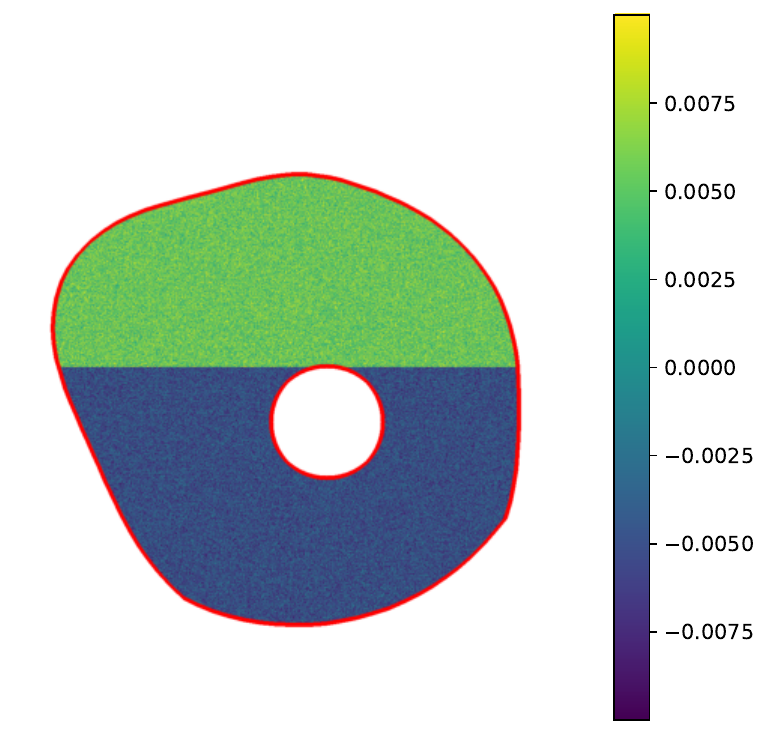}
        \caption{Noisy Input}
    \end{subfigure}
    \begin{subfigure}[H]{0.19\linewidth}
        \centering
        \includegraphics[height=\linewidth]{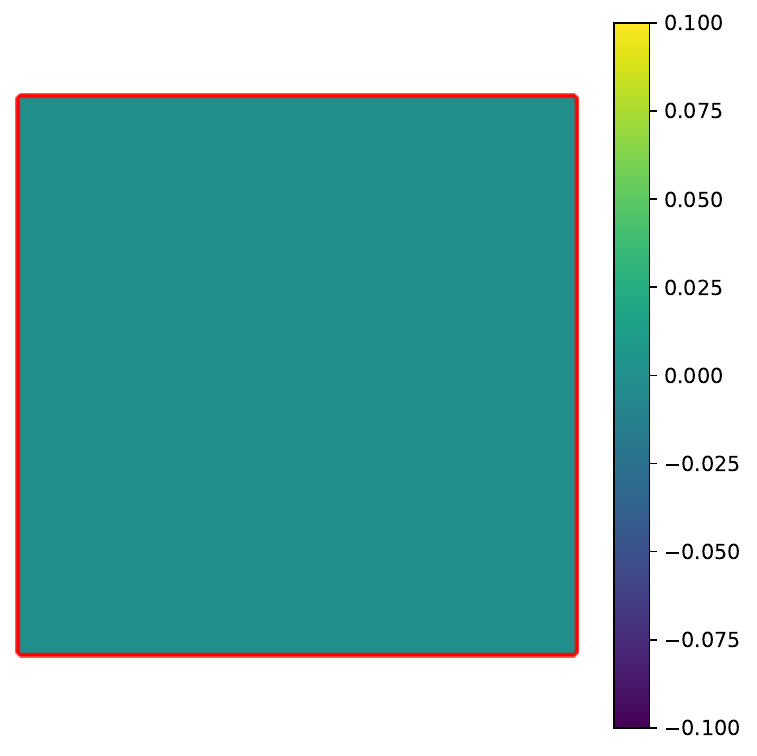}
        \caption{Lap. Square}
    \end{subfigure}
    \begin{subfigure}[H]{0.19\linewidth}
        \centering
        \includegraphics[height=\linewidth]{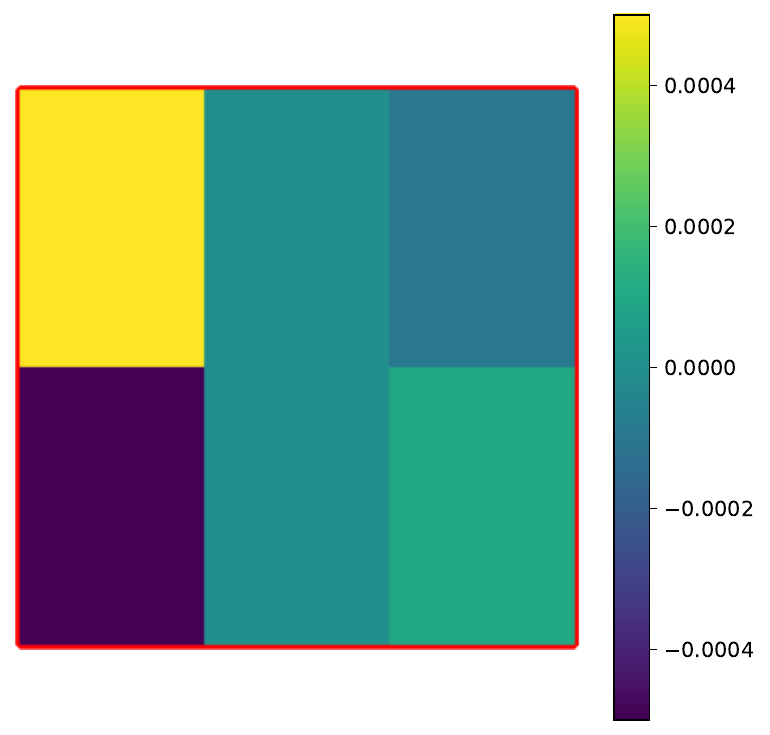}
        \caption{Poisson Square}
    \end{subfigure}
    \begin{subfigure}[H]{0.22\linewidth}
        \centering
        \includegraphics[height=100pt]{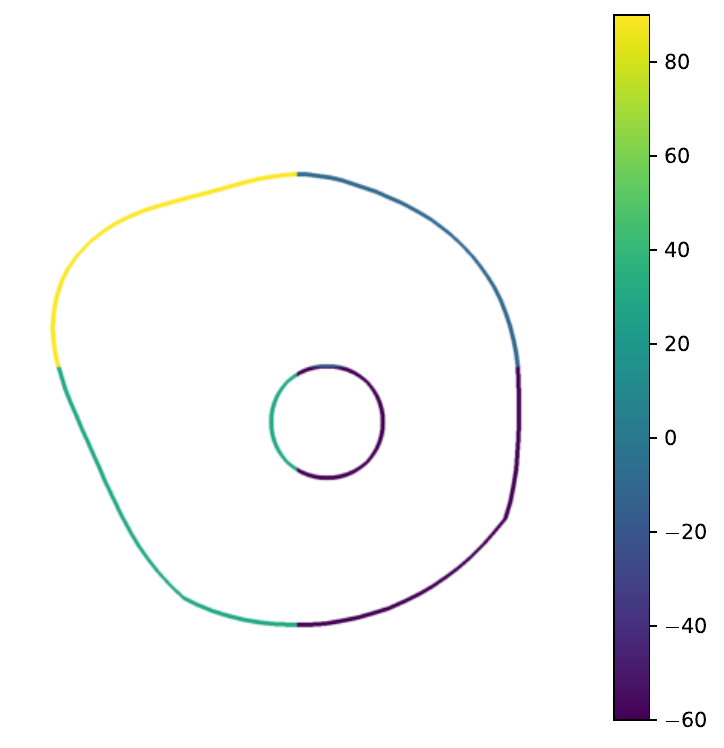}
        \caption{Boundary value for (h)}
    \end{subfigure}
    \begin{subfigure}[H]{0.33\linewidth}
        \centering
        \includegraphics[height=100pt]{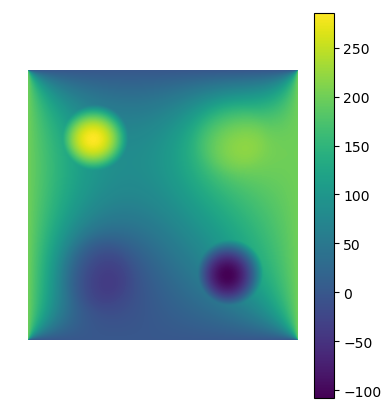}
        \caption{Ground truth ($2$D Heatmap) for (g)}
    \end{subfigure}
    \begin{subfigure}[H]{0.33\linewidth}
        \centering
        \includegraphics[height=100pt]{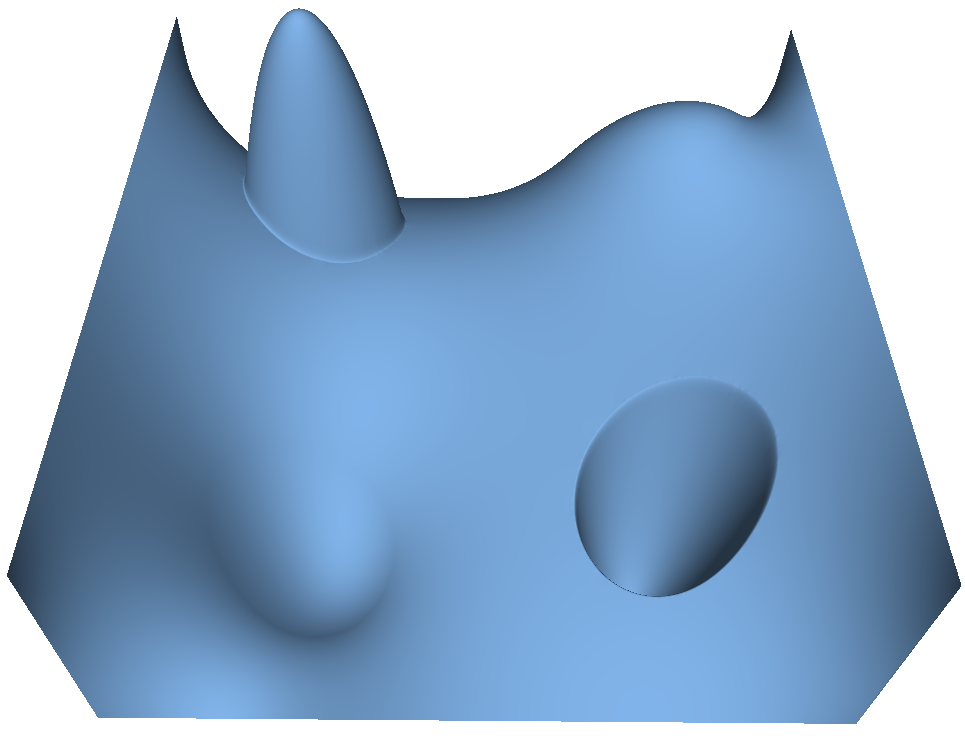}
        \caption{Ground truth ($3$D surface) for (g)}
    \end{subfigure}
    \caption{
    (a-j) illustrates the \textbf{$\vb{f}$ field 
    distributions} of our testcases. 
    The boundaries are shown in bold 
    red lines for a better view. 
    (k-m) are self-explanatory. 
    The complete illustration of \textbf{boundary values} 
    is available in 
    Sec.~\ref{sec:supplemental_material:more_specs_on_testcases}.
    }
    \label{fig:testcase}
\end{figure*}



\begin{table*}[!htb]
    \centering
    \caption{Qualitative results for 
    large-scale Poisson problems. 
    ``Time'' denotes the total time ($\mathrm{ms}$) 
    to reach relative residual errors $\le 10^{-4}$, 
    and ``Assembly'' and ``Iteration'' denotes time ($\mathrm{ms}$) for the two phases for AMG solvers;
    ``Error'' denotes 
    the final relative residual errors, 
    divided by $10^{-5}$. 
    Convergence maps are also available in the 
    appendix.
    }
    \label{tab:comparasion_1025}
    \begin{tabular}{ccccc}
        \toprule
        \textbf{Testcase} & 
        \textbf{UGrid} & 
        \textbf{AMGCL} & 
        \textbf{AmgX} & 
        \textbf{Hsieh et al.} \\
        Poisson (L) & 
        Time / Error & 
        Time / Assembly / Iteration / Error & 
        Time / Assembly / Iteration / Error & 
        Time / Error \\
        \midrule
        
        Bag & 
        \textbf{18.66} / 2.66 & 
        202.95 / 192.89 / 10.06 / 4.80 & 
        92.14 / 23.55 / 68.59 / 4.26 &
        58.09 / 420 \\ 
        
        Cat & 
        \textbf{10.09} / 2.70 & 
        270.33 / 261.75 / 8.58 / 6.63 & 
        113.92 / 25.58 / 88.35 / 6.60 & 
        49.79 / 14.6 \\ 
        
        Lock & 
        \textbf{10.55} / 9.88 & 
        140.13 / 133.06 / 7.07 / 4.05 & 
        67.60 / 17.90 / 49.69 / 4.87 & 
        49.92 / 55.78 \\ 

        N. Input & 
        \textbf{10.16} / 2.64 & 
        260.92 / 254.36 / 6.55 / 3.52 & 
        116.11 / 24.38 / 91.72 / 9.22 &
        51.07 / 2654 \\
        
        Note & 
        \textbf{10.31} / 4.06 & 
        128.96 / 121.81 / 7.16 / 3.00 & 
        64.44 / 16.55 / 47.89 / 4.80 &
        20.26 / 8.67 \\ 
        
        S. Feat. & 
        \textbf{20.01} / 3.80 & 
        424.41 / 413.66 / 10.75 / 4.13 & 
        174.76 / 38.85 / 135.92 / 3.87 &
        51.22 / 24883 \\ 
        
        L-shape & 
        \textbf{15.26} / 8.43 & 
        224.23 / 216.57 / 7.66 / 6.75 & 
        110.52 / 24.52 / 86.00 / 4.92 &
        50.53 / 96.1 \\
        
        Lap. Squ. & 
        \textbf{15.10} / 3.88 & 
        420.35 / 407.60 / 12.75 / 3.63 & 
        164.71 / 40.14 / 124.57 / 0.01 &
        31.43 / 9.03 \\ 
        
        P. Squ. & 
        \textbf{15.07} / 9.37 & 
        420.93 / 407.60 / 12.74 / 5.24 & 
        161.40 / 39.57 / 121.83 / 0.01 &
        50.57 / 974 \\
       
        Star & 
        \textbf{15.18} / 7.50 & 
        154.29 / 146.03 / 8.26 / 6.65 & 
        71.35 / 19.01 / 52.35 / 5.17 &
        50.45 / 384 \\
        
        \bottomrule
    \end{tabular}
\end{table*}

\begin{table*}[!htb]
    \centering
    \caption{Qualitative results for 
    large-scale Helmholtz problems.  
    }
    \label{tab:comparasion_1025_helmholtz}
    \begin{tabular}{cccc}
        \toprule
        \textbf{Testcase} & 
        \textbf{UGrid} & 
        \textbf{AMGCL} & 
        \textbf{AmgX} \\
        Helmholtz (L) & 
        Time / Error & 
        Time / Assembly / Iteration / Error & 
        Time / Assembly / Iteration / Error \\
        \midrule
        
        Bag & 
        \textbf{20.03} / 8.08 & 
        199.34 / 191.97 / 7.37 / 4.80 & 
        93.67 / 24.01 / 69.67 / 6.12 \\ 
        
        Cat & 
        \textbf{16.85} / 0.51 & 
        274.12 / 266.14 / 7.98 / 6.63 & 
        116.02 / 26.48 / 89.64 / 4.54 \\ 
        
        Lock & 
        \textbf{12.83} / 6.82 & 
        140.21 / 133.48 / 6.73 / 4.05 & 
        68.32 / 18.02 / 50.30 / 4.99 \\ 

        N. Input & 
        \textbf{11.98} / 6.79 & 
        247.65 / 241.00 / 6.65 / 3.51 & 
        116.98 / 24.71 / 92.28 / 9.16 \\
        
        Note & 
        \textbf{12.18} / 6.48 & 
        128.56 / 121.93 / 6.63 / 3.00 & 
        65.91 / 16.89 / 49.02 / 8.22 \\ 
        
        S. Feat. & 
        \textbf{57.68} / 9.86 & 
        412.93 / 403.02 / 9.92 / 4.13 & 
        Diverge \\ 
        
        L-shape & 
        \textbf{23.14} / 4.68 & 
        211.86 / 204.38 / 7.49 / 6.75 & 
        111.80 / 24.54 / 87.26 / 5.37 \\
        
        Lap. Squ. & 
        \textbf{44.91} / 9.98 & 
        401.24 / 388.40 / 12.84 / 3.63 & 
        165.90 / 41.32 / 124.58 / 6.91 \\ 
        
        P. Squ. & 
        \textbf{43.55} / 8.31 & 
        402.97 / 389.34 / 13.62 / 5.24 & 
        167.35 / 40.93 / 126.42 / 9.14 \\
       
        Star & 
        \textbf{15.18} / 7.50 & 
        150.84 / 142.86 / 7.98 / 6.65 & 
        72.69 / 19.53 / 53.16 / 5.78 \\
        
        \bottomrule
    \end{tabular}
\end{table*}

\begin{table*}[!htb]
    \centering
    \caption{Qualitative results for 
    large-scale steady-state 
    diffusion-convection-reaction problems.  
    }
    \label{tab:comparasion_1025_diffusion_convection}
    \begin{tabular}{cccc}
        \toprule
        \textbf{Testcase} & 
        \textbf{UGrid} & 
        \textbf{AMGCL} & 
        \textbf{AmgX} \\
        Diffusion (L) & 
        Time / Error & 
        Time / Assembly / Iteration / Error & 
        Time / Assembly / Iteration / Error \\
        \midrule
        
        Bag & 
        \textbf{41.89} / 4.77 & 
        197.55 / 190.33 / 7.22 / 5.29 & 
        104.53 / 23.83 / 80.71 / 5.12 \\ 
        
        Cat & 
        \textbf{100.68} / 9.06 & 
        273.65 / 265.80 / 7.85 / 9.21 & 
        124.11 / 26.43 / 97.69 / 5.63 \\ 
        
        Lock & 
        \textbf{58.79} / 4.78 & 
        141.53 / 134.62/ 6.91 / 4.78 & 
        141.53 / 134.62/ 6.91 / 4.78 \\ 

        N. Input & 
        \textbf{84.29} / 8.75 & 
        260.69 / 253.78 / 6.90 / 4.40 & 
        127.31 / 25.77 / 101.54 / 0.08 \\
        
        Note & 
        \textbf{25.24} / 7.42 & 
        127.19 / 121.74 / 5.45 / 6.73 & 
        61.26 / 16.73 / 44.53 / 4.78 \\ 
        
        S. Feat. & 
        \textbf{33.80} / 7.90 & 
        412.85 / 402.80 / 10.05 / 4.26 & 
        182.50 / 40.38 / 142.12 / 0.46 \\ 
        
        L-shape & 
        \textbf{30.09} / 4.70 & 
        223.98 / 216.29 / 7.69 / 6.29 & 
        111.92 / 24.76 / 87.16 / 4.57 \\
        
        Lap. Squ. & 
        \textbf{60.31} / 6.62 & 
        422.53 / 409.46 / 13.07 / 4.56 & 
        196.37 / 43.44 / 152.93 / 7.83 \\ 
        
        P. Squ. & 
        \textbf{48.60} / 7.89 & 
        418.10 / 405.11 / 12.99 / 5.11 & 
        210.96 / 48.09 / 162.86 / 5.63 \\
       
        Star & 
        \textbf{25.59} / 9.38 & 
        158.00 / 150.00 / 8.02 / 6.02 & 
        75.57 / 19.92 / 55.65 / 4.09 \\
        
        \bottomrule
    \end{tabular}
\end{table*}

\begin{table*}[!htb]
    \centering
    \caption{Qualitative results for XL-scale Poisson problems.  
    }
    \label{tab:comparasion_2049}
    \begin{tabular}{cccc}
    \toprule
        \textbf{Testcase} & \textbf{UGrid} & \textbf{AMGCL} & \textbf{AmgX} \\
        Poisson (XL) & Time / Error & Time / Assembly / Iteration / Error & Time / Assembly / Iteration / Error \\ 
        \midrule
        Bag & \textbf{36.52} / 6.95 & 893.28 / 869.92 / 23.35 / 4.64 & 280.91 / 86.22 / 194.69 / 2.80 \\ 
        Cat & \textbf{36.41} / 7.06 & 1224.02 / 1198.05 / 25.97 / 6.54 & 339.78 / 94.98 / 244.80 / 2.65 \\ 
        Lock & \textbf{36.18} / 3.09 & 642.52 / 622.94 / 19.58 / 4.47 & 191.26 / 54.01 / 137.24 / 1.99 \\ 
        N. Input & \textbf{36.29} / 6.16 & 1142.47 / 1120.98 / 21.49 / 3.61 & 322.05 / 93.25 / 228.80 / 0.04 \\ 
        Note & \textbf{36.23} / 2.70 & 568.23 / 550.23 / 18.00 / 2.62 & 168.90 / 48.42 / 120.48 / 3.17 \\ 
        S. Feat. & \textbf{89.50} / 7.32 & 1781.23 / 1748.57 / 32.67 / 5.18 & Diverge \\ 
        L-shape & \textbf{108.06} / 8.07 & 997.72 / 973.48 / 24.23 / 9.64 & 282.22 / 83.25 / 198.97 / 3.34 \\ 
        Lap. Squ. & \textbf{58.23} / 6.57 & 1782.50 / 1745.55 / 36.95 / 7.29 & 495.05 / 146.40 / 348.64 / 70.8 \\ 
        P. Squ. & \textbf{215.16} / 9.46 & 1791.26 / 1750.34 / 40.92 / 5.60 & 498.01 / 146.22 / 351.78 / 71.3 \\ 
        Star & \textbf{90.25} / 4.55 & 679.37 / 657.72 / 21.66 / 8.34 & 204.33 / 56.80 / 147.53 / 2.72 \\ 
        \bottomrule
    \end{tabular}
\end{table*}

\begin{table*}[!htb]
    \centering
    \caption{Qualitative results for XXL-scale Poisson problems.}
    \label{tab:comparasion_4097}
    \begin{tabular}{cccc}
        \toprule
        \textbf{Testcase} & \textbf{UGrid} & \textbf{AMGCL} & \textbf{AmgX} \\ 
        Poisson (XXL) & Time / Error & Time / Assembly / Iteration / Error & Time / Assembly / Iteration / Error \\ 
        \midrule
        Bag & \textbf{214.16} / 6.71 & 3646.92 / 3570.03 / 76.88 / 7.66 & 828.96 / 269.75 / 559.20 / 1.15 \\
        Cat & \textbf{210.52} / 7.56 & 4970.15 / 4881.97 / 88.18 / 8.77 & 1063.20 / 359.20 / 704.00 / 1.72 \\
        Lock & \textbf{142.23} / 4.36 & 2466.58 / 2399.21 / 67.36 / 4.89 & 617.33 / 203.30 / 414.03 / 1.28 \\
        N. Input & \textbf{145.77} / 7.62 & 4469.40 / 4396.84 / 72.55 / 4.91 & Diverge \\
        Note & \textbf{142.67} / 2.06 & 2267.48 / 2211.69 / 55.78 / 6.01 & 538.91 / 180.52 / 358.38 / 1.23 \\ 
        S. Feat. & \textbf{1131.88} / 11.2 & 7107.18 / 6960.41 / 146.76 / 4.46 & Diverge \\ 
        L-shape & \textbf{1131.77} / 14.7 & 4055.08 / 3959.96 / 95.11 / 4.17 & Diverge \\ 
        Lap. Squ. & \textbf{243.03} / 7.44 & 7194.95 / 7049.11 / 145.84 / 6.77 & 1746.84 / 598.24 / 1148.60 / 38.8 \\ 
        P. Squ. & \textbf{1137.95} / 13.9 & 7376.67 / 7183.54 / 193.13 / 5.65 & 1813.83 / 603.60 / 1210.23 / 52.3 \\ 
        Star & \textbf{1118.97} / 9.36 & 2862.27 / 2771.75 / 90.51 / 6.05 & Diverge \\ 
        \bottomrule
    \end{tabular}
\end{table*}

\textbf{Experimental Results}. 
We compare our model with 
two state-of-the-art legacy solvers, 
AMGCL \cite{Demidov18:AMGCL}, 
and NVIDIA AmgX \cite{AmgX}, 
as well as one SOTA neural solver proposed by 
\cite{Hsieh19}. 

Our testcases as shown in Fig.~\ref{fig:testcase}. 
These are all with complex geometry and topology,
and \textbf{none} of which are present in the training data, 
except the geometry of Fig.~\ref{fig:testcase} (h). 
Testcase(s) (a-f) examines the strong
generation power and robustness of UGrid
for irregular boundaries with complex
geometries and topology 
\textbf{unobserved} during the training phase; 
(g) is designed to showcase UGrid's power to handle 
both sharp and smooth features in one scene 
(note that there are two sharp-feature circles 
on the top-left and bottom-right corners, 
as well as two smooth-feature circles on 
the opposite corners);
(h) examines UGrid's robustness against noisy input 
(boundary values, boundary geometries/topology, 
and Laplacian distribution);
(i-j) are two baseline surfaces.

In Table~\ref{tab:comparasion_1025},
for Poisson problems, 
UGrid reaches the desirable precision 
$10$-$20$x faster than AMGCL, 
$5$-$10$x faster than NVIDIA AmgX, 
and $2$-$5$x faster than \cite{Hsieh19}
(note that \cite{Hsieh19} \textbf{diverged} for
most of the testcases, those cases are \textbf{not} 
counted). 
This shows the efficiency and accuracy of our method. 
Moreover, the testcase ``Noisy Input'' showcases
UGrid's robustness. 
Furthermore, among all the ten testcases, 
only (the geometry of) the ``Noisy Input'' case is observed 
in the training phase of UGrid. 
This shows that UGrid converges 
to unseen scenarios whose boundary conditions are 
of complex geometry 
(e.g., ``Cat'', ``Star'', and ``L-shape'') 
and topology (e.g., ``Note'', ``Bag'', and ``Cat''). 
On the contrary, \cite{Hsieh19} \textbf{failed} 
to converge in most of our testcases, 
which verifies one of their claimed limitations 
(i.e., \textbf{no} guarantee of convergence to unseen cases), 
and showcases the strong generalization power of our method. 
In addition, even for those converged cases,
our method is still faster than \cite{Hsieh19}. 

In Table~\ref{tab:comparasion_1025_helmholtz},
for large-scale Helmholtz problems, 
UGrid reaches the desirable precision 
$7$-$20$x faster than AMGCL and 
$5$-$10$x faster than NVIDIA AmgX. 

In Table~\ref{tab:comparasion_1025_diffusion_convection},
for large-scale steady-state 
diffusion-convection-reaction problems, 
UGrid reaches the desirable precision 
$2$-$12$x faster than AMGCL and on average 
$2$-$6$x faster than NVIDIA AmgX. 

\textbf{Scalability}. 
We further conducted two more sets of experiments
on XL and XXL Poisson problems 
\textbf{without} retraining UGrid. 
The results are available in 
Tables~\ref{tab:comparasion_2049} and \ref{tab:comparasion_4097}. 
UGrid still delivers a performance boost similar
to that observed in large-scale Poisson problems. 
These experimental results collectively validate the strong scalability of UGrid (without the need for retraining). 
Due to the page limit, results for ``small-scale'' problems are available in the appendix.

\textbf{Ablation Study}. 
The results are available in  Table~\ref{tab:comparasion_1025_ablation_main}
(more details are available in Sec.~\ref{sec:appendix:ablation_study}). 
The residual loss endows our UGrid model with 
as much as $2$x speed up versus the legacy loss. 
The residual loss also endows UGrid to
converge to the failure cases of 
its counterpart trained on legacy loss. 
These results demonstrate the claimed merits
of the residual loss. 
On the other hand, it will \textbf{diverge} 
if we naively apply 
the vanilla U-Net architecture directly to Poisson's
equations. 
This showcases the significance 
of UGrid's mathematically-rigorous
network architecture. 

\begin{table}[!htb]
    \centering
    \caption{Ablation study on large-scale Poisson problems. 
    Columns from left to right: 
    UGrid trained with residual loss, UGrid trained with legacy loss, and vanilla U-Net trained with residual loss. 
    }
    \label{tab:comparasion_1025_ablation_main}
    \begin{tabular}{cccc}
        \toprule
        \textbf{Testcase} & 
        \textbf{UGrid} & 
        \textbf{UGrid (L)} & 
        \textbf{U-Net} \\
        Poi. (L) & 
        Time / Error & 
        Time / Error & 
        Time / Error \\
        \midrule
        
        Bag & 
        \textbf{18.66} / 2.66 & 
        28.81 / 4.86 & 
        Diverge \\ 
        
        Cat & 
        \textbf{10.09} / 2.70 & 
        23.80 / 1.43 & 
        Diverge \\ 
        
        Lock & 
        \textbf{10.55} / 9.88 & 
        Diverge & 
        Diverge \\ 

        N. Input & 
        \textbf{10.16} / 2.64 & 
        20.65 / 2.42 & 
        Diverge \\ 
        
        Note & 
        \textbf{10.31} / 4.06 & 
        Diverge & 
        Diverge \\ 
        
        S. Feat. & 
        \textbf{20.01} / 3.80 & 
        31.34 / 5.14 & 
        Diverge \\ 
        
        L-shape & 
        \textbf{15.26} / 8.43 & 
        Diverge & 
        Diverge  \\ 
        
        Lap. Squ. & 
        \textbf{15.10} / 3.88 & 
        30.72 / 2.76 & 
        Diverge \\ 
        
        P. Squ. & 
        \textbf{15.07} / 9.37 & 
        31.52 / 3.33 & 
        Diverge \\ 
       
        Star & 
        \textbf{15.18} / 7.50 & 
        Diverge & 
        Diverge \\ 
        
        \bottomrule
    \end{tabular}
\end{table}

\textbf{Limitations}. 
UGrid is currently designed 
for linear PDEs \textit{only}, 
as Theorem.~\ref{thm:mask_iterator_convergence} 
does not hold for non-linear PDEs. 
Another limitation lies in the fact that 
there is no mathematical guarantee on 
\textit{how fast} UGrid will converge. 
As a consequence,  
UGrid does not necessarily 
converge faster on small-scale testcases. 
}

%
%
\section{Conclusion and Future Work}
\label{sec:conclusion}

This paper has articulated a novel 
efficient-and-rigorous neural PDE solver
built upon the U-Net and the Multigrid method, 
naturally combining the mathematical backbone 
of correctness and convergence as well as
the knowledge gained from data observations. 
Extensive experiments validate 
all the claimed advantages of our proposed approach.
Our future research efforts will be  
extending the current work to non-linear PDEs. 
The critical algorithmic barrier between our approach 
and non-linear PDEs is the limited expressiveness of 
the convolution semantics. We would like to explore 
more alternatives with stronger expressive power. 

\newpage
\textbf{Acknowledgements}. 
This work was partially supported by: (1) USA NSF IIS-1715985 and USA NSF IIS-1812606 (awarded to Hong QIN); and (2) The National Key R\&D  Program of China under Grant 2023YFB3002901, and the Basic Research Project of ISCAS under Grant ISCAS-JCMS-202303.

\section*{Impact Statement}
This paper presents work whose goal is to advance the field of Machine Learning. There are many potential societal consequences of our work, 
none which we feel must be specifically highlighted here.

\bibliography{bib}
\bibliographystyle{icml2024}

\newpage
\appendix
\onecolumn

\section{Appendix}
This supplemental material is provided to 
readers in the interest of 
our paper's theoretical and experimental completeness.

\subsection{More Specifications on Our Training Data}

\begin{figure*}[ht]
    \centering
    \begin{subfigure}[H]{0.19\linewidth}
        \centering
        \includegraphics[height=\linewidth]{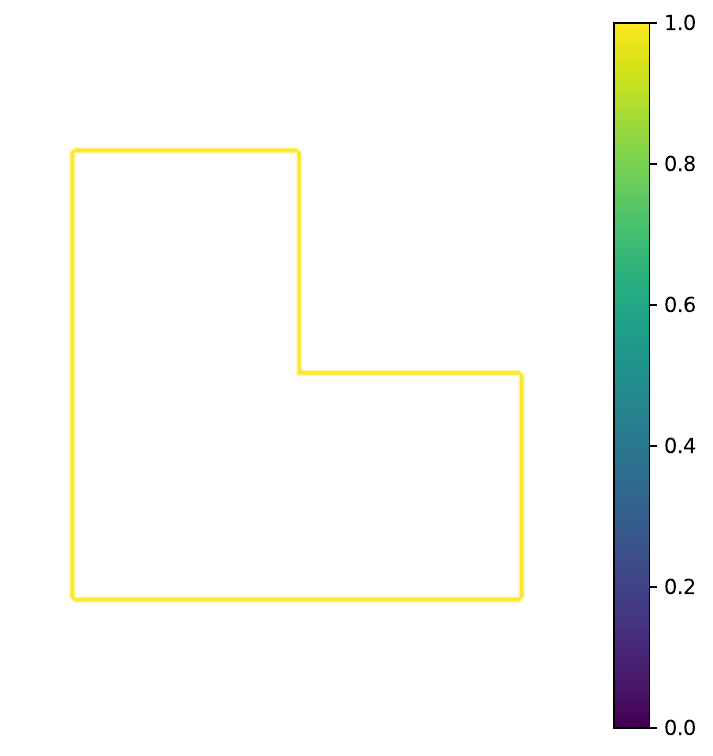}
        \caption{L-Shape}
    \end{subfigure}
    \begin{subfigure}[H]{0.19\linewidth}
        \centering
        \includegraphics[height=\linewidth]{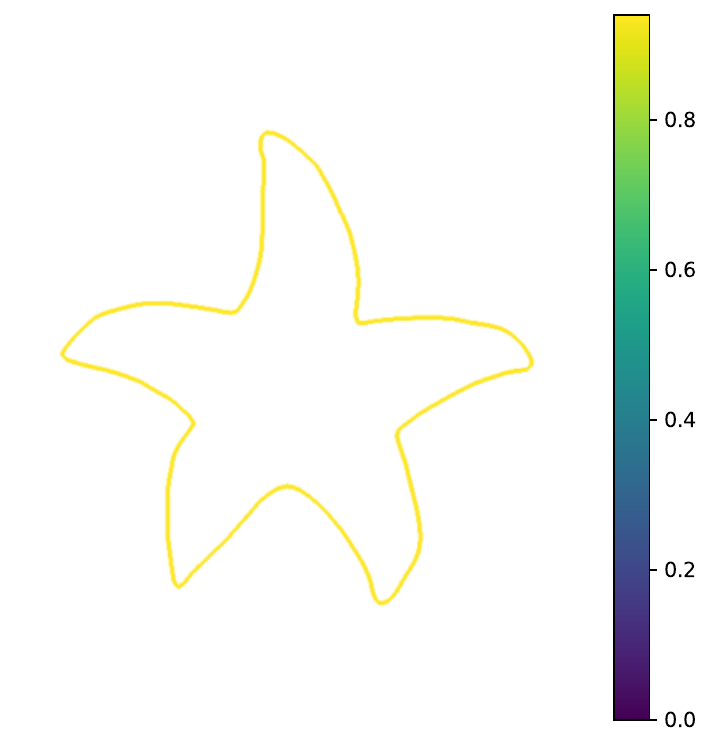}
        \caption{Star}
    \end{subfigure}
    \begin{subfigure}[H]{0.19\linewidth}
        \centering
        \includegraphics[height=\linewidth]{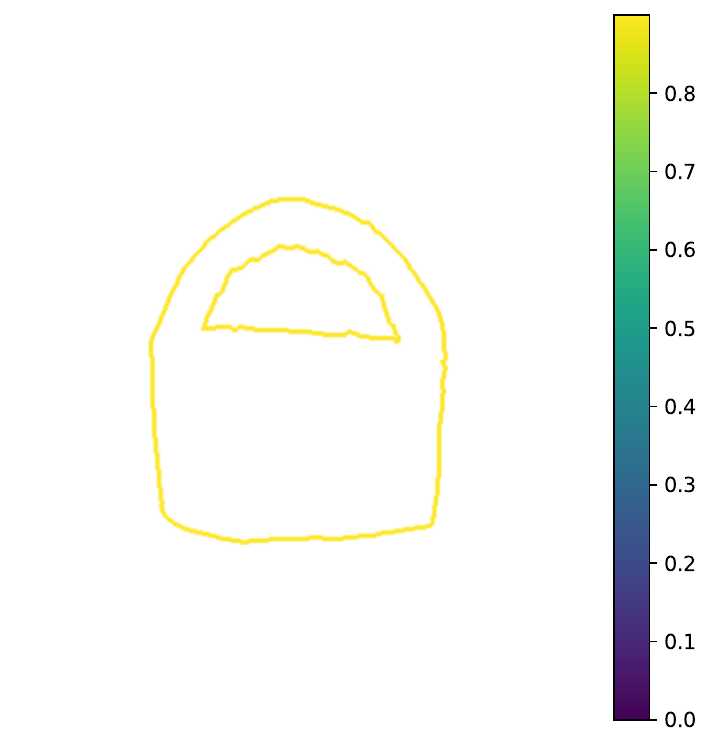}
        \caption{Lock}
    \end{subfigure}
    \begin{subfigure}[H]{0.19\linewidth}
        \centering
        \includegraphics[height=\linewidth]{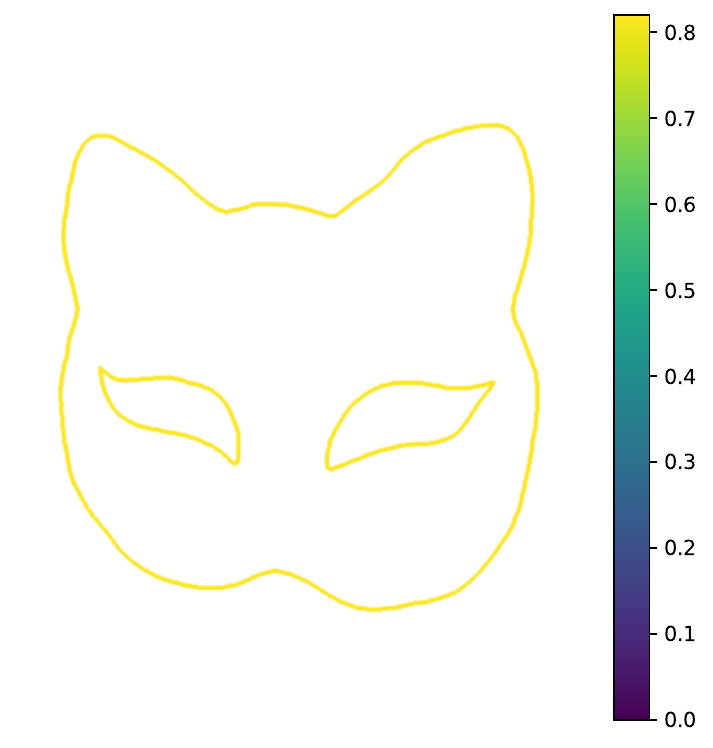}
        \caption{Cat}
    \end{subfigure}
    \begin{subfigure}[H]{0.19\linewidth}
        \centering
        \includegraphics[height=\linewidth]{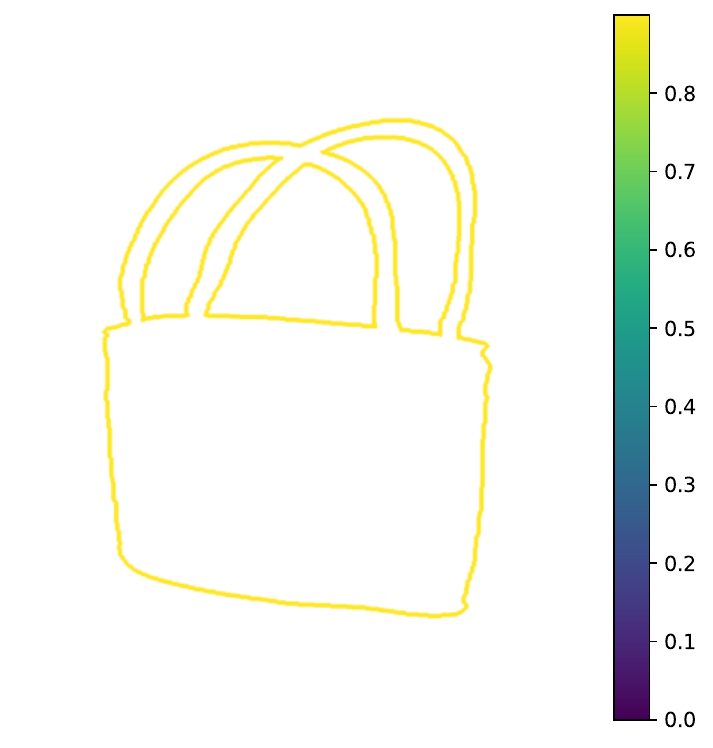}
        \caption{Bag}
    \end{subfigure}
    \begin{subfigure}[H]{0.19\linewidth}
        \centering
        \includegraphics[height=\linewidth]{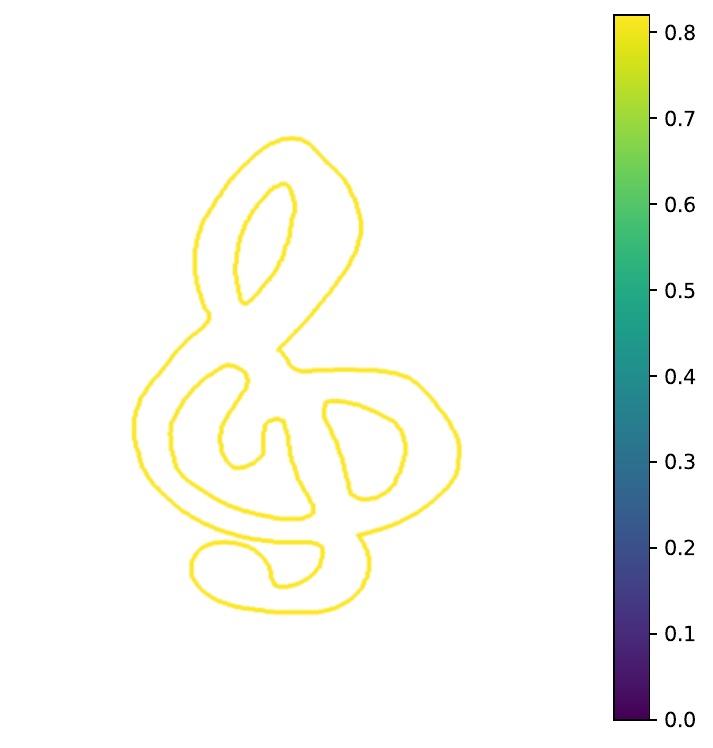}
        \caption{Note}
    \end{subfigure}
    \begin{subfigure}[H]{0.19\linewidth}
        \centering
        \includegraphics[height=\linewidth]{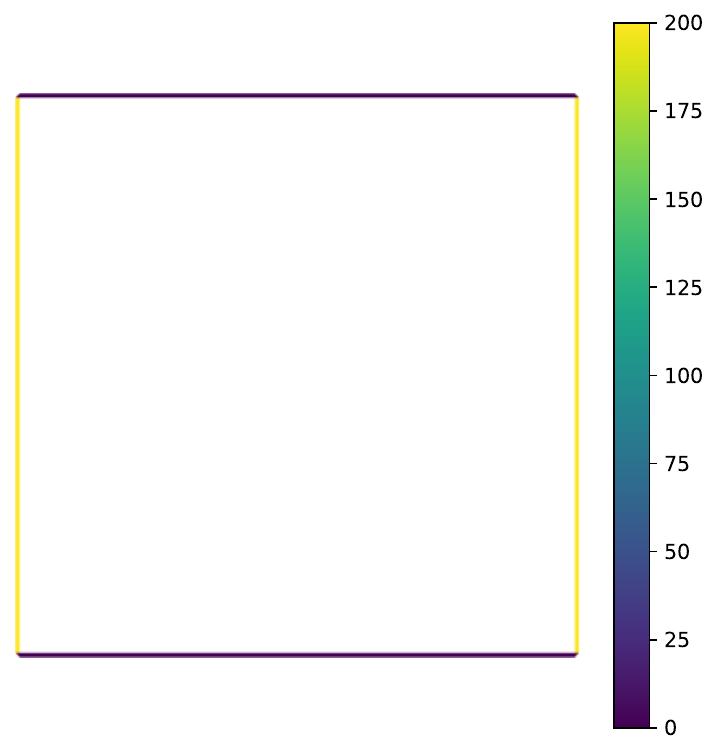}
        \caption{Sharp Feature}
    \end{subfigure}
    \begin{subfigure}[H]{0.19\linewidth}
        \centering
        \includegraphics[height=\linewidth]{imgs/testcase_boundaries/punched_curve_1025.pdf}
        \caption{Noisy Input}
    \end{subfigure}
    \begin{subfigure}[H]{0.19\linewidth}
        \centering
        \includegraphics[height=\linewidth]{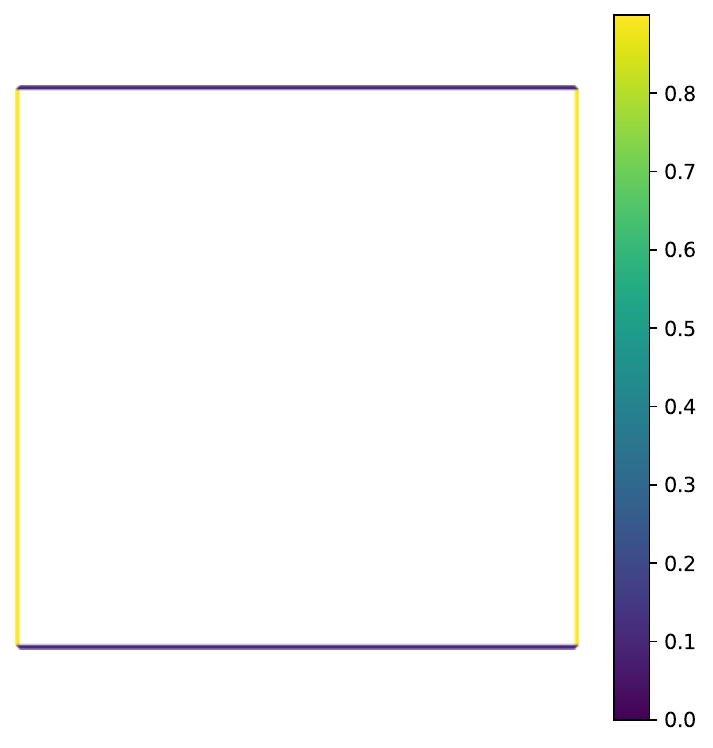}
        \caption{Lap. Square}
    \end{subfigure}
    \begin{subfigure}[H]{0.19\linewidth}
        \centering
        \includegraphics[height=\linewidth]{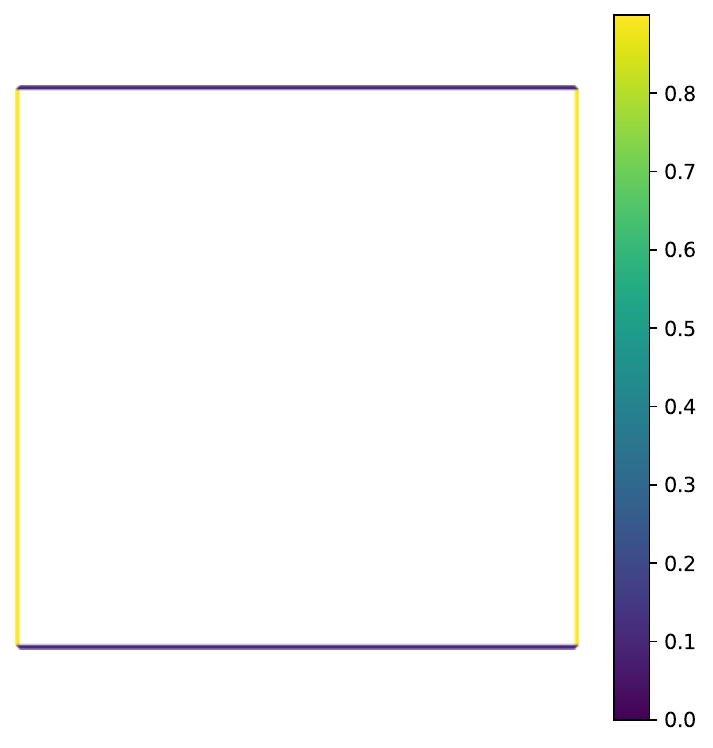}
        \caption{Poisson Square}
    \end{subfigure}
    \caption{
    Illustration of the Dirichlet boundary values
    of our ten testcases. 
    Again, all boundaries are shown in bold for a better view. 
    Note that these boundaries are \textbf{not} required 
    to be constant and could have \textit{discontinuities}, 
    which could be observed in testcases (g-j). 
    }
    \label{fig:testcase_boundaries}
\end{figure*}

UGrid is trained \textbf{only} with pairs of boundary masks
and boundary values as shown in 
Fig.~\ref{fig:testcase_boundaries} (h). 
To be specific, throughout the whole training phase, 
UGrid is exposed only to \textit{zero} $\vb{f}$-fields and  
piecewise-constant Dirichlet boundary conditions
with the ``Donut-like'' geometries.  
UGrid is \textbf{unaware} of all other 
complex geometries, topology, as well as the irregular/noisy 
distribution of boundary values/Laplacians
observed in our testcases. 
This showcases the strong generalization power of
our UGrid neural solver. 

\subsection{More Specifications on Our Testcases}
\label{sec:supplemental_material:more_specs_on_testcases}

\subsection{More Specifications on Our Qualitative Evaluations}
\label{sec:supplemental_material:more_specs_on_quals}

\textbf{Baseline Specifications}.
AMGCL is a C++ multigrid library with multiple 
GPU backends, we are comparing with its CUDA backend, 
with CG solver at the coarsest level. 
AmgX is part of NVIDIA's CUDA-X GPU-Accelerated Library, 
and we adopt its official \texttt{AMGX\_LEGACY\_CG} 
configuration. 
\cite{Hsieh19}'s code is available at
\url{https://github.com/ermongroup/Neural-PDE-Solver}.
They did not release a pre-trained model, 
so we train their model 
with configurations as-is in their 
training and data-generation scripts, 
with minimal changes to make the program run. 

\textbf{Testing}.
All of our testcases are tested for $100$ times
and the results are averaged. 
For UGrid and \cite{Hsieh19},
we set the maximum number of iterations as $64$, 
and the iteration is terminated immediately upon 
reaching this threshold, no matter whether 
the current numerical solution has reached 
the desirable precision. 
AmgX has no direct support 
for relative residual errors, 
so we set tolerance on absolute residual errors 
case-by-case to achieve similar precision. 

In our qualitative results, 
the “time” columns for AMGCL and AmgX include MG hierarchy building time
and do not include training time for UGrid. 
This is because: 
(1) Training is required only once for one type of PDEs, 
and could be ignored over the solver's lifespan; 
and (2) AMGCL/AmgX must reconstruct their MG hierarchy 
when input grid or boundary geometry changes; 
UGrid doesn't need retraining for these cases. 
When MG hierarchy is constructed already, and only RHS changes, 
AMGCL/AmgX performs better. 
However, in fields like PDE-based CAD, 
grid and boundary-geometry changes as frequently as RHS, 
and UGrid will be a better choice.

For experimental completeness, for AMGCL and AMGX, 
we also separately report the MG hierachy building time
as well as iteration time in Tables~\ref{tab:comparasion_1025}, 
\ref{tab:comparasion_1025_helmholtz} 
and \ref{tab:comparasion_1025_diffusion_convection}. 
The results validate that AMGX is slower than UGrid even without assembly time, 
and compared to UGrid, AMGCL has slightly better efficiency 
(without the assembly phase), 
yet its assembly phase takes $10$-$20$x time 
compared to UGrid's overall time consumption.

\textbf{Convergence Maps}.
For the experimental completeness of this paper, 
we also provide readers with the convergence maps of UGrid
and the three SOTA solvers we compare with. 
The convergence maps are plotted for all of our ten 
testcases, each for its two different scales. 
We further plot these convergence maps as functions of time, and
functions of iterations. 

The convergence maps for large-scale problems are as follows:

\begin{figure*}[ht]
    \centering
    \begin{minipage}[H]{\textwidth}
        \centering
            \begin{subfigure}[t]{0.19\textwidth}
            \centering
            \includegraphics[height=0.98\textwidth]{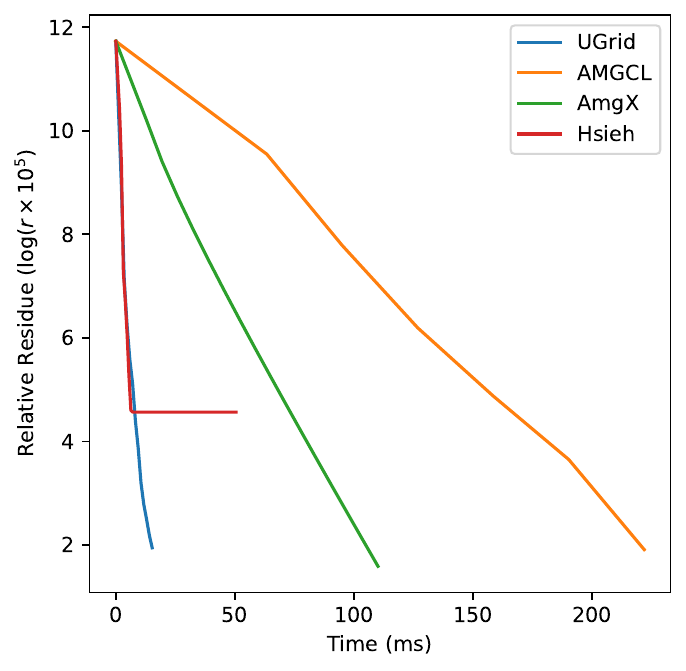}
            \caption{L-Shape}
        \end{subfigure}
        \begin{subfigure}[t]{0.19\textwidth}
            \centering
            \includegraphics[height=0.98\textwidth]{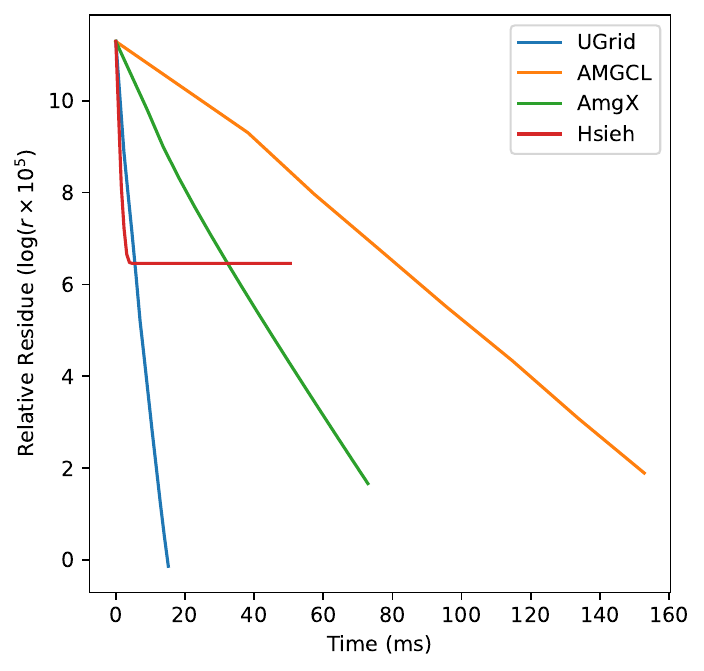}
            \caption{Star}
        \end{subfigure}
        \begin{subfigure}[t]{0.19\textwidth}
            \centering
            \includegraphics[height=0.98\textwidth]{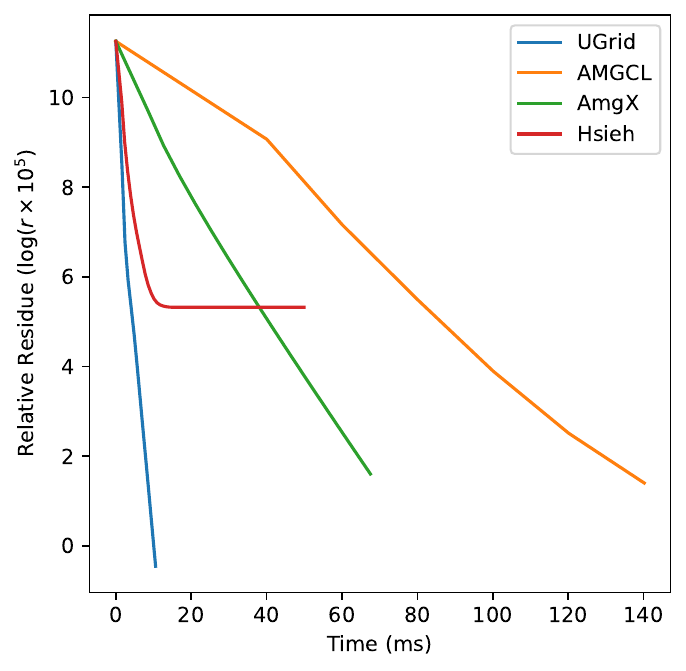}
            \caption{Lock}
        \end{subfigure}
        \begin{subfigure}[t]{0.19\textwidth}
            \centering
            \includegraphics[height=0.98\textwidth]{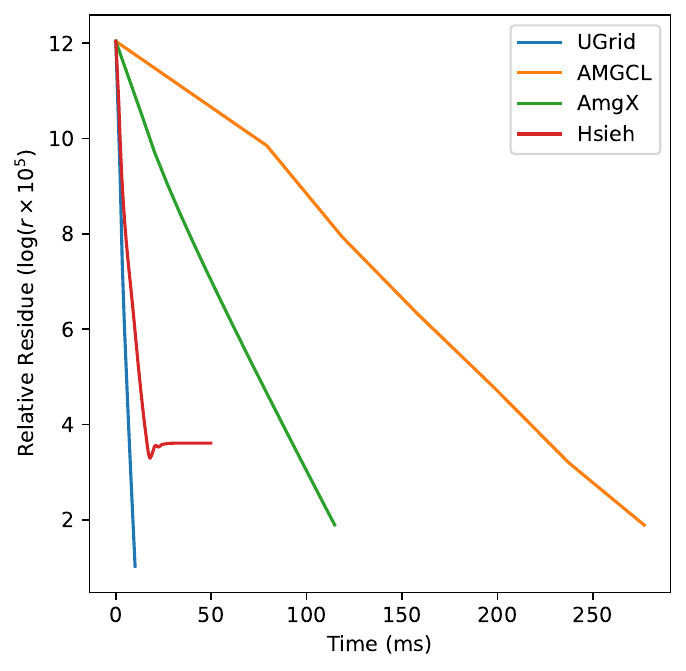}
            \caption{Cat}
        \end{subfigure}
        \begin{subfigure}[t]{0.19\textwidth}
            \centering
            \includegraphics[height=0.98\textwidth]{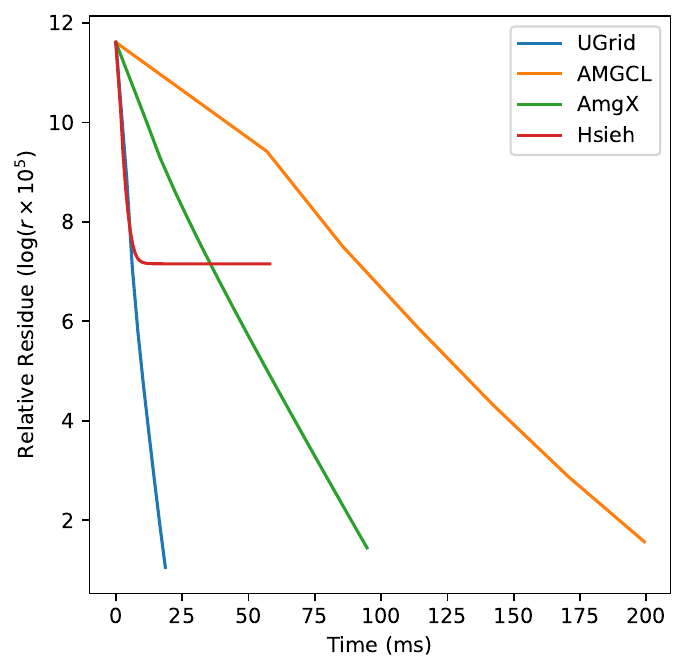}
            \caption{Bag}
        \end{subfigure}
        \begin{subfigure}[t]{0.19\textwidth}
            \centering
            \includegraphics[height=0.98\textwidth]{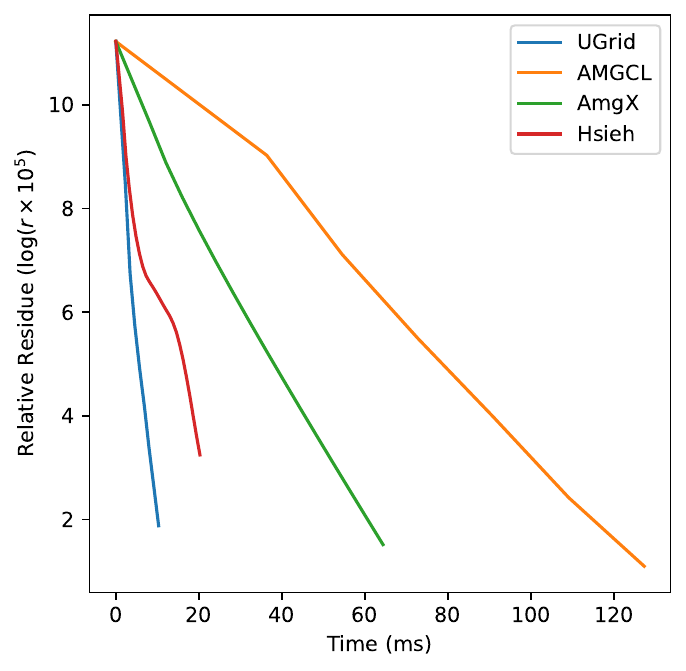}
            \caption{Note}
        \end{subfigure}
        \begin{subfigure}[t]{0.19\textwidth}
            \centering
            \includegraphics[height=0.98\textwidth]{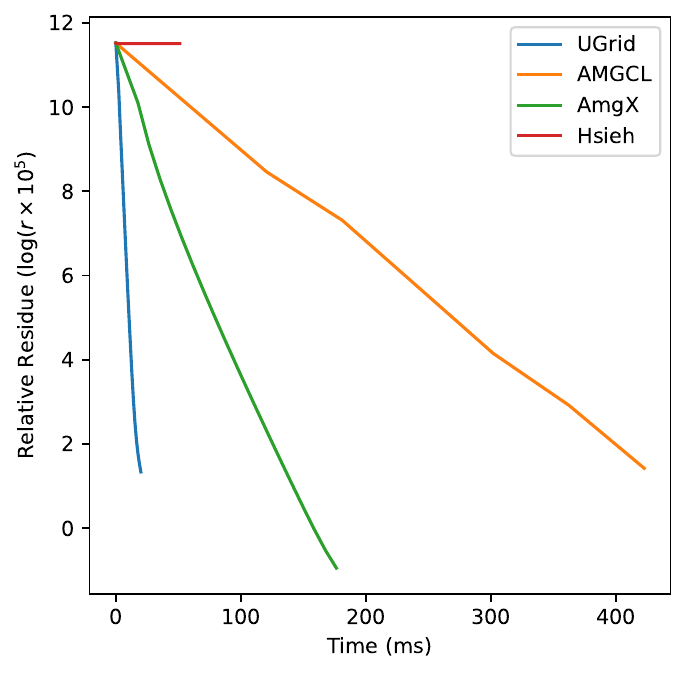}
            \caption{Sharp Feature}
        \end{subfigure}
        \begin{subfigure}[t]{0.19\textwidth}
            \centering
            \includegraphics[height=0.98\textwidth]{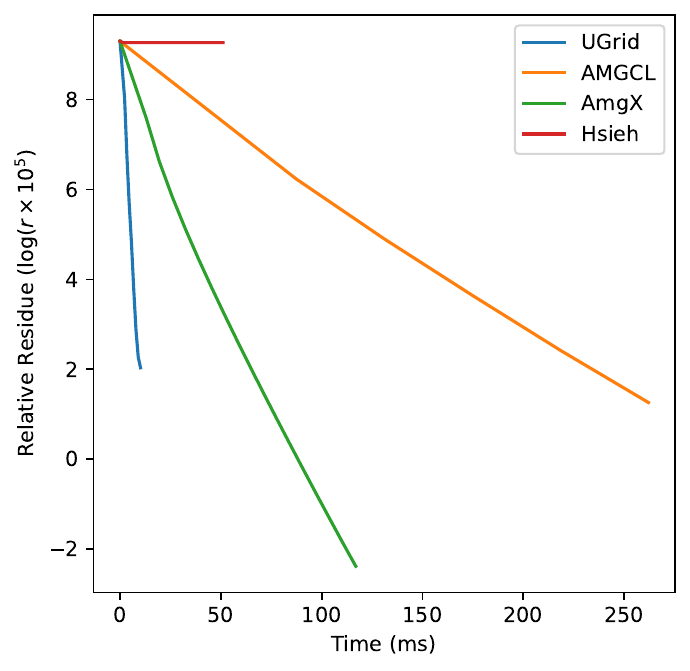}
            \caption{Noisy Input}
        \end{subfigure}
        \begin{subfigure}[t]{0.19\textwidth}
            \centering
            \includegraphics[height=0.98\textwidth]{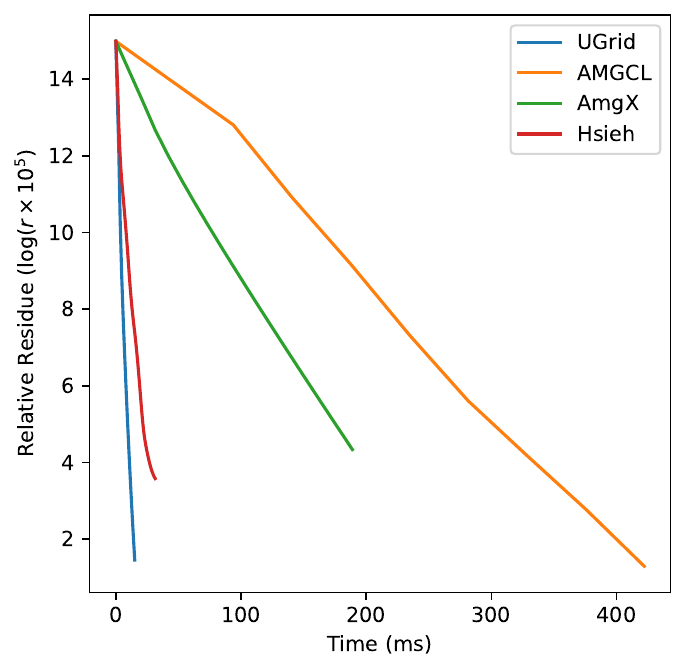}
            \caption{Lap. Square}
        \end{subfigure}
        \begin{subfigure}[t]{0.19\textwidth}
            \centering
            \includegraphics[height=0.98\textwidth]{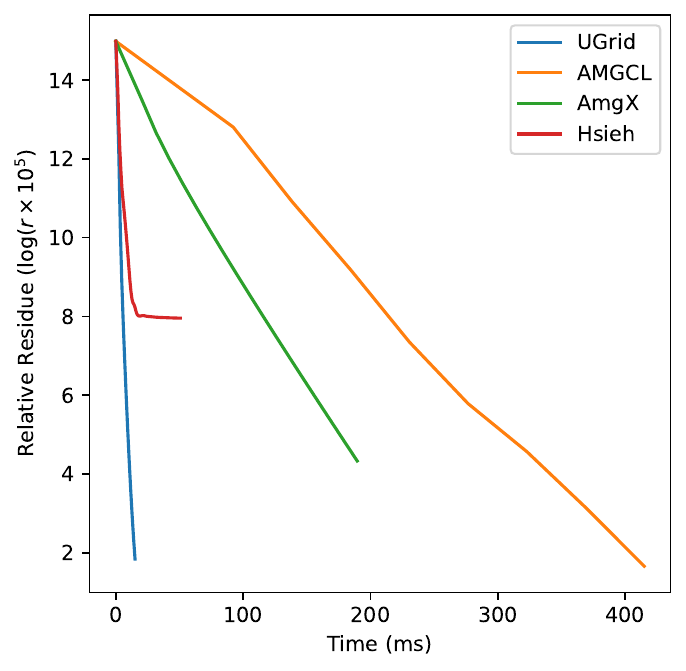}
            \caption{Poisson Square}
        \end{subfigure}
        \caption{
        Convergence map on large-scale Poisson problem. 
        The $x$ coordinates are time(s), shown in $\mathrm{ms}$; 
        the $y$ coordinates are the relative residual errors, 
        shown in logarithm ($\log(r \times 10^{5})$) for
        a better view. 
        }
        \label{fig:convergence_map_1025}
    \end{minipage}
    \begin{minipage}[H]{\textwidth}
        \centering
            \begin{subfigure}[t]{0.19\textwidth}
            \centering
            \includegraphics[height=0.98\textwidth]{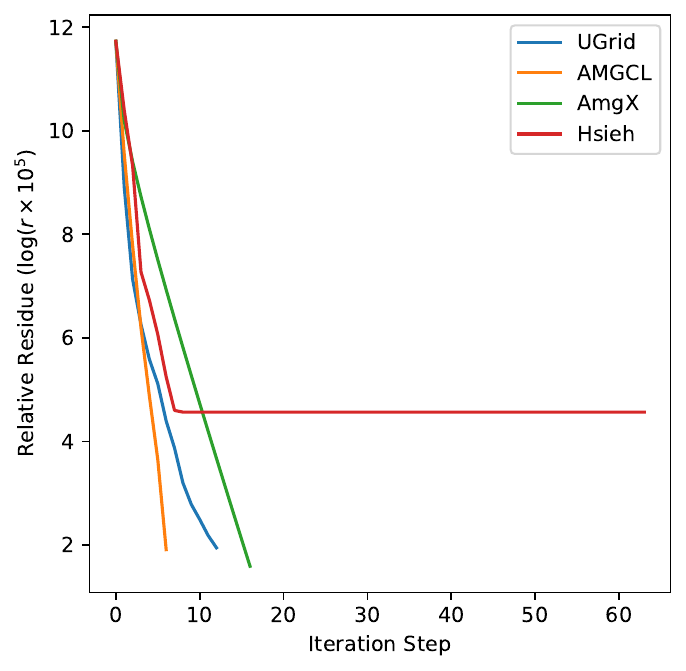}
            \caption{L-Shape}
        \end{subfigure}
        \begin{subfigure}[t]{0.19\textwidth}
            \centering
            \includegraphics[height=0.98\textwidth]{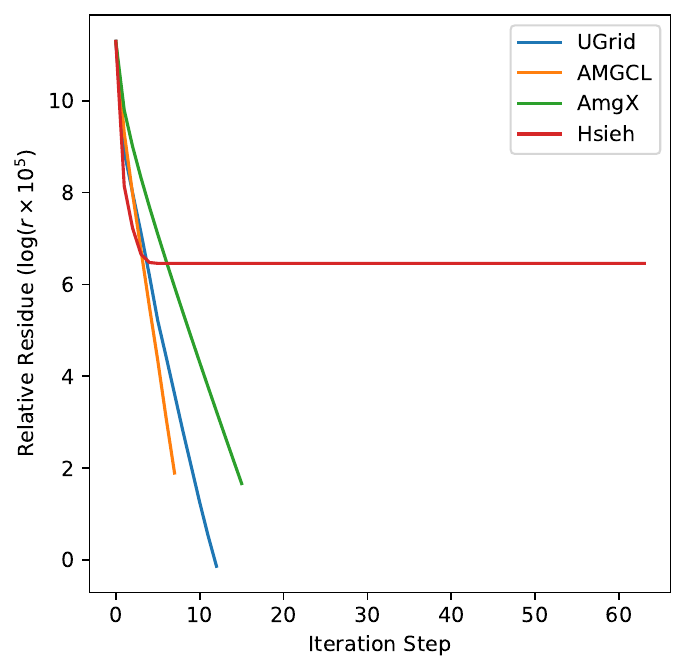}
            \caption{Star}
        \end{subfigure}
        \begin{subfigure}[t]{0.19\textwidth}
            \centering
            \includegraphics[height=0.98\textwidth]{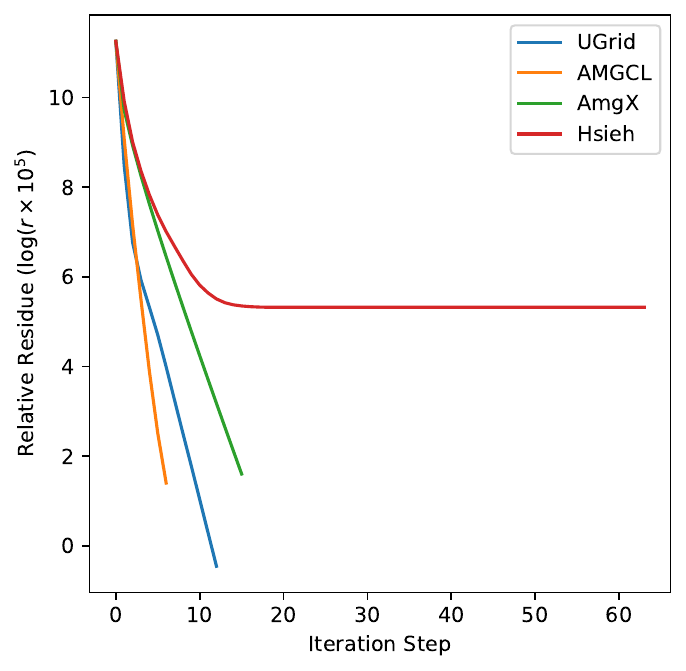}
            \caption{Lock}
        \end{subfigure}
        \begin{subfigure}[t]{0.19\textwidth}
            \centering
            \includegraphics[height=0.98\textwidth]{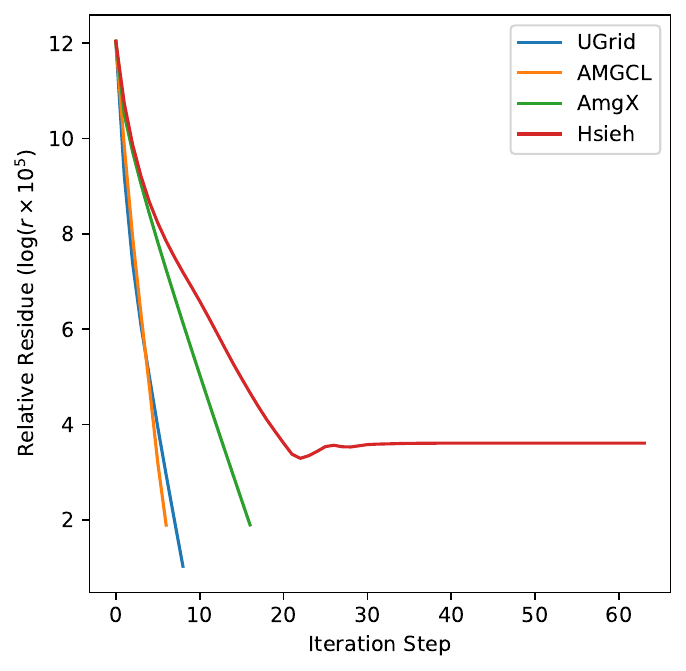}
            \caption{Cat}
        \end{subfigure}
        \begin{subfigure}[t]{0.19\textwidth}
            \centering
            \includegraphics[height=0.98\textwidth]{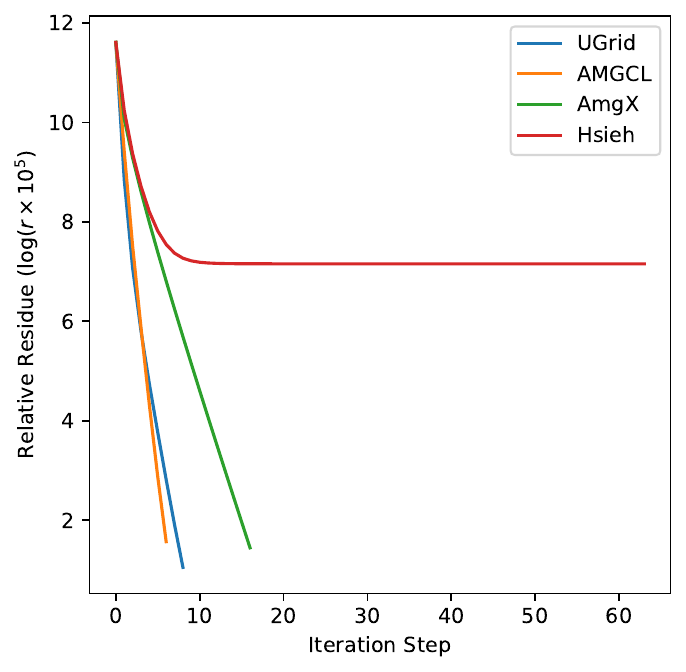}
            \caption{Bag}
        \end{subfigure}
        \begin{subfigure}[t]{0.19\textwidth}
            \centering
            \includegraphics[height=0.98\textwidth]{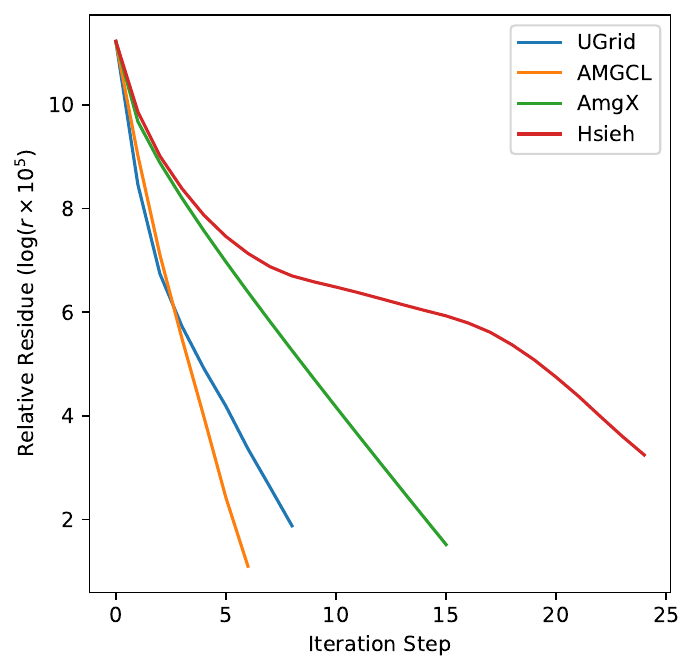}
            \caption{Note}
        \end{subfigure}
        \begin{subfigure}[t]{0.19\textwidth}
            \centering
            \includegraphics[height=0.98\textwidth]{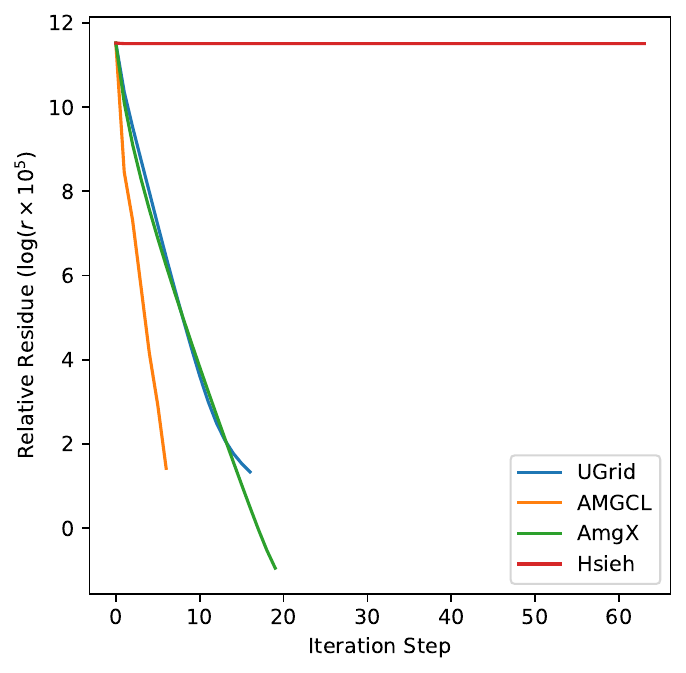}
            \caption{Sharp Feature}
        \end{subfigure}
        \begin{subfigure}[t]{0.19\textwidth}
            \centering
            \includegraphics[height=0.98\textwidth]{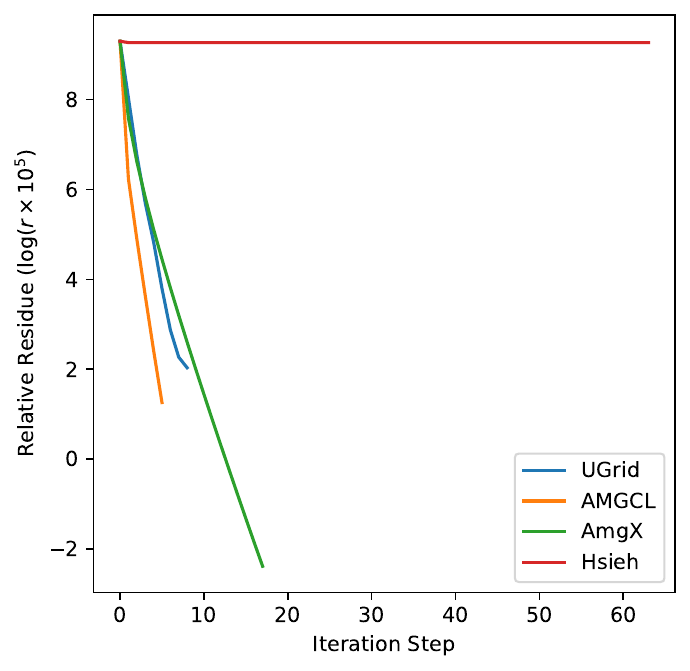}
            \caption{Noisy Input}
        \end{subfigure}
        \begin{subfigure}[t]{0.19\textwidth}
            \centering
            \includegraphics[height=0.98\textwidth]{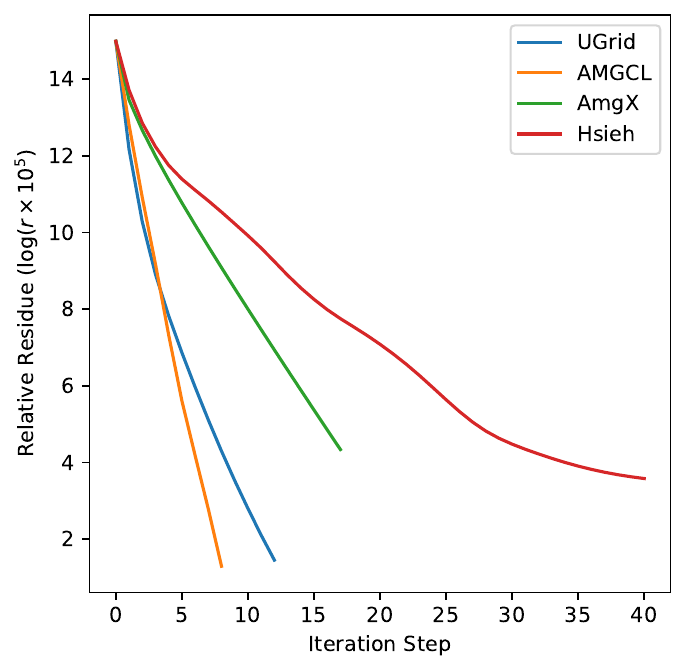}
            \caption{Lap. Square}
        \end{subfigure}
        \begin{subfigure}[t]{0.19\textwidth}
            \centering
            \includegraphics[height=0.98\textwidth]{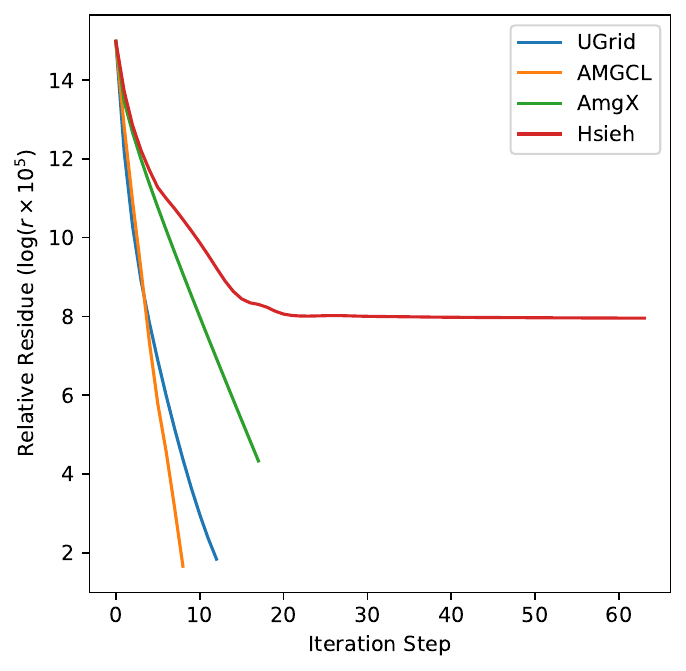}
            \caption{Poisson Square}
        \end{subfigure}
        \caption{
        Convergence map on large-scale Poisson problem. 
        The $x$ coordinates are the iteration steps; 
        the $y$ coordinates are the relative residual errors, 
        shown in logarithm ($\log(r \times 10^{5})$) for
        a better view. 
        }
        \label{fig:convergence_map_2_1025}
    \end{minipage}
\end{figure*}

\newpage
We also provide the convergence maps for small-scale problems as follows:

\begin{figure*}[ht]
    \centering
    \begin{minipage}[H]{\textwidth}
        \centering
        \begin{subfigure}[H]{0.19\linewidth}
            \centering
            \includegraphics[height=0.98\linewidth]{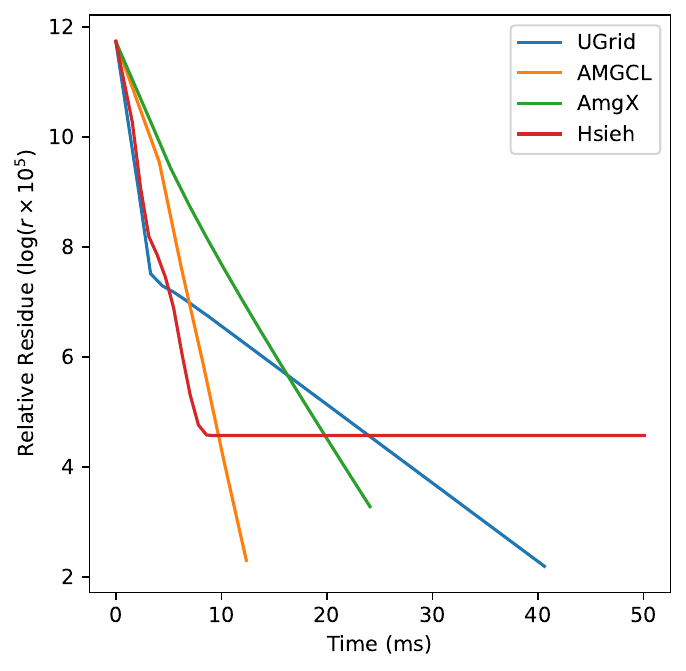}
            \caption{L-Shape}
        \end{subfigure}
        \begin{subfigure}[H]{0.19\linewidth}
            \centering
            \includegraphics[height=0.98\linewidth]{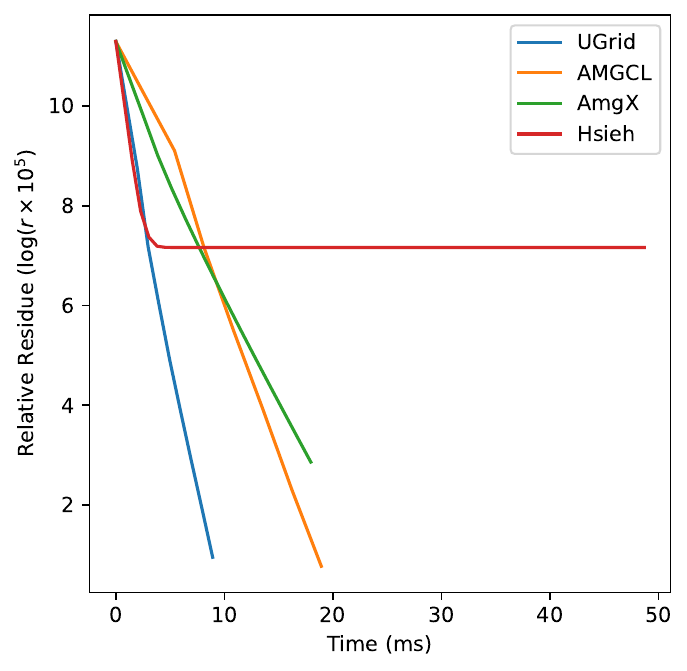}
            \caption{Star}
        \end{subfigure}
        \begin{subfigure}[H]{0.19\linewidth}
            \centering
            \includegraphics[height=0.98\linewidth]{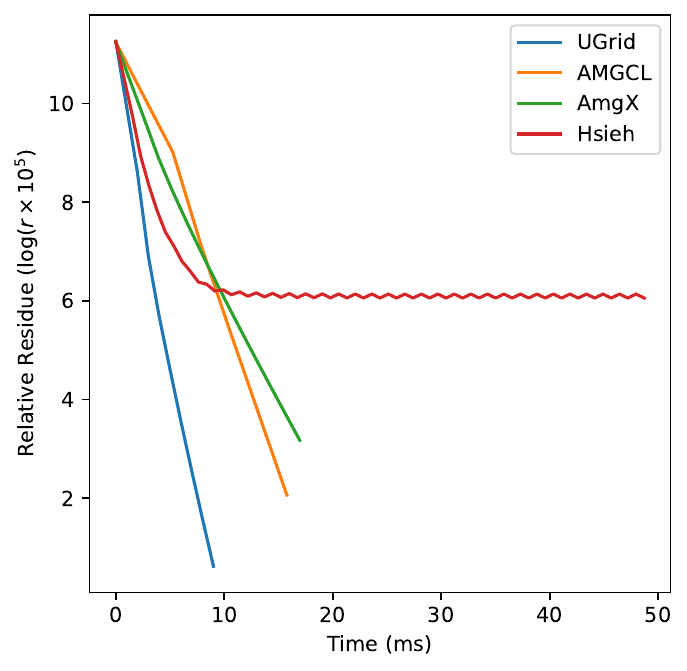}
            \caption{Lock}
        \end{subfigure}
        \begin{subfigure}[H]{0.19\linewidth}
            \centering
            \includegraphics[height=0.98\linewidth]{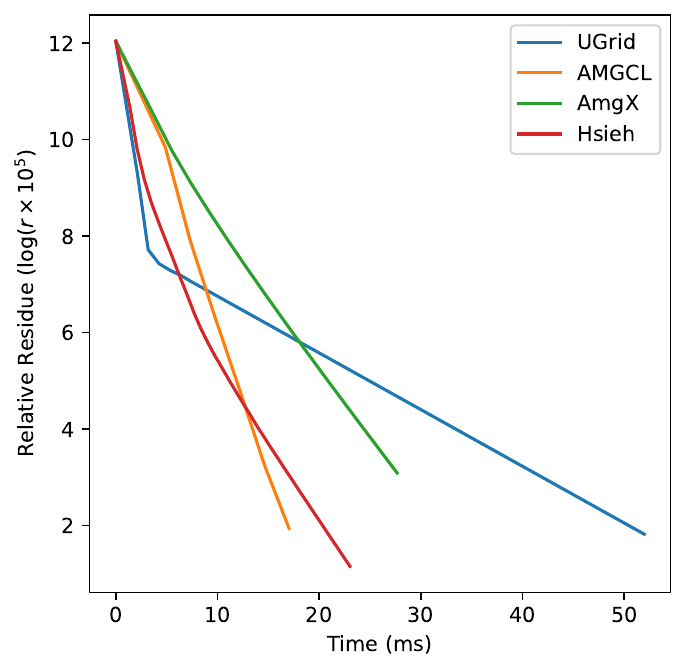}
            \caption{Cat}
        \end{subfigure}
        \begin{subfigure}[H]{0.19\linewidth}
            \centering
            \includegraphics[height=0.98\linewidth]{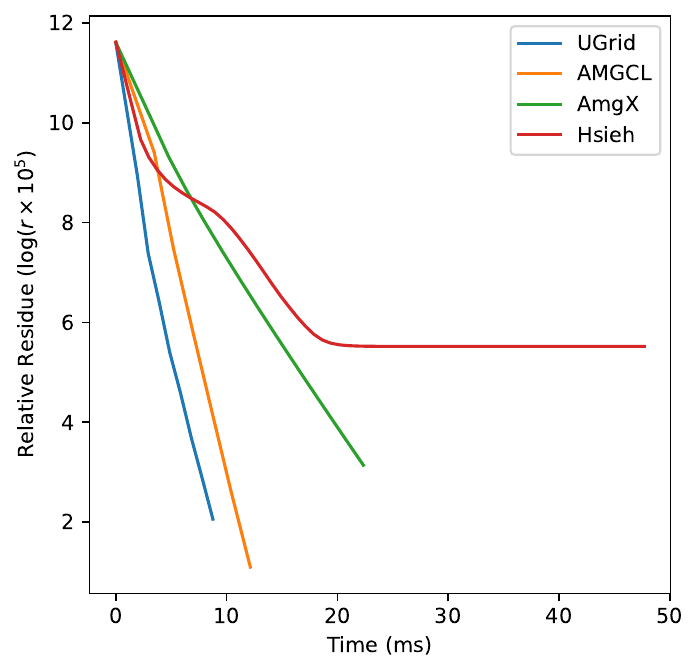}
            \caption{Bag}
        \end{subfigure}
        \begin{subfigure}[H]{0.19\linewidth}
            \centering
            \includegraphics[height=0.98\linewidth]{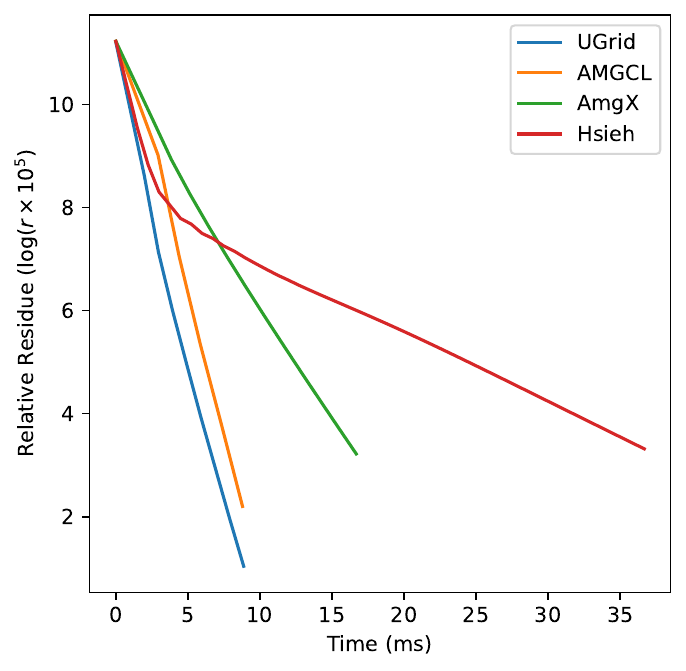}
            \caption{Note}
        \end{subfigure}
        \begin{subfigure}[H]{0.19\linewidth}
            \centering
            \includegraphics[height=0.98\linewidth]{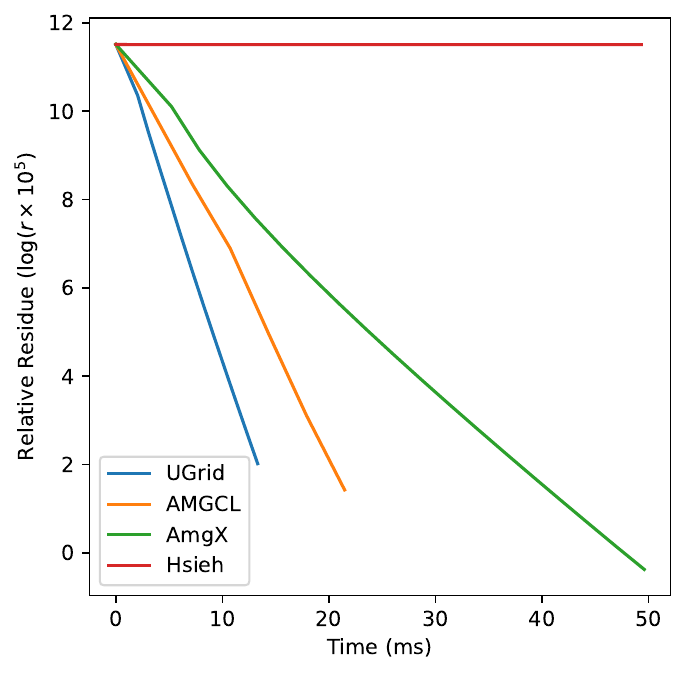}
            \caption{Sharp Feature}
        \end{subfigure}
        \begin{subfigure}[H]{0.19\linewidth}
            \centering
            \includegraphics[height=0.98\linewidth]{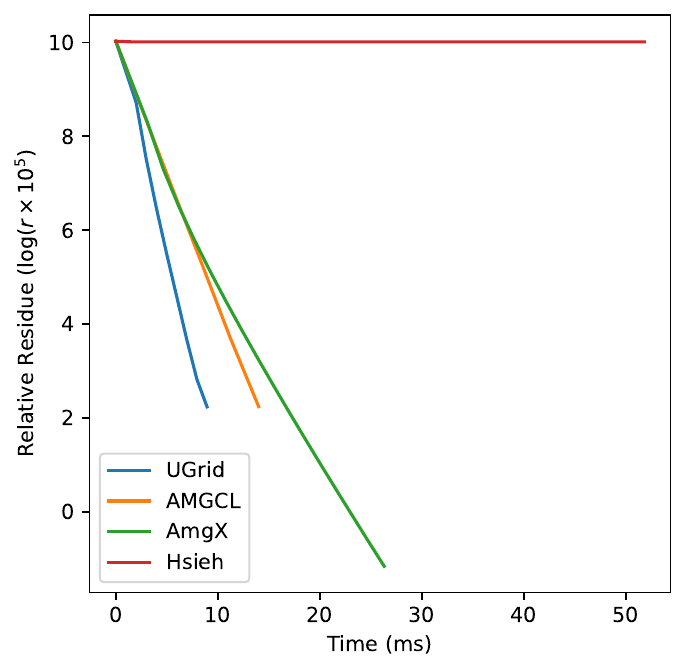}
            \caption{Noisy Input}
        \end{subfigure}
        \begin{subfigure}[H]{0.19\linewidth}
            \centering
            \includegraphics[height=0.98\linewidth]{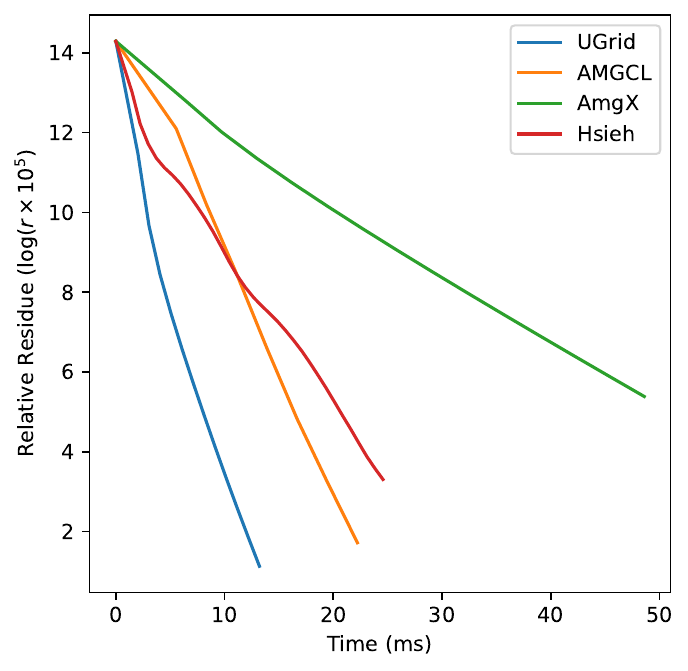}
            \caption{Lap. Square}
        \end{subfigure}
        \begin{subfigure}[H]{0.19\linewidth}
            \centering
            \includegraphics[height=0.98\linewidth]{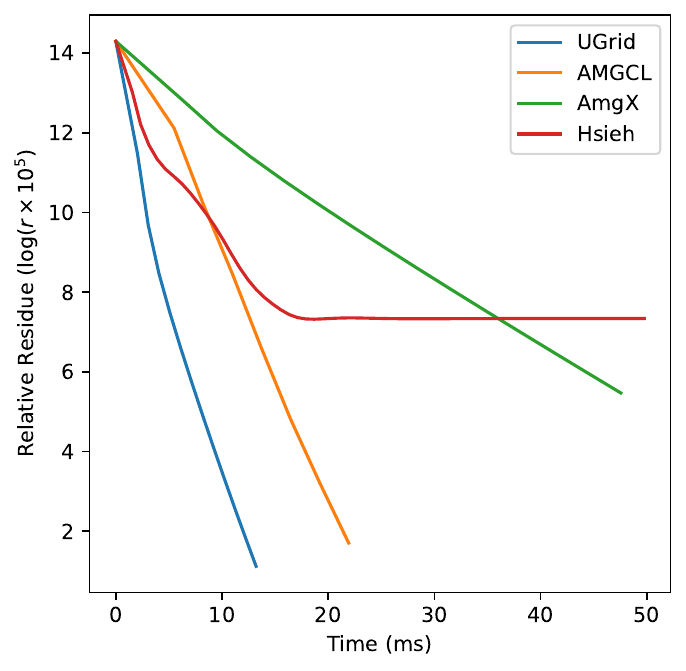}
            \caption{Poisson Square}
        \end{subfigure}
        \caption{
        Convergence map on small-scale Poisson problem. 
        The $x$ coordinates are time(s), shown in $\mathrm{ms}$; 
        the $y$ coordinates are the relative residual errors, 
        shown in logarithm ($\log(r \times 10^{5})$) for
        a better view. 
        }
        \label{fig:convergence_map_257}
    \end{minipage}
    \begin{minipage}[H]{\textwidth}
        \centering
        \begin{subfigure}[H]{0.19\linewidth}
            \centering
            \includegraphics[height=0.98\linewidth]{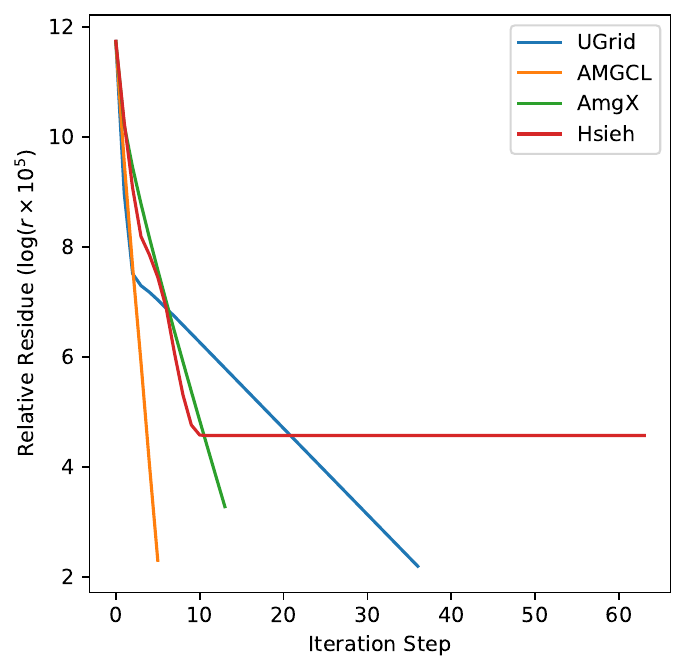}
            \caption{L-Shape}
        \end{subfigure}
        \begin{subfigure}[H]{0.19\linewidth}
            \centering
            \includegraphics[height=0.98\linewidth]{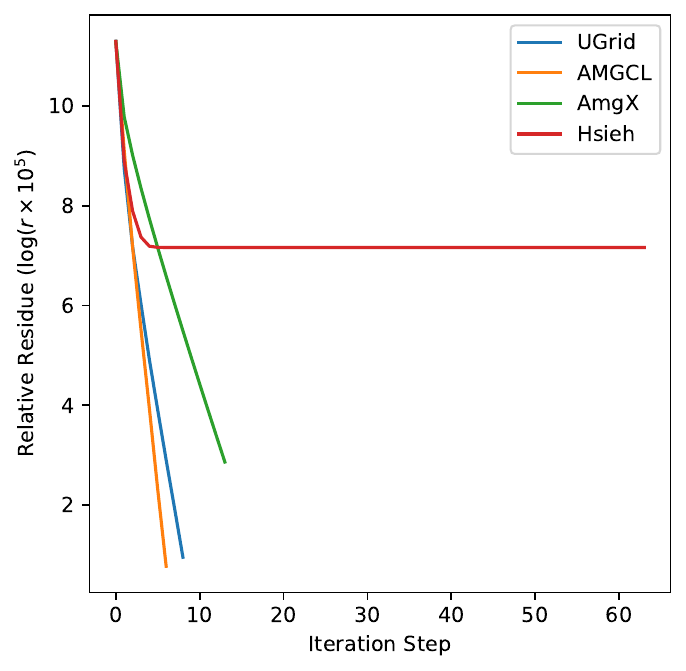}
            \caption{Star}
        \end{subfigure}
        \begin{subfigure}[H]{0.19\linewidth}
            \centering
            \includegraphics[height=0.98\linewidth]{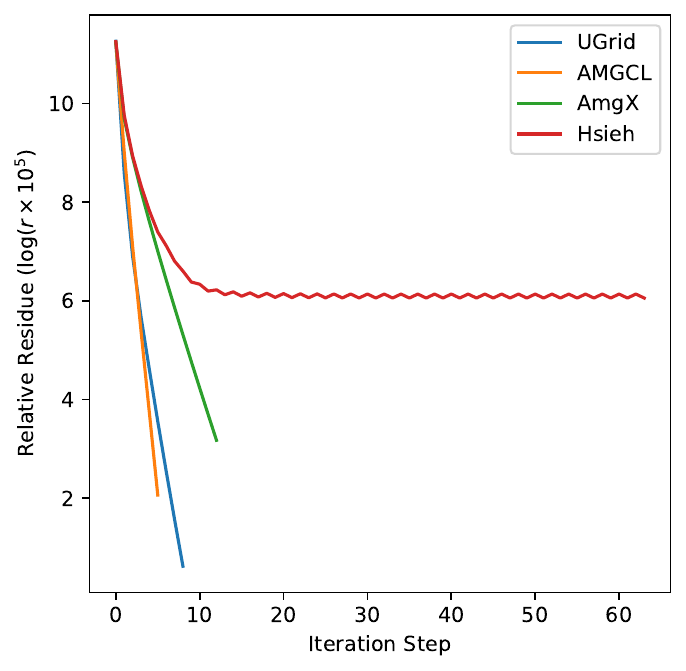}
            \caption{Lock}
        \end{subfigure}
        \begin{subfigure}[H]{0.19\linewidth}
            \centering
            \includegraphics[height=0.98\linewidth]{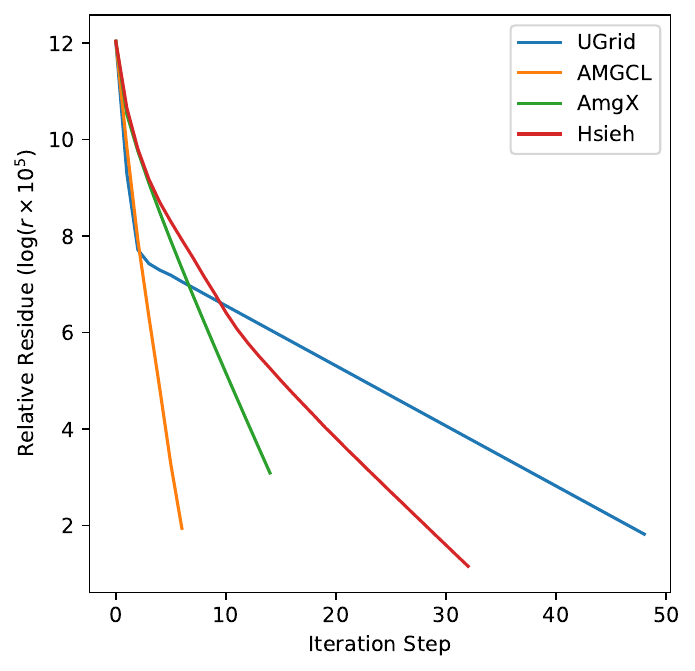}
            \caption{Cat}
        \end{subfigure}
        \begin{subfigure}[H]{0.19\linewidth}
            \centering
            \includegraphics[height=0.98\linewidth]{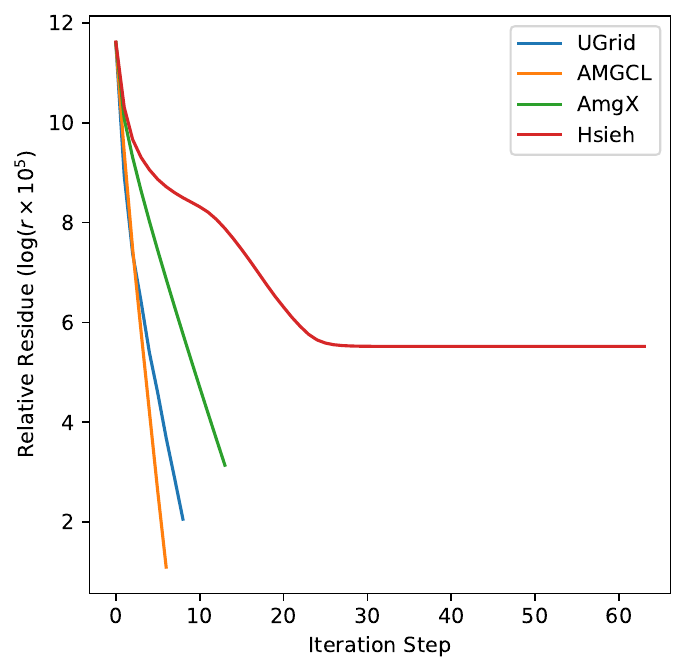}
            \caption{Bag}
        \end{subfigure}
        \begin{subfigure}[H]{0.19\linewidth}
            \centering
            \includegraphics[height=0.98\linewidth]{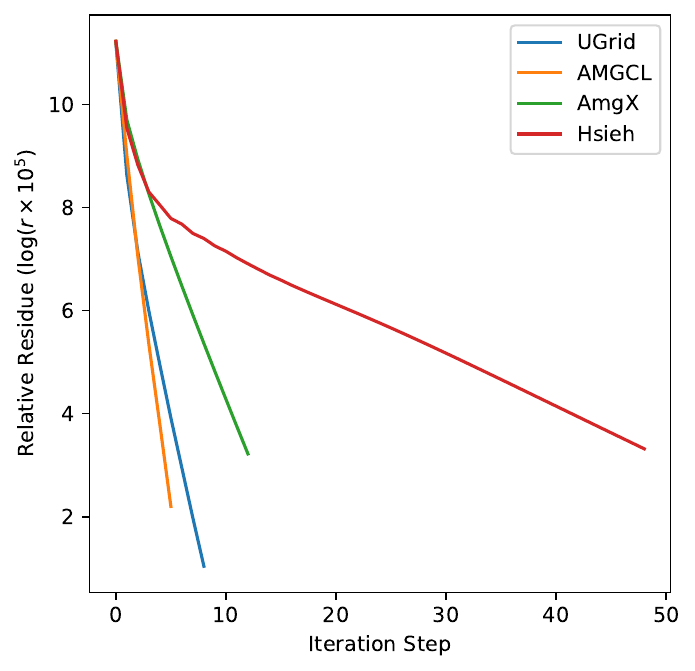}
            \caption{Note}
        \end{subfigure}
        \begin{subfigure}[H]{0.19\linewidth}
            \centering
            \includegraphics[height=0.98\linewidth]{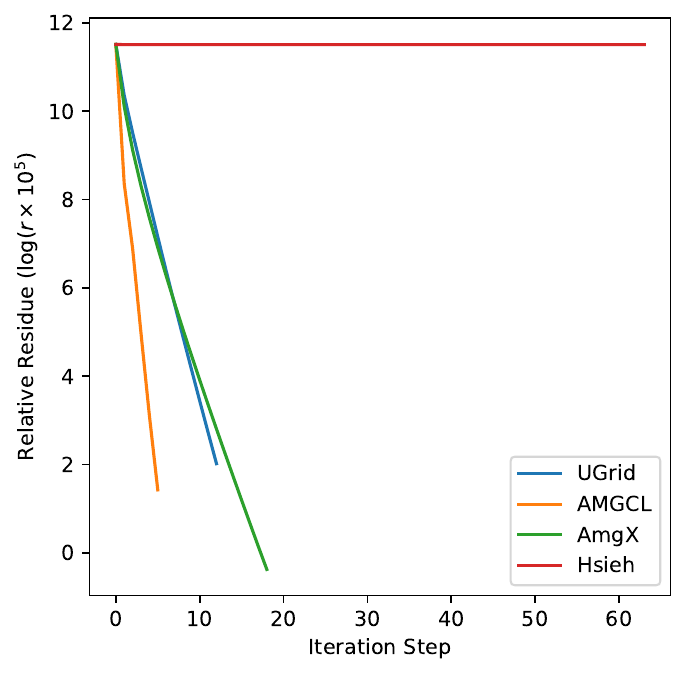}
            \caption{Sharp Feature}
        \end{subfigure}
        \begin{subfigure}[H]{0.19\linewidth}
            \centering
            \includegraphics[height=0.98\linewidth]{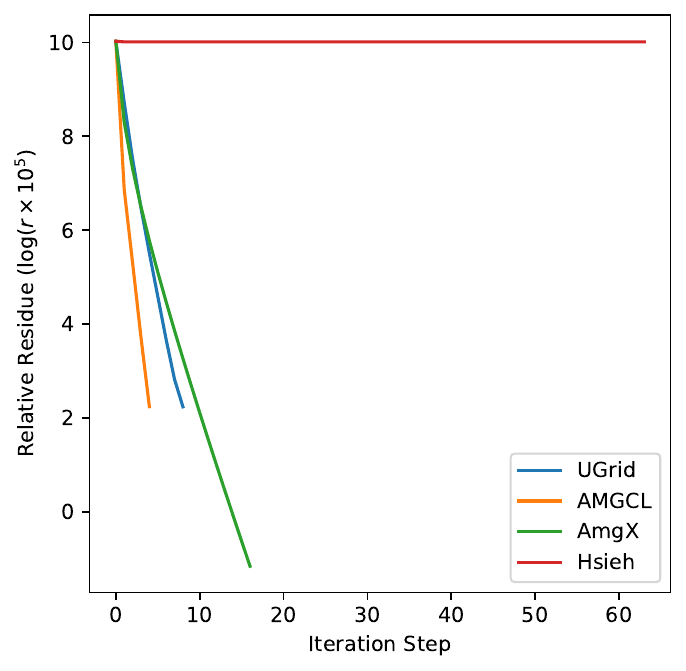}
            \caption{Noisy Input}
        \end{subfigure}
        \begin{subfigure}[H]{0.19\linewidth}
            \centering
            \includegraphics[height=0.98\linewidth]{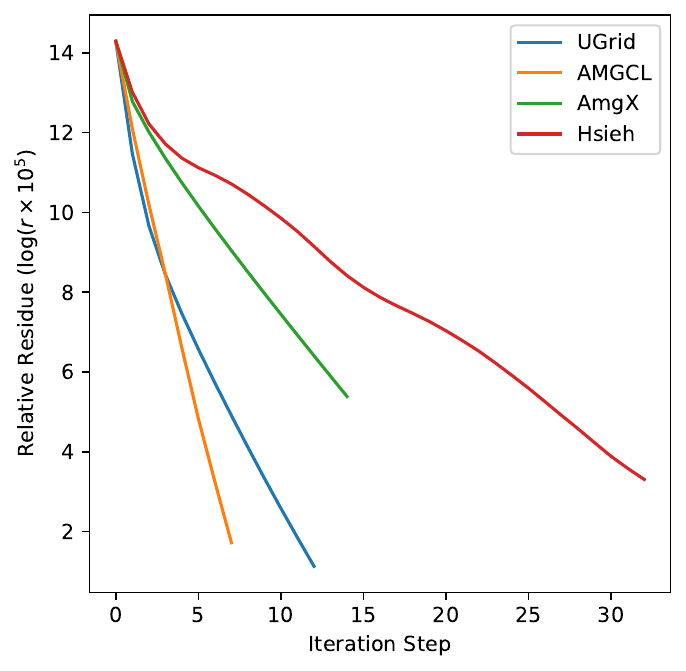}
            \caption{Lap. Square}
        \end{subfigure}
        \begin{subfigure}[H]{0.19\linewidth}
            \centering
            \includegraphics[height=0.98\linewidth]{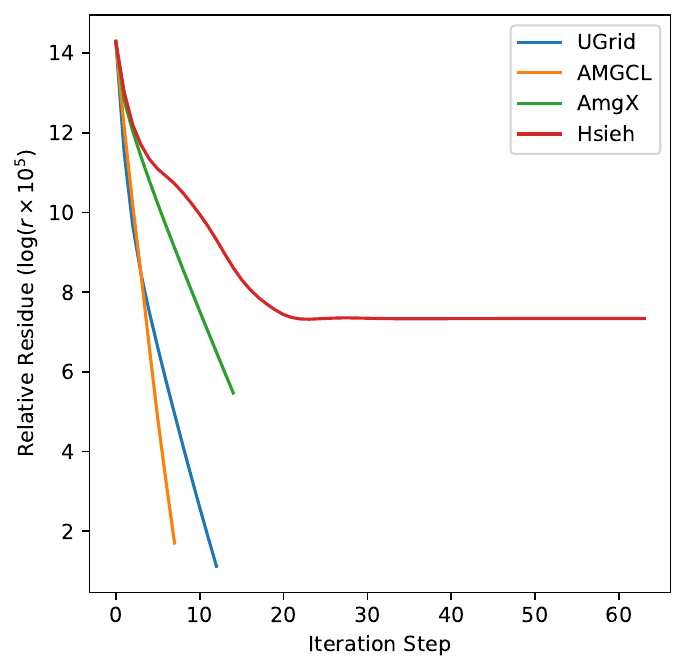}
            \caption{Poisson Square}
        \end{subfigure}
        \caption{
        Convergence map on small-scale Poisson problem. 
        The $x$ coordinates are the iteration steps; 
        the $y$ coordinates are the relative residual errors, 
        shown in logarithm ($\log(r \times 10^{5})$) for
        a better view. 
        }
        \label{fig:convergence_map_2_257}
    \end{minipage}
\end{figure*}

\subsection{Proof of Correntness of Eq.~\ref{equ:masked_iterator}}

The proof of correctness of Eq.~\ref{equ:masked_iterator} 
originates from \cite{Hsieh19}. 
For the completeness of this paper, we will also go through the mathematical proof. We start with the following Lemmas and Theorems: 

\begin{lemma}
    For a fixed linear iterator in the form of
    \begin{equation}
        \mathbf{u}_{k + 1} = \mathbf{G} \cdot \mathbf{u}_{k} + \mathbf{c} \ ,
        \label{equ:linear_iterator}
    \end{equation}
    with a square update matrice $\mathbf{G}$ 
    having a spectral radius $\rho(\mathbf{G}) < 1$, 
    $\mathbf{I} - \mathbf{G}$ is non-singular, 
    and Eq.~\ref{equ:linear_iterator} 
    converges for any constant $\mathbf{c}$ 
    and initial guess $\mathbf{u}_0$. 
    Conversely, 
    if Eq.~\ref{equ:linear_iterator} 
    converges for any $\mathbf{c}$ and $\mathbf{u}_0$, 
    then $\rho(\mathbf{G}) < 1$. 
    \label{thm:regular_iterator_convergence}
\end{lemma}
\begin{proof}
    Proved as Theorem 4.1 in \cite{Saad03}. 
\end{proof}
\begin{lemma}
    For all operator norms $\norm{\cdot}_k$, 
    $k = 1, 2, \dots, \infty$, 
    the spectral radius of a matrix $\mathbf{G}$ 
    satisfies $\rho(\mathbf{G}) \le \norm{\mathbf{G}}_k$.  
    \label{thm:operator_norm}
\end{lemma}
\begin{proof}
    Proved as Lemma 6.5 in \cite{Demmel97}. 
\end{proof}
\begin{theorem}
    For a PDE in the form of 
    Eq.~\ref{equ:masked_linear_pde}, 
    the masked iterator Eq.~\ref{equ:masked_iterator} 
    converges to its ground-truth solution
    when its prototype Jacobi iterator converges
    and $\mathbf{P}$ is full-rank diagonal. 
    \label{thm:mask_iterator_convergence_sm}
\end{theorem}
\begin{proof}
    To prove this Theorem, we only need to prove: 
    (1) 
    Eq.~\ref{equ:masked_iterator} 
    converges to a fixed point 
    $\mathbf{u}$; 
    and (2) 
    The fixed point $\mathbf{u}$ 
    satisfies Eq.~\ref{equ:masked_linear_pde}. 
    
    To prove (1), 
    we only need to prove that for the update matrix 
    \[
        \vb{G} = 
        \qty(\vb{I} - \vb{M})
        \qty(\vb{I} - \vb{P}^{-1} \vb{A})
        \ ,
    \]
    its spectral radius $\rho(\mathbf{G}) < 1$. 
    Given that the prototype Jacobi iterator
    converges, we have 
    $\rho \paren{\mathbf{I} - \mathbf{P}^{-1} \mathbf{A}} < 1$. 
    From Lemma~\ref{thm:operator_norm}, 
    taking the \textit{spectral norm} 
    $\norm{\cdot}_2$ (i.e., $k = 2$),
    we have
    \[
        \rho(\vb{G}) \le 
        \norm{ 
            \qty(\vb{I} - \vb{M})
            \qty(\vb{I} - \vb{P}^{-1} \vb{A}) 
        }_2 \le 
        \norm{ \vb{I} - \vb{M} }_2 \, 
        \norm{ \vb{I} - \vb{P}^{-1} \vb{A} }_2
        \ .
    \]
    Furthermore, 
    because $\mathbf{I} - \mathbf{P}^{-1} \mathbf{A}$ 
    is symmetric, we have 
    $\rho \paren{ \mathbf{I} - \mathbf{P}^{-1} \mathbf{A} } = 
    \norm{ \mathbf{I} - \mathbf{P}^{-1} \mathbf{A} }_2$. 
    On the other hand, because 
    $\mathbf{I} - \mathbf{M} \in \{ 0, 1 \}^{n^2 \times n^2}$ 
    is a binary diagonal matrix, we have
    $\norm{\mathbf{I} - \mathbf{M}}_2 = 1$. 
    This yields $\rho(\mathbf{G}) < 1$.
    
    To prove (2), we first notice that the fixed point 
    $\mathbf{u} = \mathbf{u}_{k + 1} = \mathbf{u}_k$ 
    of Eq.~\ref{equ:masked_iterator} satisfies
    \begin{align*}
        \mathbf{u} & = 
        \paren{ \mathbf{I} - \mathbf{M} }  
        \paren{ \paren{ \mathbf{I} - \mathbf{P}^{-1} \mathbf{A} } 
                \mathbf{u} +   
                \mathbf{P}^{-1} \mathbf{f} } + 
        \mathbf{M} \mathbf{b} \notag
        \ , \quad \text{i.e., } \\   
        \paren{ \mathbf{I} - \mathbf{M} } \mathbf{u} + 
        \mathbf{M} \mathbf{u} & = 
        \paren{ \mathbf{I} - \mathbf{M} }  
        \paren{ \paren{ \mathbf{I} - \mathbf{P}^{-1} \mathbf{A} } 
                \mathbf{u} +   
                \mathbf{P}^{-1} \mathbf{f} } + 
        \mathbf{M} \mathbf{b} \ .
    \end{align*}
    Again, since 
    $\mathbf{M} \in \{ 0, 1 \}^{n^2 \times n^2}$ 
    is a binary diagonal matrix, we have
    \begin{equation}
        \left\{
        \begin{aligned}
            \paren{ \mathbf{I} - \mathbf{M} } \mathbf{u} & = 
            \paren{ \mathbf{I} - \mathbf{M} } 
            \paren{ \paren{ \mathbf{I} - \mathbf{P}^{-1} \mathbf{A} } 
                    \mathbf{u} +   
                    \mathbf{P}^{-1} \mathbf{f} } \\ 
            \mathbf{M} \mathbf{u} & =
            \mathbf{M} \mathbf{b} 
        \end{aligned}
        \right. \ .
        \label{equ:compare_elements}
    \end{equation}
    The second equation in Eq.~\ref{equ:compare_elements} 
    is essentially the second equation
    in Eq.~\ref{equ:masked_linear_pde}. 
    Furthermore, the first equation in 
    Eq.~\ref{equ:compare_elements} 
    could be simplified into
    $\qty(\vb{I} - \vb{M}) 
    \vb{P}^{-1} 
    \qty(\vb{A} \vb{u} - \vb{f})
    = \vb{0}$.
    Since $\vb{P}$ is full-rank diagonal, 
    $\vb{P}^{-1}$ should also be full-rank diagonal. 
    Then we have
    $\qty(\vb{I} - \vb{M}) \qty(\vb{A} \vb{u} - \vb{f}) = \vb{0}$,
    which means that $\vb{u}$ also satisfies 
    the first equation in Eq.~\ref{equ:masked_linear_pde}. 
\end{proof}

\subsection{Proof of Theorem.~\ref{thm:legacy_loss}}
\label{sec:supplemental_material:legacy_loss}

\begin{proof}
    Denote $\varepsilon_{\mathbf{x}}$ as $\mathbf{x}$'s 
    relative residual error, then we have: 
    \begin{align*}
        \varepsilon_{\mathbf{x}} 
        & = \dfrac{\norm{\widetilde{\mathbf{f}} - 
                         \widetilde{\mathbf{A}} \, \mathbf{x}}}
                  {\norm{\widetilde{\mathbf{f}}}} 
        \approx \dfrac{\norm{\widetilde{\mathbf{f}} - 
                               \widetilde{\mathbf{A}} \, 
                               \paren{\mathbf{y} \pm 
                                      l_{\mathrm{max}} \,
                                      \mathbf{y}}}}
                        {\norm{\widetilde{\mathbf{f}}}} 
        = \dfrac{\norm{\paren{\widetilde{\mathbf{f}} - 
                          \widetilde{\mathbf{A}} \, 
                          \mathbf{y}}
                         \mp
                         \paren{l_{\mathrm{max}} \, 
                                \widetilde{\mathbf{A}} \, 
                                \mathbf{y}}}}
                  {\norm{\widetilde{\mathbf{f}}}} \notag \\ 
        & \le \dfrac{\norm{\widetilde{\mathbf{f}} - 
               \widetilde{\mathbf{A}} \, \mathbf{y}}}
              {\norm{\widetilde{\mathbf{f}}}}
              + 
              \dfrac{\norm{l_{\mathrm{max}} \, 
                         \widetilde{\mathbf{A}} \, 
                         \mathbf{y}}}
                    {\norm{\widetilde{\mathbf{f}}}} 
        = \varepsilon_{\mathbf{y}}
              + 
              \dfrac{\norm{l_{\mathrm{max}} \, 
                         \widetilde{\mathbf{A}} \, 
                         \mathbf{y}}}
                    {\norm{\widetilde{\mathbf{f}}}} \ , 
    \end{align*} 
    where $\varepsilon_{\mathbf{y}}$ denotes the relative residual
    error of the ``ground-truth'' value $\mathbf{y}$. 
    $\varepsilon_{\mathbf{y}} \ne 0$ because in most cases, 
    a PDE's ground-truth solution could only be 
    a numerical approximation with errors. 
    The upper bound is input-dependent because 
    $\widetilde{\mathbf{A}}$ and $\widetilde{\mathbf{f}}$
    are input-dependent. 

    Similarly, we could prove that the
    lower bound of $\varepsilon_{\mathbf{x}}$ is also input-dependent:
    \begin{align*}
        \varepsilon_{\mathbf{x}} 
        & = \dfrac{\norm{\widetilde{\mathbf{f}} - 
                         \widetilde{\mathbf{A}} \, \mathbf{x}}}
                  {\norm{\widetilde{\mathbf{f}}}} 
        \approx \dfrac{\norm{\widetilde{\mathbf{f}} - 
                               \widetilde{\mathbf{A}} \, 
                               \paren{\mathbf{y} \pm 
                                      l_{\mathrm{max}} \,
                                      \mathbf{y}}}}
                        {\norm{\widetilde{\mathbf{f}}}} 
        = \dfrac{\norm{\paren{\widetilde{\mathbf{f}} - 
                          \widetilde{\mathbf{A}} \, 
                          \mathbf{y}}
                         \mp
                         \paren{l_{\mathrm{max}} \, 
                                \widetilde{\mathbf{A}} \, 
                                \mathbf{y}}}}
                  {\norm{\widetilde{\mathbf{f}}}} \notag \\ 
        & \ge \abs{
              \dfrac{\norm{\widetilde{\mathbf{f}} - 
               \widetilde{\mathbf{A}} \, \mathbf{y}}}
              {\norm{\widetilde{\mathbf{f}}}}
              - 
              \dfrac{\norm{l_{\mathrm{max}} \, 
                         \widetilde{\mathbf{A}} \, 
                         \mathbf{y}}}
                    {\norm{\widetilde{\mathbf{f}}}} 
              }
        = \abs{
          \varepsilon_{\mathbf{y}}
              - 
              \dfrac{\norm{l_{\mathrm{max}} \, 
                         \widetilde{\mathbf{A}} \, 
                         \mathbf{y}}}
                    {\norm{\widetilde{\mathbf{f}}}}
          } 
        > 0 \ . 
    \end{align*}
\end{proof}

\subsection{Ablation Study}
\label{sec:appendix:ablation_study}

For experimental completeness, 
we conduct ablation studies with respect to
our residual loss and the UGrid architecture itself. 
In addition to the UGrid model trained with residual loss
(as proposed in the main content), 
we also train another UGrid model with legacy loss, 
as well as one vanilla U-Net model
with residual loss. 
This U-Net model has the same number of layers as UGrid, 
and has non-linear layers as proposed in 
\cite{Ronneberger15}.  
We let the U-Net directly regress 
the solutions to Poisson's equations. 
All these models are trained with the same data
in the same manner (except for the loss metric), 
as detailed in Section~\ref{sec:experiments}. 


We conduct qualitative experiments on the same set of
testcases as detailed in Section~\ref{sec:experiments}, 
and the results are as follows: 

\begin{table}[!htb]
    \centering
    \caption{Ablation study on large-scale Poisson problems. 
    }
    \label{tab:comparasion_1025_ablation}
    \begin{tabular}{cccc}
        \toprule
        \textbf{Testcase} & 
        \textbf{UGrid} & 
        \textbf{UGrid (L)} & 
        \textbf{U-Net} \\
        Poisson (L) & 
        Time / Error & 
        Time / Error & 
        Time / Error \\
        \midrule
        
        Bag & 
        \textbf{18.66} / 2.66 & 
        28.81 / 4.86 & 
        81.71 / 1384131 \\ 
        
        Cat & 
        \textbf{10.09} / 2.70 & 
        23.80 / 1.43 & 
        70.09 / 2539002 \\ 
        
        Lock & 
        \textbf{10.55} / 9.88 & 
        Diverge & 
        70.92 / 1040837 \\ 

        N. Input & 
        \textbf{10.16} / 2.64 & 
        20.65 / 2.42 & 
        73.05 / 21677 \\ 
        
        Note & 
        \textbf{10.31} / 4.06 & 
        Diverge & 
        69.97 / 614779 \\ 
        
        S. Feat. & 
        \textbf{20.01} / 3.80 & 
        31.34 / 5.14 & 
        70.08 / 222020 \\ 
        
        L-shape & 
        \textbf{15.26} / 8.43 & 
        Diverge & 
        74.67 / 1800815  \\ 
        
        Lap. Squ. & 
        \textbf{15.10} / 3.88 & 
        30.72 / 2.76 & 
        72.24 / 30793035 \\ 
        
        P. Squ. & 
        \textbf{15.07} / 9.37 & 
        31.52 / 3.33 & 
        71.74 / 31043896 \\ 
       
        Star & 
        \textbf{15.18} / 7.50 & 
        Diverge & 
        70.01 / 1138821 \\ 
        
        \bottomrule
    \end{tabular}
\end{table}

In Table~\ref{tab:comparasion_1025_ablation}, 
the residual loss endows our UGrid model with 
as much as $2$x speed up versus the legacy loss. 
The residual loss also endows UGrid to
converge to the failure cases of 
its counterpart trained on legacy loss. 
These results demonstrate the claimed merits
of the residual loss. 
On the other hand, it will \textbf{diverge} if we naively apply 
the vanilla U-Net architecture directly to Poisson's
equations. For experimental completeness only, we list 
the diverged results in the last column. 
(The ``time'' column measures the time taken for $64$ iterations; 
the iterators are shut down once they reach this threshold.)
This showcases the significance 
of UGrid's mathematically-rigorous network architecture. 


\begin{table}[!htb]
    \centering
    \caption{Ablation study on small-scale Poisson problems. 
    }
    \label{tab:comparasion_257_ablation}
    \begin{tabular}{cccc}
        \toprule
        \textbf{Testcase} & 
        \textbf{UGrid} & 
        \textbf{UGrid (L)} & 
        \textbf{U-Net} \\
        Poisson (S) & 
        Time / Error & 
        Time / Error &
        Time / Error \\
        \midrule
        
        Bag & 
        \textbf{8.76} / 8.05 & 
        17.89 / 4.50 & 
        71.86 / 678141 \\ 
        
        Cat & 
        \textbf{51.96} / 6.21 & 
        Diverge & 
        68.89 / 1317465 \\ 
        
        Lock & 
        \textbf{9.00} / 2.11 & 
        18.32 / 2.83 & 
        69.47 / 189412 \\ 

        Noisy Input & 
        \textbf{8.94} / 6.00 & 
        17.88 / 6.58 & 
        69.54 / 21666 \\
        
        Note & 
        \textbf{8.87} / 2.75 & 
        17.79 / 3.06 & 
        69.59 / 24715 \\ 
        
        Sharp Feature & 
        \textbf{13.31} / 7.52 & 
        26.64 / 1.91 &
        70.57 / 191499 \\ 
        
        L-shape & 
        \textbf{40.60} / 7.09 & 
        Diverge & 
        69.71 / 1011364 \\
        
        Laplacian Square & 
        \textbf{13.21} / 3.27 & 
        22.23 / 9.55 & 
        73.80 / 15793109 \\ 
        
        Poisson Square & 
        \textbf{13.21} / 2.88 & 
        22.13 / 9.76 & 
        71.56 / 15393069 \\
       
        Star & 
        \textbf{8.92} / 2.36 & 
        17.60 / 5.69 & 
        73.72 / 502993 \\
        
        \bottomrule
    \end{tabular}
\end{table}

In Table~\ref{tab:comparasion_257_ablation}, 
for small-scale problems, 
the residual loss still endows UGrid with
as much as $2$x speedup and 
stronger generalization power
against its counterpart trained with legacy loss. 
Once again, the vanilla U-Net model diverged
for all testcases, and we list its diverged
results for experimental completeness only.

\subsection{Qualitative Evaluations on 
Small-scale Poisson Problems}
\label{sec:supplemental_material:eval_small_poisson}

\begin{table*}[!htb]
    \centering
    \caption{Comparison of UGrid 
    and state-of-the-art on small-scale Poisson problems.}
    \label{tab:comparasion_257}
    \begin{tabular}{ccccc}
        \toprule
        \textbf{Testcase} & 
        \textbf{UGrid} & 
        \textbf{AMGCL} & 
        \textbf{AmgX} &
        \textbf{Hsieh et al.} \\
        Poisson (S) & 
        Time / Error & 
        Time / Error &
        Time / Error & 
        Time / Error \\
        \midrule
        
        Bag & 
        \textbf{8.76} / 8.05 & 
        12.14 / 3.00 & 
        22.34 / 8.20 &
        47.69 / 252 \\ 
        
        Cat & 
        51.96 / 6.21 & 
        \textbf{17.03} / 6.98 & 
        27.66 / 4.83 &
        23.02 / 9.95 \\ 
        
        Lock & 
        \textbf{9.00} / 2.11 & 
        15.77 / 7.89 & 
        16.96 / 9.36 & 
        48.72 / 117.9 \\ 

        Noisy Input & 
        \textbf{8.94} / 6.00 & 
        14.00 / 9.39 & 
        26.30 / 3.14 &
        51.79 / 5576 \\
        
        Note & 
        8.87 / 2.75 & 
        \textbf{8.79} / 9.02 & 
        16.68 / 7.23 &
        36.66 / 8.28 \\ 
        
        Sharp Feature & 
        \textbf{13.31} / 7.52 & 
        21.47 / 4.15 & 
        49.59 / 6.85 &
        49.31 / 24876 \\ 
        
        L-shape & 
        40.60 / 7.09 & 
        \textbf{12.36} / 9.97 & 
        24.08 / 9.35 &
        50.06 / 96.44 \\
        
        Laplacian Square & 
        \textbf{13.21} / 3.27 & 
        22.22 / 5.60 & 
        48.60 / 3.98 &
        24.57 / 6.54 \\ 
        
        Poisson Square & 
        \textbf{13.21} / 2.88 & 
        21.93 / 5.51 & 
        47.56 / 4.03 &
        49.77 / 473 \\
       
        Star & 
        \textbf{8.92} / 2.36 & 
        18.93 / 2.17 & 
        17.96 / 9.42 &
        48.68 / 456 \\
        
        \bottomrule
    \end{tabular}
\end{table*}

In Table~\ref{tab:comparasion_257}, 
even on small-scale problems that hinder our solver 
with a compact multigrid-like hierarchy 
from delivering its full power, 
the UGrid model is still faster than or 
exhibits comparable efficiency
with respect to the three SOTA legacy/neural solvers.  
Again, this shows the high efficiency 
as well as the strong generalization power of our new method. 
The testcases ``Cat'' and ``L-shape'' showcase that the
generalization power (in terms of problem size)
does come with a price of potentially downgraded efficiency.  
Thus, for the sake of the best efficiency, 
we still recommend re-training UGrid 
for problems of different sizes. 

\subsection{Specifications on Evaluations of
Inhomogeneous Helmholtz Problems 
with Spatially-varying Wavenumbers}
\label{sec:supplemental_material:eval_helmholtz}

\label{sec:appendix:helmholtz}
We train UGrid with the same training data 
and residual loss as mentioned
in Section~\ref{sec:experiments}. 
As one exception, 
we also input randomly-sampled $k^2$ 
during training, evaluation, and testing. 

\begin{figure*}[!htb]
    \centering
    \begin{subfigure}[H]{0.19\linewidth}
        \centering
        \includegraphics[height=\linewidth]{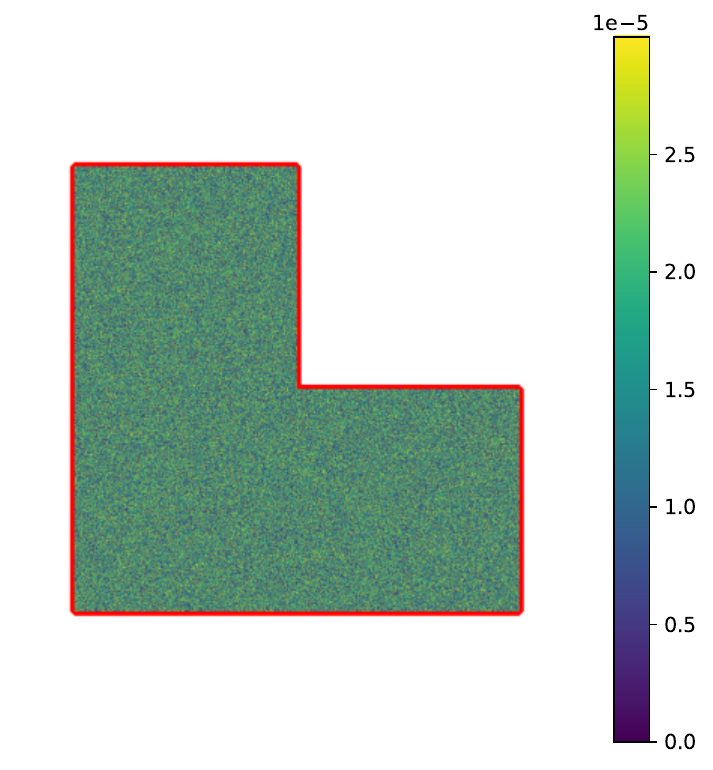}
        \caption*{(a) L-Shape}
    \end{subfigure}
    \begin{subfigure}[H]{0.19\linewidth}
        \centering
        \includegraphics[height=\linewidth]{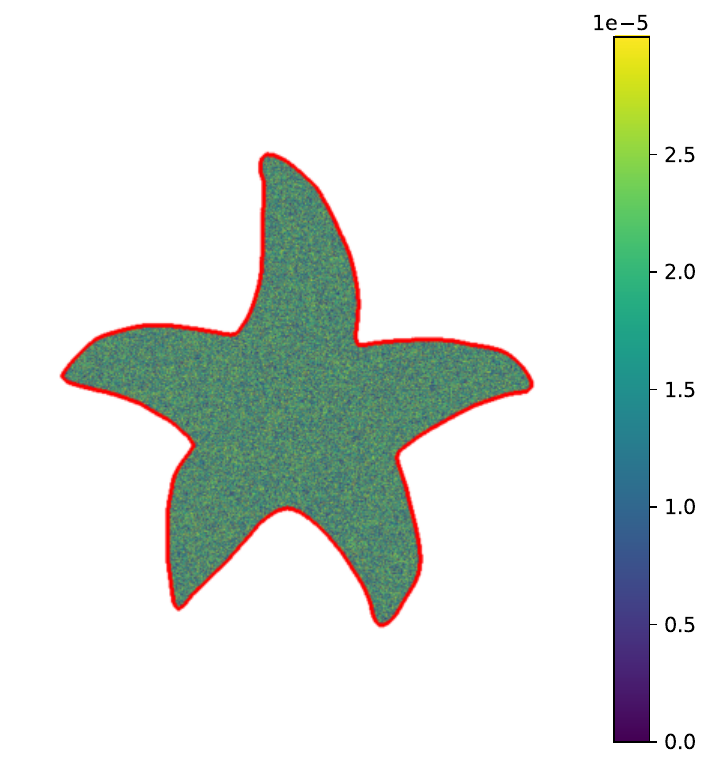}
        \caption*{(b) Star}
    \end{subfigure}
    \begin{subfigure}[H]{0.19\linewidth}
        \centering
        \includegraphics[height=\linewidth]{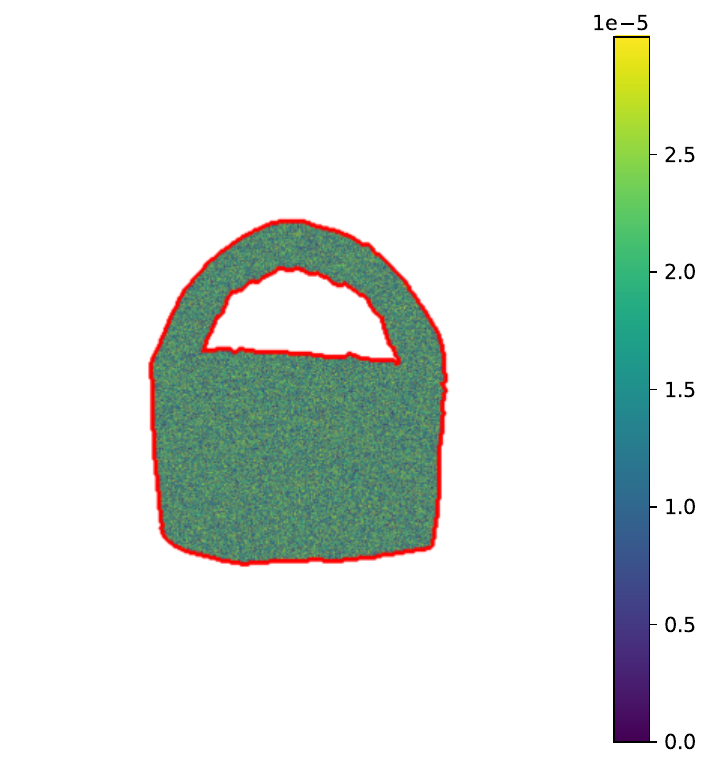}
        \caption*{(c) Lock}
    \end{subfigure}
    \begin{subfigure}[H]{0.19\linewidth}
        \centering
        \includegraphics[height=\linewidth]{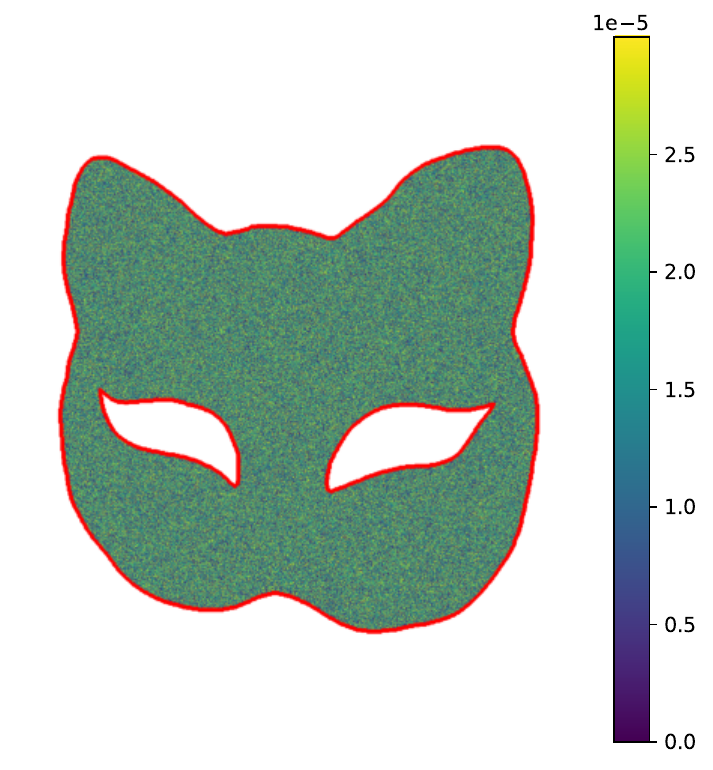}
        \caption*{(d) Cat}
    \end{subfigure}
    \begin{subfigure}[H]{0.19\linewidth}
        \centering
        \includegraphics[height=\linewidth]{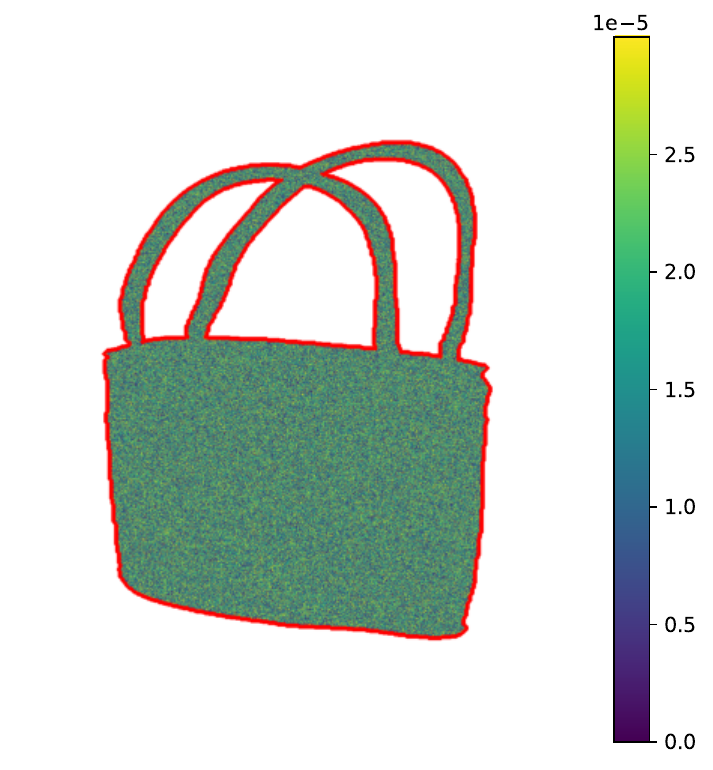}
        \caption*{(e) Bag}
    \end{subfigure}
    \begin{subfigure}[H]{0.19\linewidth}
        \centering
        \includegraphics[height=\linewidth]{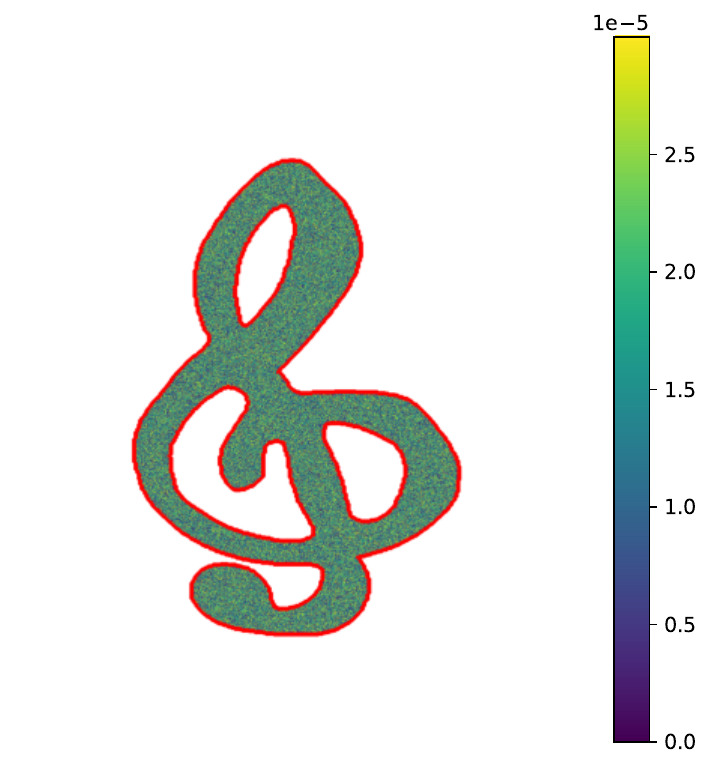}
        \caption*{(f) Note}
    \end{subfigure}
    \begin{subfigure}[H]{0.19\linewidth}
        \centering
        \includegraphics[height=\linewidth]{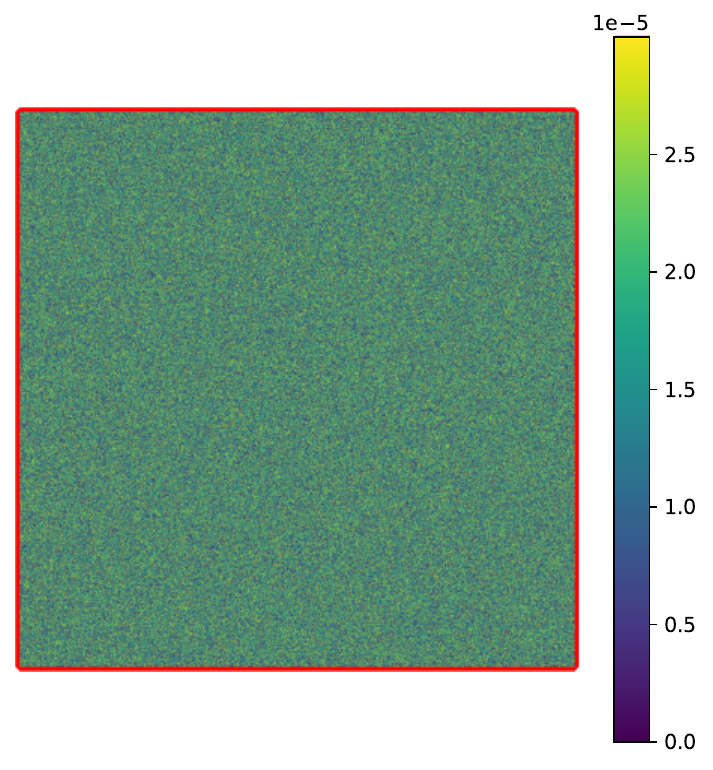}
        \caption*{(g) Sharp Feature}
    \end{subfigure}
    \begin{subfigure}[H]{0.19\linewidth}
        \centering
        \includegraphics[height=\linewidth]{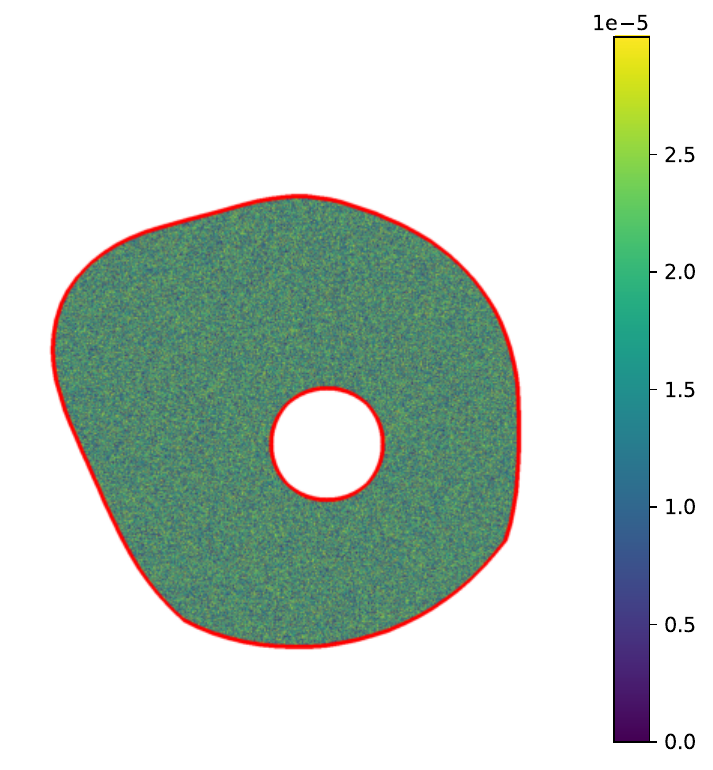}
        \caption*{(h) Noisy Input}
    \end{subfigure}
    \begin{subfigure}[H]{0.19\linewidth}
        \centering
        \includegraphics[height=\linewidth]{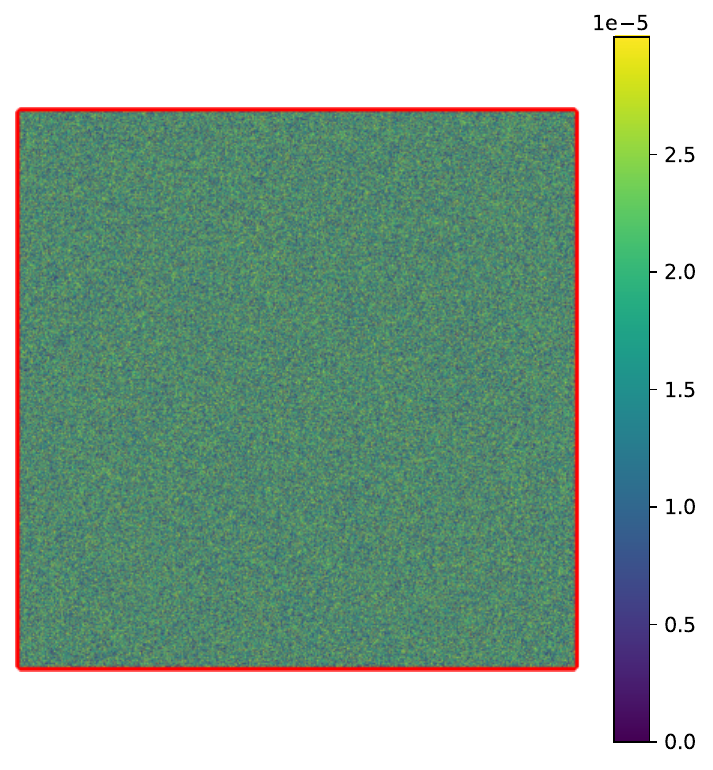}
        \caption*{(i) Lap. Square}
    \end{subfigure}
    \begin{subfigure}[H]{0.19\linewidth}
        \centering
        \includegraphics[height=\linewidth]{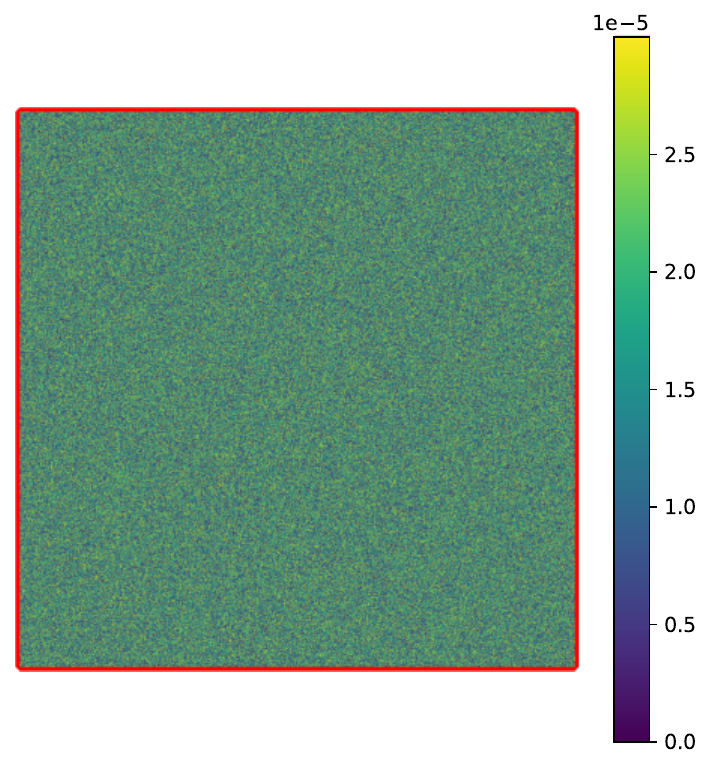}
        \caption*{(j) Poisson Square}
    \end{subfigure}
    \caption{
    Illustration of the wavenumber distributions
    of our testcases. 
    The boundaries are shown in bold 
    red lines for a better view. 
    (Boundary values are shown in 
    Fig.~\ref{fig:testcase_boundaries}.)
    }
    \label{fig:testcase_k2}
\end{figure*}

The randomly-sampled $k^2$s we used
are illustrated in Fig.~\ref{fig:testcase_k2}. 
For qualitative experiments, 
we use the same boundary conditions and Laplacian distributions
as shown in Fig.~\ref{fig:testcase_boundaries}
and Fig.~\ref{fig:testcase}, and we randomly initialize the
wavenumber field $k^2$ 
across the whole domain, resulting in a noisy distribution. 

\begin{table}[!htb]
    \centering
    \caption{Qualitative results for small-scale Helmholtz problems. 
    }
    \label{tab:comparasion_257_helmholtz}
    \begin{tabular}{cccc}
        \toprule
        \textbf{Testcase} & 
        \textbf{UGrid} & 
        \textbf{AMGCL} & 
        \textbf{AmgX} \\
        Helmholtz (S) & 
        Time / Error & 
        Time / Error & 
        Time / Error \\
        \midrule
        
        Bag & 
        14.70 / 6.29 & 
        \textbf{12.92} / 3.00 & 
        25.82 / 1.79 \\ 
        
        Cat & 
        \textbf{16.63} / 7.86 & 
        16.82 / 6.98 & 
        28.37 / 3.03 \\ 
        
        Lock & 
        \textbf{9.78} / 5.87 & 
        16.23 / 7.88 & 
        19.75 / 2.02 \\ 

        Noisy Input & 
        14.95 / 0.76 & 
        \textbf{14.34} / 9.40 & 
        28.92 / 0.04 \\
        
        Note & 
        14.37 / 8.28 & 
        \textbf{9.01} / 9.02 & 
        18.76 / 2.55 \\ 
        
        Sharp Feature & 
        \textbf{19.46} / 1.18 & 
        21.37 / 4.21 & 
        52.82 / 0.13 \\ 
        
        L-shape & 
        14.64 / 0.88 & 
        \textbf{12.29} / 9.99 & 
        26.90 / 2.10 \\
        
        Laplacian Square & 
        \textbf{14.60} / 4.60 & 
        22.43 / 5.59 & 
        43.68 / 17.8 \\ 
        
        Poisson Square & 
        \textbf{15.27} / 6.53 & 
        22.35 / 5.50 & 
        43.57 / 17.2 \\
       
        Star & 
        \textbf{9.77} / 5.96 & 
        19.09 / 2.16 & 
        20.89 / 2.18 \\
        
        \bottomrule
    \end{tabular}
\end{table}

The qualitative results for large-scale problems 
are available in the main paper as Table~\ref{tab:comparasion_1025_helmholtz}. 
Qualitative results for small-scale problems are available in 
Table~\ref{tab:comparasion_257_helmholtz}. 
We could observe that even on small-scale problems that hinder
our solver with a compact multigrid-like hierarchy 
from delivering its full power, 
UGrid is still faster than or 
exhibits comparable efficiency
with respect to the SOTA.

\newpage
\subsection{Specifications on Evaluations of
Inhomogeneous Steady-state Convection-diffusion-reaction Problems}
\label{sec:supplemental_material:eval_diffusion}

We train UGrid in the same manner as for
Helmholtz equations. 
As one exception, 
we input randomly-sampled $\vb{v}$s, 
$\alpha$s, and $\beta$s
during training, evaluation, and testing. 
These values are sampled using the same routine 
as for Helmholtz equations, 
resulting in noisy velocity fields 
like Fig.~\ref{fig:testcase_k2}
as well as randomized $\alpha$, $\beta$ coefficients
($4 \alpha + \beta \ne 0$). 
The qualitative results 
for large-scale problems are available
in the main paper as Table~\ref{tab:comparasion_1025_diffusion_convection}. 
Qualitative results for the small-scale problem are available 
in Table~\ref{tab:comparasion_257_diffusion_convection}. 

\begin{table}[!htb]
    \centering
    \caption{Qualitative results for small-scale Diffusion problems. 
    }
    \label{tab:comparasion_257_diffusion_convection}
    \begin{tabular}{cccc}
        \toprule
        \textbf{Testcase} & 
        \textbf{UGrid} & 
        \textbf{AMGCL} & 
        \textbf{AmgX} \\
        Diffusion (S) & 
        Time / Error & 
        Time / Error & 
        Time / Error \\
        \midrule
        
        Bag & 
        16.99 / 3.79 & 
        \textbf{12.34} / 9.94 & 
        22.58 / 3.66 \\ 
        
        Cat & 
        69.12 / 8.76 & 
        \textbf{16.57} / 9.79 & 
        27.44 / 6.74 \\ 
        
        Lock & 
        17.07 / 1.43 & 
        \textbf{16.07} / 6.34 & 
        18.34 / 4.06 \\ 

        Noisy Input & 
        22.84 / 6.47 & 
        \textbf{15.94} / 2.74 & 
        24.35 / 0.42 \\ 
        
        Note & 
        22.76 / 1.17 & 
        \textbf{9.38} / 2.67 & 
        19.25 / 3.79 \\ 
        
        Sharp Feature & 
        \textbf{17.05} / 5.15 & 
        21.41 / 4.29 & 
        41.67 / 0.72 \\ 
        
        L-shape & 
        35.18 / 6.25 & 
        \textbf{13.02} / 2.30 & 
        25.51 / 3.97 \\ 
        
        Laplacian Square & 
        90.73 / 64.5 & 
        \textbf{22.14} / 8.17 & 
        50.49 / 3.50 \\ 
        
        Poisson Square & 
        50.95 / 5.01 & 
        \textbf{22.05} / 7.34 & 
        50.06 / 3.28 \\ 
       
        Star & 
        \textbf{17.06} / 3.69 & 
        18.70 / 7.55 & 
        18.88 / 4.71 \\ 
        
        \bottomrule
    \end{tabular}
\end{table}

\newpage
In Table~\ref{tab:comparasion_257_diffusion_convection}, 
again, even on small-scale problems that hinder
our solver with a compact multigrid-like hierarchy 
from delivering its full power, 
UGrid is still 
exhibits comparable efficiency
with respect to the SOTA.  
This demonstrates UGrid's generalization power over
problem sizes, though possibly at the price of 
relatively lower efficiency compared to the size
it is trained on. 

\subsection{Additional Experiments on Poisson Problems}

We have conducted four additional experiments to showcase the generalization power of the proposed UGrid solver. The four testcases are illustrated as follows: 

\begin{figure*}[ht]
    \centering
    \begin{subfigure}[H]{0.49\linewidth}
        \centering
        \includegraphics[height=0.49\linewidth]{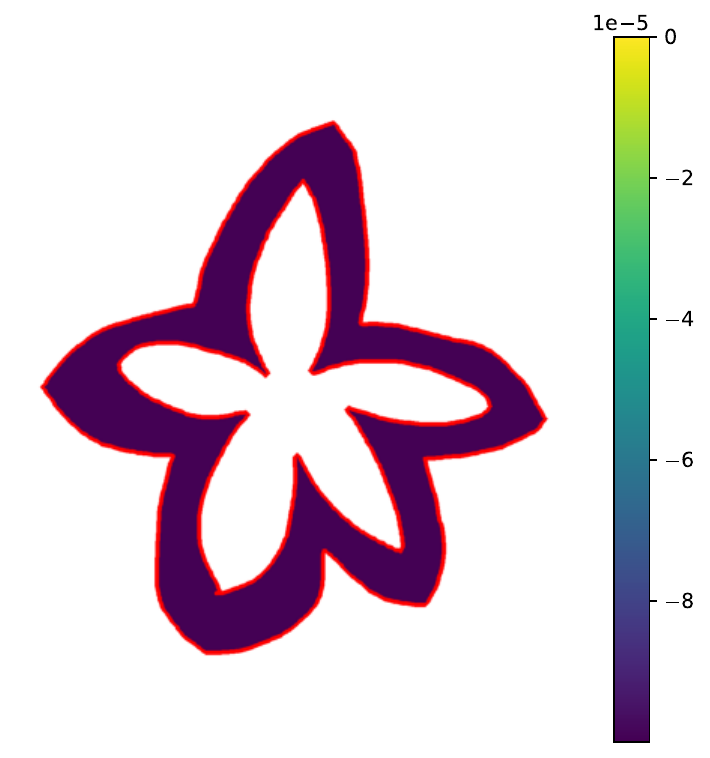}
        \includegraphics[height=0.49\linewidth]{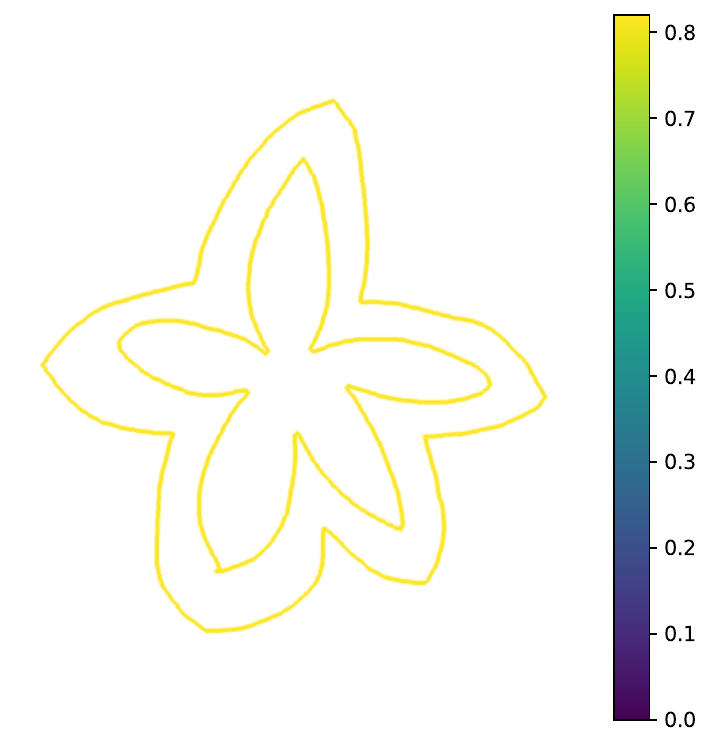}
        \caption{Flower (Left: Laplacian; Right: Boundary)}
    \end{subfigure}
    \begin{subfigure}[H]{0.49\linewidth}
        \centering
        \includegraphics[height=0.49\linewidth]{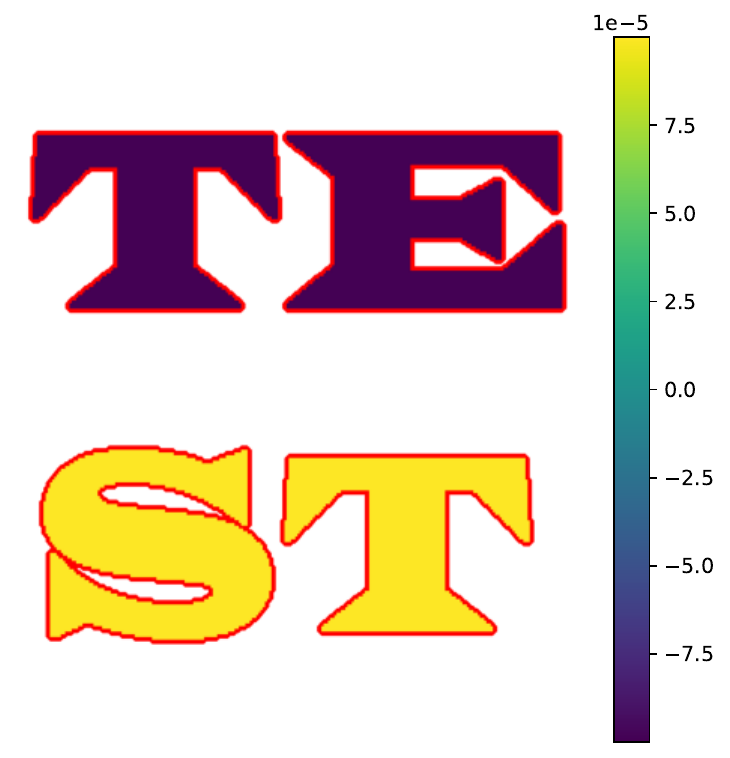}
        \includegraphics[height=0.49\linewidth]{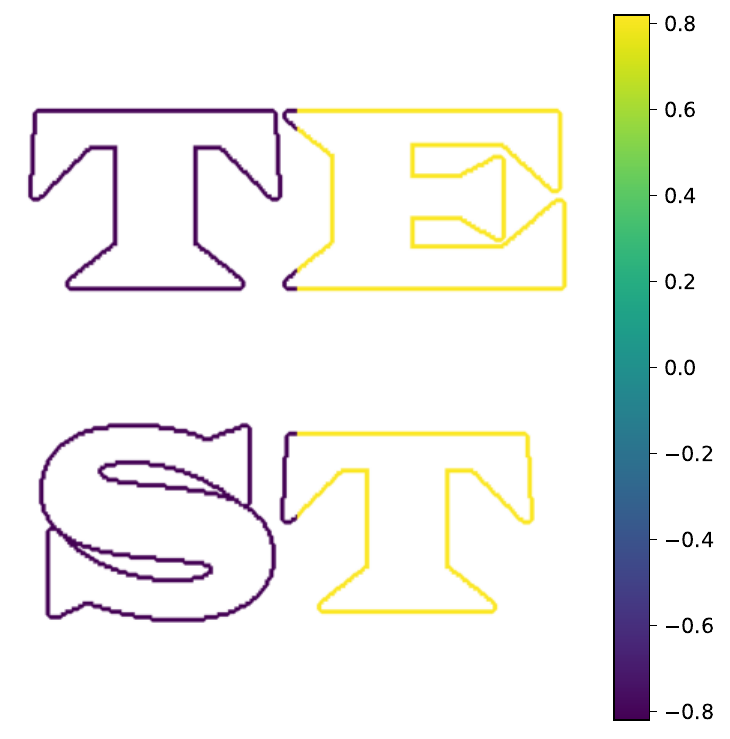}
        \caption{Word (Left: Laplacian; Right: Boundary)}
    \end{subfigure}
    \begin{subfigure}[H]{0.49\linewidth}
        \centering
        \includegraphics[height=0.49\linewidth]{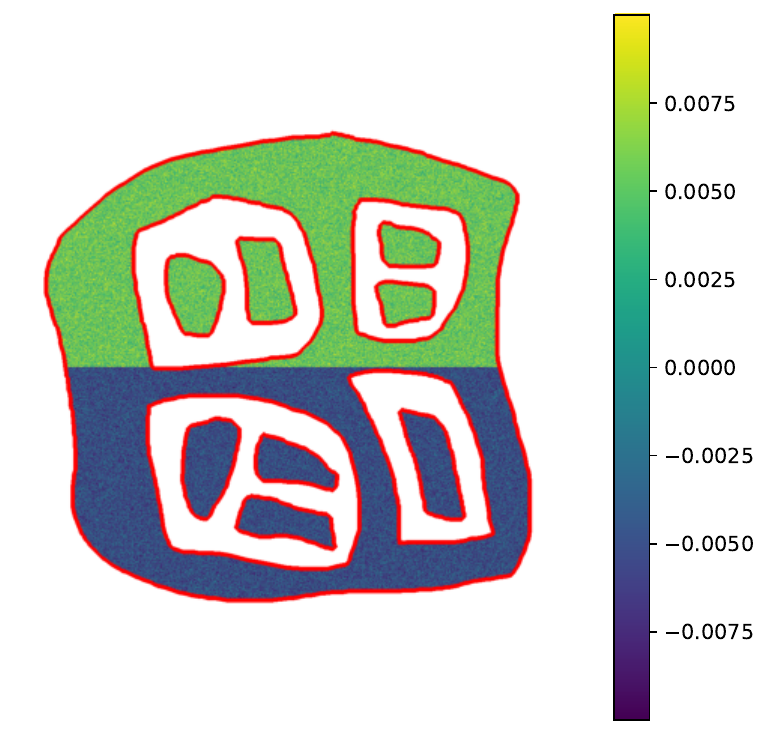}
        \includegraphics[height=0.49\linewidth]{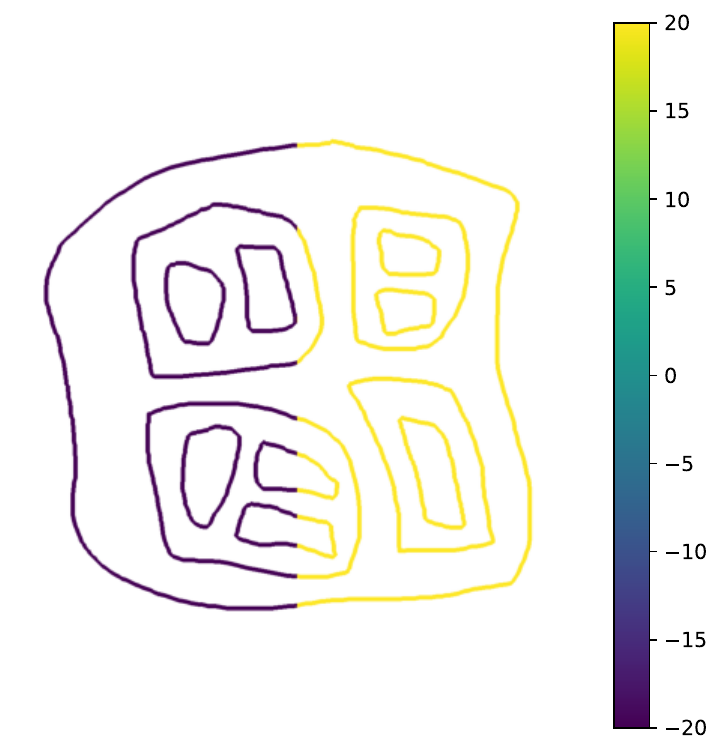}
        \caption{Topo (Left: Laplacian; Right: Boundary)}
    \end{subfigure}
    \begin{subfigure}[H]{0.49\linewidth}
        \centering
        \includegraphics[height=0.49\linewidth]{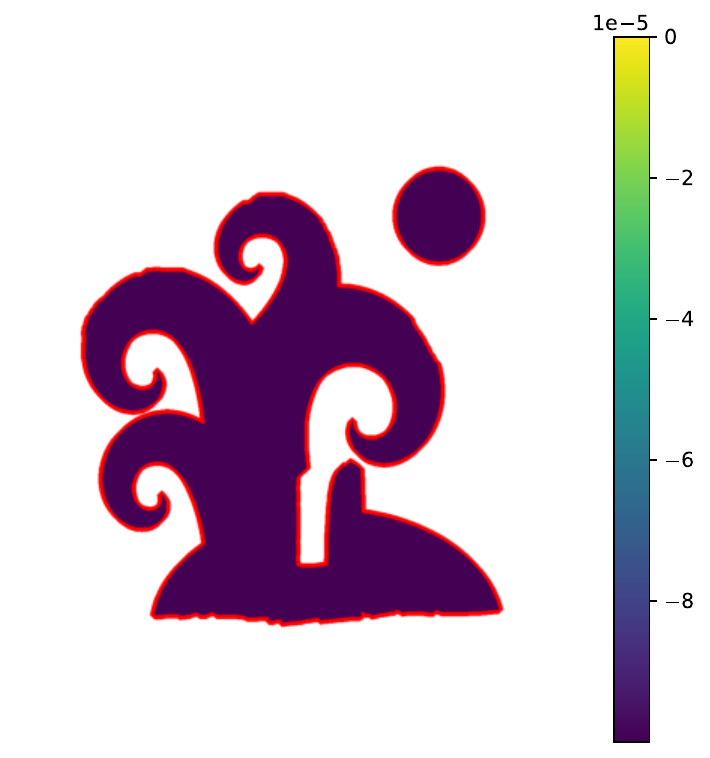}
        \includegraphics[height=0.49\linewidth]{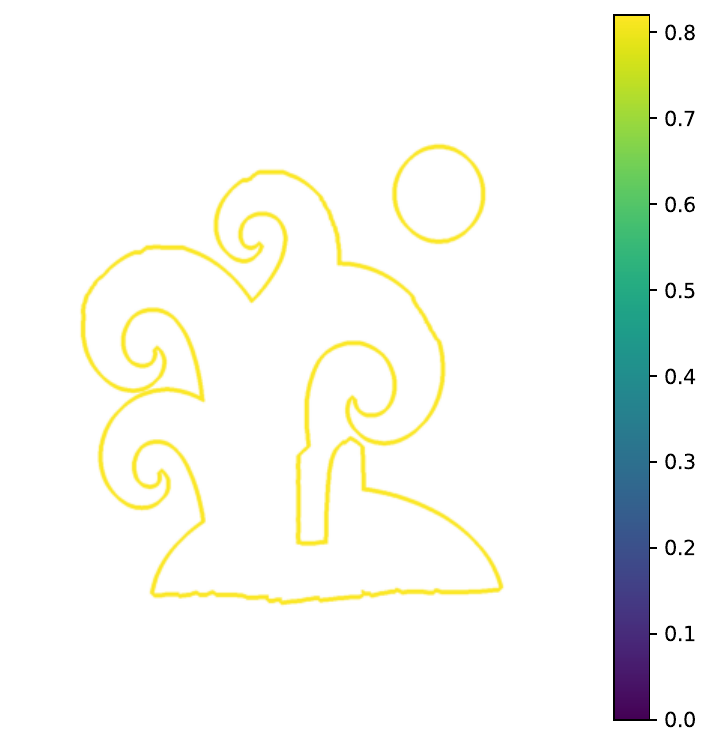}
        \caption{Tree (Left: Laplacian; Right: Boundary)}
    \end{subfigure}
    \caption{Illustration of four additional testcases.
    }
    \label{fig:testcase_extra}
\end{figure*}

\newpage
The qualitative results are as follows: 

\begin{table*}[ht]
    \centering
    \caption{Qualitative results for 
    additional large-scale Poisson problems. 
    ``Time'' denotes the total time ($\mathrm{ms}$) 
    to reach relative residual errors $\le 10^{-4}$, 
    and ``Assembly'' and ``Iteration'' denotes time ($\mathrm{ms}$) for the two phases for AMG solvers;
    ``Error'' denotes 
    the final relative residual errors, 
    divided by $10^{-5}$. 
    }
    \label{tab:comparasion_1025_extra}
    \begin{tabular}{ccccc}
        \toprule
        \textbf{Testcase} & 
        \textbf{UGrid} & 
        \textbf{AMGCL} & 
        \textbf{AmgX} & 
        \textbf{Hsieh et al.} \\
        Poisson (L) & 
        Time / Error & 
        Time / Assembly / Iteration / Error & 
        Time / Assembly / Iteration / Error & 
        Time / Error \\
        \midrule
        
        Flower & 
        \textbf{18.52} / 1.60 & 
        144.91 / 138.00 / 6.91 / 3.07 & 
        69.14 / 18.11 / 51.03 / 5.41 &
        27.23 / 7.16 \\ 
        
        Word & 
        \textbf{41.14} / 5.60 & 
        188.59 / 178.34 / 10.25 / 4.15 & 
        108.51 / 21.44 / 87.06 / 4.62 & 
        50.76 / 2612 \\ 
        
        Topo & 
        \textbf{10.30} / 4.29 & 
        210.06 / 203.86 / 6.20 / 4.86 & 
        97.22 / 24.29 / 72.93 / 2.94 & 
        52.29 / 3933 \\ 

        Tree & 
        \textbf{10.05} / 2.52 & 
        173.70 / 166.86 / 6.84 / 4.61 & 
        78.73 / 20.98 / 57.76 / 6.88 &
        16.05 / 7.60 \\
        
        \bottomrule
    \end{tabular}
\end{table*}

In Table~\ref{tab:comparasion_1025_extra}, UGrid still delivers similar efficiency, accuracy and generalization power like the ten testcases covered in the main contents of the paper. Note that Hsieh et al. \textbf{diverged} for testcases "Word" and "Topo"; their time for these two cases is the time to reach the maximum number of iterations.

\end{document}